\documentclass[12pt,a4paper]{article}
\usepackage[utf8x]{inputenc}
\usepackage[english,russian]{babel}
\usepackage{amsmath, amsfonts}
\newcommand{\mx}{\operatorname{\mathfrak M}}
\newcommand{\mn}{\operatorname{\mathfrak m}}
\renewcommand{\Im}{\operatorname{Im}}
\newcommand{\Ker}{\operatorname{Ker}}
\newcommand{\loc}{\operatorname{loc}}
\newcommand{\e}{\operatorname{\mathfrak e}}
\renewcommand{\span}{\operatorname{span}}
\newcommand{\close}{\operatorname{close}}
\newcommand{\s}{\operatorname{\mathcal s}}
\newcommand{\card}{\operatorname{card}}
\newcommand{\diam}{\operatorname{diam}}
\newcommand{\supp}{\operatorname{supp}}
\newcommand{\supvrai}{\operatornamewithlimits{sup\,vrai}}
\newcommand{\N}{\mathbb N}
\newcommand{\Z}{\mathbb Z}
\newcommand{\R}{\mathbb R}
\renewcommand{\C}{\mathbb C}
\newcommand{\Nu}{\mathcal N}

\newcommand{\mes}{\operatorname{mes}}

\begin{document}

\author{ С. Н. Кудрявцев }
\title{Аналог теоремы Литтлвуда-Пэли для ортопроекторов
на подпространства всплесков}
\date{}
\maketitle
УДК 517.5
\begin{abstract}
В статье доказано утверждение, представляющее собой
аналог теоремы Литтлвуда-Пэли для ортопроекторов на
подпространства всплесков, соответствующие неизотропному
кратно-масштабному анализу, порожденному тензорным произведением
гладких достаточно быстро стремящихся к нулю на бесконечности  масштабирующих функций
одной переменной.
\end{abstract}

Ключевые слова: ортопроектор, подпространства всплесков,
масштабирующая функция, кратно-масштабный анализ,
теорема Литтлвуда-Пэли
\bigskip

\centerline{Введение}

Как известно (см., например, [1], [2] и др. работы), важное значение для
вывода порядковых оценок
точности приближения в $L_p$ способами,
основанными на применении кратных тригонометрических рядов, классов
периодических функций нескольких переменных с условиями на смешанные
производные (разности) имеют теорема Литтлвуда-Пэли для кратных рядов
Фурье и следствия из нее.
По поводу теоремы Литтлвуда-Пэли для кратных рядов Фурье
см., например, [3, п. 1.5.2], [4] и приведенную там литературу.
Для средств приближения классов (непериодических) функций, Заданных
на кубе $ I^d, $ подчиненных условиям на смешанные разности, аналог теоремы
Литтлвуда-Пэли установлен автором в [5].
Для получения соответствующих оценок точности приближения классов
непериодических функций смешанной гладкости, заданных на всем
пространстве $ \R^d, $ полезно иметь аналоги этих утверждений для средств
приближения таких классов функций.
С этой целью в работе доказан аналог теоремы Литтлвуда-Пэли для
семейства ортопроекторов $ \{\mathcal E_\kappa, \kappa \in \Z_+^d \}, $
(см. п. 2.4.) на подпространства всплесков, соответствующие неизотропному
кратно-масштабному анализу (КМА), порожденному тензорным произведением
гладких достаточно быстро стремящихся к нулю на бесконечности масштабирующих
функций одной переменной, а именно, показано, что при $ 1 < p < \infty $
для любой функции $ f \in L_p(\R^d) $ выполняются неравенства
\begin{equation*} \tag{1}
c_{1} \| f \|_{L_p(\R^d)} \le
\biggl( \int_{\R^d} (\sum_{\kappa \in \Z_+^d} |(\mathcal E_\kappa f)(x)|^2 )^{p/2}
dx\biggr)^{1/p} \le c_{2} \| f \|_{L_p(\R^d)}
\end{equation*}
с некоторыми константами $ c_1, c_2 > 0, $ зависящими от $ d,p $ и функций,
порождающих КМА.
Задача построения масштабирующих функций и всплесков не является предметом
рассмотрения настоящей работы, и по этому поводу см., например, [6], [7].
Известные автору варианты аналогов теоремы Литтлвуда-Пэли, относящиеся к
теории всплесков, не рассматривают неизотропный случай КМА (см. [8, гл. 6, \S 2], [9, п. 6.4]). Кроме того,
в [8] и [9] соответствующие утверждения описывают связь между свойствами
масштабирующей и/или всплеск-функции с одной стороны и свойствами системы
коэффициентов разложения функций в соответствующем базисе всплесков
с другой стороны. В то время как в настоящей работе устанавливается связь
между свойствами масштабирующей и двойственной ей функций
и свойствами семейства ортопроекторов на соответствующие подпространства
всплесков. Причем, условия при которых доказаны отмеченные выше утверждения в
[8], [9] не совпадают с условиями, при которых установлено соотношение (1).
Отметим еще, что схема и средства проведения доказательства (1) отличаются
от тех, что использовались в [8], [9] для получения имеющихся там вариантов
аналогов теоремы Литтлвуда-Пэли.

Работа состоит из введения и двух параграфов. В \S 1 приведены
предварительные сведения, используемые при доказательстве основных результатов
работы. В п. 2.4. \S 2 на основании фактов, полученных в предыдущих пунктах,
устанавливается справедливость (1). Перейдем к точным формулировкам
и доказательствам.
\bigskip

\centerline{ \S 1. Предварительные сведения и вспомогательные утверждения}

 1.1. В этом пункте вводятся обозначения, используемые
в настоящей работе, а также приводятся некоторые факты, необходимые
в дальнейшем.

Для $ d \in \N $ через $ \Z_+^d $ обозначим множество
$$
\Z_+^d =\{\lambda =(\lambda_1,  \ldots, \lambda_d) \in \Z^d:
\lambda_j \ge0, j=1, \ldots, d\}.
$$
 Обозначим также при  $ d \in \N $ для $ l \in \Z_+^d $ через
$ \Z_+^d(l) $ множество
 $$
 \Z_+^d(l) =\{ \lambda  \in \Z_+^d: \lambda_j \le l_j, j=1, \ldots,  d\}.
 $$

Для комплексного числа $ z $ через $ \overline z $ обозначается
комплексно-сопряженное число.

Через $ C^1(D) $ обозначается пространство комплекснозначных непрерывно
дифференцируемых функций в области $ D \subset \R^d, $
а через $ C_0^\infty(D) $ -- пространство бесконечно дифференцируемых
финитных функций, носители которых лежат в $ D. $

Для измеримого по Лебегу множества $ D \subset \R^d$
при $ 1 \le p \le \infty $ через  $  L_p(D),$  как обычно,
обозначается пространство  всех   комплекснозначных измеримых на $ D $
функций $f,$ для которых определена норма
\begin{equation*}
\|f\|_{L_p(D)}= \begin{cases} (\int_D |f(x)|^p
dx)^{1/p}, 1 \le p<\infty;\\
\supvrai_{x \in D}|f(x)|, p=\infty. \end{cases}
\end{equation*}
$ L_2(D) $ является гильбертовым пространством со скалярным произведением
$ \langle f, g \rangle = \int_D f \overline g dx, f,g \in L_2(D). $

Обозначим еще через $ L_1^{\loc}(\R^d) $ линейное пространство всех локально
суммируемых в $ \R^d $ функций, принимающих комплексные значения.

Введем еще следующие обозначения.

Для $ x,y \in \R^d $ положим $ xy =x \cdot y = (x_1 y_1, \ldots,
x_d y_d), $ а для $ x \in \R^d $ и $ A \subset \R^d $ определим
$$
x A = x \cdot A = \{xy: y \in A\}.
$$

Будем обозначать
$$
(x,y) = \sum_{j=1}^d x_j y_j, x, y \in \R^d,
$$
$$
x^{\lambda}  = x_1^{\lambda_1} \ldots x_d^{\lambda_d}, x \in \R^d,
\lambda \in \Z_+^d. $$

При $ d \in \N $ для $ t \in \R^d $ через $ 2^t $ будем обозначать
вектор $ 2^t = (2^{t_1}, \ldots, 2^{t_d}). $

Обозначим через $ \R_+^d $ множество $ x \in \R^d, $ для которых $
x_j >0 $ при $ j=1,\ldots,d,$ и для $ a \in \R_+^d, x \in \R^d $
положим $ a^x = a_1^{x_1} \ldots a_d^{x_d}.$

При $ d \in \N $ определим множества
\begin{multline*}
I^d  =  \{x \in \R^d: 0 < x_j < 1,j=1,\ldots,d\},\\
\overline I^d  =  \{x \in \R^d: 0 \le x_j \le 1,j=1,\ldots,d\},\\
B^d  =  \{x \in \R^d: -1 \le x_j \le 1,j=1,\ldots,d\}.
\end{multline*}

Через $ \e $ будем обозначать вектор в $ \R^d, $ задаваемый
равенством $ \e =(1,\ldots,1). $

При $ d \in \N $ для $ x \in \R^d $ и $ J = \{j_1,\ldots,j_k\}
\subset \N: 1 \le j_1 < j_2 < \ldots < j_k \le d, $ через $ x^J $
обозначим вектор $ x^J = (x_{j_1},\ldots,x_{j_k}) \in \R^k, $ а
для множества $ A \subset \R^d $ положим $ A^J = \{x^J: x \in A\}. $

При $ d \in \N $ для $ x \in \R^d $ через $\s(x) $ обозначим
множество $ \s(x) = \{j =1,\ldots,d: x_j \ne 0\}. $

Будем также обозначать через $ \chi_A $ характеристическую функцию
множества $ A \subset \R^d. $

Для $ d \in \N, x,y \in \R^d $ будем писать $ x \le y (x < y), $
если для каждого $ j=1,\ldots,d $ выполняется неравенство $ x_j
\le y_j (x_j < y_j). $

Для банахова пространства $ X $ (над $ \C$) обозначим $ B(X) = \{x
\in X: \|x\|_X \le 1\}. $

Для банаховых пространств $ X,Y $ через $ \mathcal B(X,Y) $
обозначим банахово пространство, состоящее из непрерывных линейных
операторов $ T: X \mapsto Y, $ с нормой
$$
\|T\|_{\mathcal B(X,Y)} = \sup_{x \in B(X)} \|Tx\|_Y.
$$
Отметим, что если $ X=Y,$ то $ \mathcal B(X,Y) $ является
банаховой алгеброй. Отметим еще, что при $ Y =  \C $ пространство
$ \mathcal B(X, \C) $ обозначается также $ X^*. $

В завершении этого пункта приведем сведения о кратных
рядах, которыми будем пользоваться в дальнейшем.

При $ d \in \N $ для $ y \in \R^d $ обозначим
\begin{eqnarray*}
\mx(y) &  = & \max_{j=1,\ldots,d} y_j,\\
\mn(y) & = & \min_{j=1,\ldots,d} y_j
\end{eqnarray*}
и для банахова пространства $ X, $ вектора $ x \in X $ и семейства
$ \{x_\kappa \in X, \kappa \in \Z_+^d\} $ будем писать $ x =
\lim_{ \mn(\kappa) \to \infty} x_\kappa, $ если для любого $
\epsilon >0 $ существует $ n_0 \in \N $ такое, что для любого $
\kappa \in \Z_+^d, $ для которого $ \mn(\kappa) > n_0, $
справедливо неравенство $ \|x -x_\kappa\|_X  < \epsilon. $

Пусть $ X $ -- банахово пространство (над $ \C $), $ d \in \N $ и
$ \{ x_\kappa \in X: \kappa \in \Z_+^d\} $ -- семейство векторов.
Тогда под суммой ряда $ \sum_{\kappa \in \Z_+^d} x_\kappa $ будем
понимать вектор $ x \in X, $ для которого выполняется равенство $
x = \lim_{\mn(k) \to \infty} \sum_{\kappa \in \Z_+^d(k)} x_\kappa. $

При $ d \in \N $ через $ \Upsilon^d $ обозначим множество
$$
\Upsilon^d = \{ \epsilon \in \Z^d: \epsilon_j \in \{0,1\},
j=1,\ldots,d\}.
$$

Имеет место

   Предложение  1.1.1

Пусть $ X $ -- банахово пространство, а вектор $ x \in X $ и
семейство $ \{x_\kappa \in X: \kappa \in \Z_+^d\} $ таковы, что $
x = \lim_{ \mn(\kappa) \to \infty} x_\kappa, $ Тогда для семейства $
\{ \mathcal X_\kappa \in X, \kappa \in \Z_+^d \}, $ определяемого
равенством
$$
\mathcal X_\kappa = \sum_{\epsilon \in \Upsilon^d: \s(\epsilon)
\subset \s(\kappa)} (-\e)^\epsilon x_{\kappa -\epsilon}, \kappa
\in \Z_+^d,
$$
справедливо равенство
$$
x = \sum_{\kappa \in \Z_+^d} \mathcal X_\kappa.
$$

Предложение является следствием того, что при $ k \in \Z_+^d $
выполняется равенство
\begin{equation*} \tag{1.1.1}
\sum_{\kappa \in \Z_+^d(k)}
\mathcal X_\kappa = x_k \text{(см. [10])}.
\end{equation*}

Замечание.

Легко заметить, что для любого семейства чисел $ \{x_\kappa \in
\R: x_\kappa \ge 0, \kappa \in \Z_+^d\} , $ если ряд $
\sum_{\kappa \in \Z_+^d} x_\kappa $ сходится, т.е. существует
предел $ \lim\limits_{\mn(k) \to \infty} \sum_{\kappa \in
\Z_+^d(k)} x_\kappa, $ то для любой последовательности подмножеств
$ \{Z_n \subset \Z_+^d, n \in \Z_+\}, $
 таких, что $ \card Z_n < \infty,
Z_n \subset Z_{n+1},  n \in \Z_+, $
и $ \cup_{ n \in \Z_+} Z_n = \Z_+^d, $
справедливо равенство
$ \sum_{\kappa \in \Z_+^d} x_\kappa =
\lim_{ n \to \infty} \sum_{\kappa \in Z_n} x_\kappa. $
Отсюда несложно понять, что если для семейства векторов
$ \{x_\kappa \in X, \kappa \in \Z_+^d\}  $ банахова пространства $ X $
ряд $ \sum_{\kappa \in \Z_+^d} \| x_\kappa \|_X $ сходится,
то для любой последовательности подмножеств
$ \{Z_n \subset \Z_+^d, n \in \Z_+\}, $
 таких, что $ \card Z_n < \infty,
Z_n \subset Z_{n+1},  n \in \Z_+, $
и $ \cup_{ n \in \Z_+} Z_n = \Z_+^d, $
в $ X $ соблюдается равенство
$ \sum_{\kappa \in \Z_+^d} x_\kappa =
\lim_{ n \to \infty} \sum_{\kappa \in Z_n} x_\kappa. $

\bigskip

1.2. В этом пункте описываются некоторые свойства операторов проектирования,
рассматриваемых в работе. Все сведения, приводимые в этом пункте без
доказательвта, можно найти в [6], [7].

Прежде всего, отметим, что если $ H:  X \mapsto Y $ -- изоморфизм линейного
пространства $ X $ на линейное пространство $ Y $ и
$ P: Y \mapsto Y $ -- проекционный оператор, то $ H^{-1} P H $ --
проектор, действующий в $ X $ и
\begin{equation*} \tag{1.2.1}
\Im H^{-1} P H = H^{-1} \Im P,
\end{equation*}
\begin{equation*} \tag{1.2.2}
\Ker H^{-1} P H = H^{-1} \Ker P.
\end{equation*}

Лемма 1.2.1

Пусть $ P_0, P_1 $ -- операторы проектирования, действующеие в
линейном пространстве $ X, $ такие, что
\begin{equation*} \tag{1.2.3}
\Im P_0 \subset \Im P_1,
\Ker P_1 \subset \Ker P_0.
\end{equation*}
Тогда $ P_1 -P_0 $ является оператором
проектирования, причем,
\begin{equation*} \tag{1.2.4}
\Im (P_1 -P_0) = \Im P_1 \cap \Ker P_0;
\Ker (P_1 -P_0) = \Im P_0 +\Ker P_1.
\end{equation*}

В справедливости этого утверждения можно убедиться, обратившись, например, к
[6], [7] или [11].

При $ \delta \in \R_+ $ обозначим через $ h_\delta $ отображение, которое каждой
(комплекснозначной) функции $ f  $, заданной на некотором множестве
$ S  \subset \R, $ ставит в соответствие функцию $ h_\delta f, $ определяемую
на множестве $ \{ x \in \R: \delta x \in S\} $ равенством
$ (h_\delta f)(x) = f(\delta x). $
Ясно, что отображение $ h_\delta $
является биекцией на себя  множества всех функций с  областью  определения
в $ \R. $ При этом,  для $ \delta, \sigma \in \R_+ $ имеет место равенство
$ h_{\delta \sigma} = h_\delta h_\sigma $ и
обратное отображение $ h_\delta^{-1} $  задается
равенством
$ h_\delta^{-1}  = h_{\delta^{-1}}.$

Отметим, что при $ 1 \le p  \le \infty $ отображение $ h_\delta, $
суженное на $ L_p(\R), $ является линейным гомеоморфизмом
пространства $ L_p(\R) $ на себя и для $ f \in L_p(\R) $ выполняется равенство
\begin{equation*}
\| h_\delta f\|_{L_p(\R)} =
\delta^{-p^{-1}} \|f\|_{L_p(\R)},
\end{equation*}
а, следовательно, для $ f \in L_p(\R) $ соблюдается равенство
\begin{equation*}
\|h_\delta^{-1}f\|_{L_p(\R)} =
\delta^{p^{-1}} \|f\|_{L_p(\R)}.
\end{equation*}

Кроме того, при $ \delta \in \R_+, 1 \le p \le \infty $ для $ f \in L_p(\R),
g \in L_{p^\prime}(\R) (p^\prime = p/(p -1)) $  справедливо соотношение
$$
\int_\R (h_\delta f) \overline {(h_\delta g)} dx = \int_\R f(\delta x)
\overline {g(\delta x)} dx =
\delta^{-1} \int_\R f(x) \overline {g(x)} dx
$$
а, следовательно,
\begin{equation*} \tag{1.2.5}
\int_\R (h_\delta f) \overline {(h_\delta g)} dx = 0 \iff
\int_\R f \overline g dx =0.
\end{equation*}

Для семейства $ \{ x_\nu \in X, \nu \in \Z \} $ векторов банахова пространства
$ X $ и подмножества $ \Nu \subset \Z $ под суммой $ \sum_{\nu \in \Nu} x_\nu $ будем понимать
$$
\lim_{m, n \in \Z_+, \min(m, n) \to \infty} \sum_{\nu \in \Nu: -m \le \nu \le n} x_\nu.
$$

Через $ l_2, $ как обычно, обозначается гильбертово пространство векторов
$ x = \{x_\nu \in \C, \nu \in \Z\}, $ для которых ряд
$ \sum_{\nu \in \Z} | x_\nu |^2 $ сходится, со скалярным произведением
$$
\langle x, y \rangle = \sum_{\nu \in \Z } x_\nu \overline {y_\nu}, x, y \in l_2.
$$

Пусть $ X_0 $ -- замкнутое линейное подпространство в $ L_2(\R), $
для которого существует функция $ \phi \in X_0 $
такая, что система функций
$ \{ \phi(\cdot -\nu), \nu \in \Z \} $ является базисом Рисса в $ X_0, $
т.е. замыкание в $ L_2(\R) $ линейной оболочки этой системы
\begin{equation*} \tag{1.2.6}
\close_{L_2(\R)} (\span \{\phi(\cdot -\nu), \nu \in \Z \}) = X_0
\end{equation*}
и существуют константы $ A, B \in \R_+, $
обладающие тем свойством, что для любого семейства
$ \{ c_\nu, \nu \in \Z \} \in l_2 $ соблюдаются неравенства
\begin{equation*} \tag{1.2.7}
A \| \{c_\nu, \nu \in \Z \} \|_{l_2} \le
\| \sum_{\nu \in \Z} c_\nu \phi(\cdot -\nu) \|_{L_2(\R)} \le
B \| \{c_\nu, \nu \in \Z \} \|_{l_2} \text{(см. [6], [7]).}
\end{equation*}
Условия (1.2.6), (1.2.7) эквивалентны тому, что отображение
$ A­_\phi: l_2 \mapsto X_0, $ которое каждому
вектору $ \{c_\nu, \nu \in \Z \} \in l_2 $ ставит в соответствие
функцию $ \sum_{\nu \in \Z } c_\nu \phi(\cdot -\nu) \in X_0 \cap L_2(\R) $
является линейным гомеоморфизмом пространства $ l_2 $ на
подпространство $ X_0 \cap L_2(\R). $

Для функции $ \phi \in X_0, $ удовлетворяющей условиям (1.2.6), (1.2.7),
найдем двойственную ей функцию $ \tilde \phi \in X_0$ (см. гл. 1 из [6]),
т.е. такую, что при $ \nu, \mu \in \Z $ соблюдаются равенства
\begin{equation*} \tag{1.2.8}
\int_\R \phi(x -\nu) \overline {\tilde \phi(x -\mu)} dx = \begin{cases} 0,
\text{ при } \nu \in \Z \setminus \{\mu\}; \\
1, \text{при } \nu = \mu.
\end{cases}
\end{equation*}

При этом семейство функций $ \{\tilde \phi(\cdot -\nu), \nu \in \Z \} $ также является базисом Рисса в $ X_0 $ и
для каждого $ f \in X_0 $ в $ L_2(\R) $ имеет место представление
\begin{equation*} \tag{1.2.9}
f = \sum_{\nu \in \Z}  c_\nu \phi(\cdot -\nu),
\end{equation*}
где семейство $ \{ c_\nu, \nu \in \Z \} \in l_2, $ а для каждого
$ \nu \in \Z $ выполняется равенство
\begin{equation*} \tag{1.2.10}
c_\nu = \int_\R f(x) \overline {\tilde \phi(x -\nu)} dx.
\end{equation*}

Обозначим через $ U_0 $ оператор ортогонального проектирования
пространства $ L_2(\R) $ на подпространство $ X_0, $
т.е. $ U_0: L_2(\R) \mapsto L_2(\R) $ -- непрерывный проектор, у которого
\begin{equation*} \tag{1.2.11}
\Im U_0 = X_0,
\end{equation*}
\begin{equation*} \tag{1.2.12}
\Ker U_0 = \Im (E -U_0) = \{ f \in L_2(\R): \int_\R f \overline g\, dx =0\
\ \forall g \in X_0 \},
\end{equation*}
где $ E $ -- тождественный оператор в $ L_2(\R). $

Из (1.2.9), (1.2.10), (1.2.11), (1.2.12) вытекает, что для $ f \in L_2(\R) $
в $ L_2(\R) $ справедливо равенство
$$
U_0 f = \sum_{ \nu \in \Z} \biggl(\int_\R (U_0 f)(x) \overline {\tilde \phi(x -\nu)} dx\biggr)
\phi(\cdot -\nu),
$$
причем,
\begin{multline*}
\int_\R (U_0 f)(x) \overline {\tilde \phi(x -\nu)} dx =
\int_\R ((U_0 f)(x) -f(x) +f(x)) \overline {\tilde \phi(x -\nu)} dx \\ =
\int_\R ((U_0 f)(x) -f(x)) \overline {\tilde \phi(x -\nu)} dx +
\int_\R f(x) \overline {\tilde \phi(x -\nu)} dx =
\int_\R f(x) \overline {\tilde \phi(x -\nu)} dx,
\end{multline*}
и, значит, в $ L_2(\R) $ имеет место равенство
\begin{equation*} \tag{1.2.13}
U_0 f = \sum_{ \nu \in \Z} \biggl(\int_\R f(x) \overline {\tilde \phi(x -\nu)} dx\biggr)
\phi(\cdot -\nu).
\end{equation*}

Теперь при $ k \in \N $ положим
\begin{equation*} \tag{1.2.14}
X_k = h_2^k ( X_0 ) = h_2 (h_2^{k -1} X_0) = h_2 ( X_{k -1}),
\end{equation*}
\begin{equation*} \tag{1.2.15}
U_k = h_2^k U_0 (h_2^k)^{-1} =
h_{2^k} U_0 (h_{2^k})^{-1} = h_{2^k} U_0 h_{2^{-k}}.
\end{equation*}

Сопоставляя (1.2.1), (1.2.2), (1.2.14), (1.2.15), (1.2.11), (1.2.12), (1.2.5), (1.2.6),
заключаем, что $ U_k: L_2(\R) \mapsto L_2(\R) $ является непрерывным
проектором, у которого
\begin{multline*} \tag{1.2.16}
\Im U_k = X_k = h_{2^k} (\close_{L_2(\R)} (\span \{\phi(\cdot -\nu),
\nu \in \Z \})) \\ =
\close_{L_2(\R)} (h_{2^k} (\span \{\phi(\cdot -\nu),
\nu \in \Z \})) \\ =
\close_{L_2(\R)} (\span \{ h_{2^k} (\phi(\cdot -\nu)),
\nu \in \Z \}) \\ =
\close_{L_2(\R)} (\span \{ \phi(2^k \cdot -\nu)),
\nu \in \Z \}),
\end{multline*}
\begin{multline*} \tag{1.2.17}
\Ker U_k = h_{2^k} (\{ f \in L_2(\R): \int_\R f \overline g dx =0
\ \forall g \in X_0 \}) \\ =
\{ h_{2^k}(f): f \in L_2(\R), \int_\R h_{2^k}(f)
\overline { h_{2^k}(g)} dx =0
\ \forall g \in X_0 \}) \\ =
\{ F \in L_2(\R): \int_\R F \overline G dx =0 \ \forall G \in X_k \}.
\end{multline*}

Принимая во внимание (1.2.16), (1.2.17) и учитывая, что для $ f \in L_2(\R) $
справедливо соотношение
\begin{multline*}
\int_\R f(x) \overline {g(x)} dx =0 \ \forall g \in
\close_{L_2(\R)} (\span \{ \phi(2^k \cdot -\nu), \nu \in \Z \}) \\ \iff
\int_\R f(x) \overline {g(x)} dx =0 \ \forall g \in
\span \{ \phi(2^k \cdot -\nu), \nu \in \Z \} \\ \iff
\int_\R f(x) \overline {\phi(2^k x -\nu)} dx =0 \ \forall \nu \in \Z,
\end{multline*}
находим, что при  $ k \in \Z_+ $ верно равенство
\begin{equation*}
\Ker U_k = \{ f \in L_2(\R): \int_\R f(x) \overline {\phi(2^k x -\nu)} dx =0
\ \forall \nu \in \Z \}.
\end{equation*}

Отметим, что ввиду (1.2.15), (1.2.13) для $ f \in L_2(\R) $
при $ k \in \Z_+ $ в $ L_2(\R) $ выполняется равенство
\begin{multline*} \tag{1.2.18}
U_k f = h_{2^k} \biggl(\sum_{ \nu \in \Z} \biggl(\int_\R ( h_{2^{-k}} f)(x)
\overline {\tilde \phi(x -\nu)} dx\biggr) \phi(\cdot -\nu)\biggr) \\ =
\sum_{ \nu \in \Z} \biggl(\int_\R f(2^{-k} x)
\overline {\tilde \phi(x -\nu)} dx\biggr) \phi(2^k \cdot -\nu) \\ =
\sum_{ \nu \in \Z} 2^k \biggl(\int_\R f(x)
\overline {\tilde \phi(2^k x -\nu)} dx\biggr) \phi(2^k \cdot -\nu).
\end{multline*}

Предложение 1.2.2

Пусть $ \phi $ -- измеримая функция на $ \R, $ для которой
существует неотрицательная суммируемая на $ \R $ функция $ \Phi $ такая, что
почти для всех $ x \in \R $ выполняется неравенство
\begin{equation*} \tag{1.2.19}
| \phi(x)| \le \int_{(1/2) B^1} \Phi(x -u) du,
\end{equation*}
и функция $ \tilde \phi \in L_p(\R) $ при $ 1 \le p \le \infty. $
Тогда при $ 1 \le p \le \infty $ функция
\begin{equation*} \tag{1.2.20}
\phi \in L_p(\R),
\end{equation*}
и для $ \kappa \in \Z_+, f \in L_p(\R) $ и любого подмножества
$ \Nu \subset \Z $ ряд
\begin{equation*}
\sum_{ \nu \in \Nu} 2^\kappa \biggl(\int_\R f(y)
\overline {\tilde \phi(2^\kappa y -\nu)} dy\biggr) \phi(2^\kappa x -\nu)
\end{equation*}
абсолютно сходится почти в каждой точке $ x \in \R, $
причем, существует константа $ c_1(\Phi, \tilde \phi, p) >0 $ такая, что
при $ 1 \le p \le \infty $ для $ f \in L_p(\R), \kappa \in \Z_+, \Nu \subset \Z $
почти для всех $ x \in \R $ соблюдается неравенство
\begin{equation*} \tag{1.2.21}
\biggl| \sum_{ \nu \in \Nu} 2^\kappa \biggl(\int_\R f(y)
\overline {\tilde \phi(2^\kappa y -\nu)} dy\biggr) \phi(2^\kappa x -\nu) \biggr| \le
c_1 2^{\kappa /p} \| f \|_{L_p(\R)},
\end{equation*}
а при $ \kappa \in \Z_+, 1 \le p \le \infty $ линейный оператор
$ E_\kappa^p = E_\kappa^{\phi,\tilde \phi, p}: L_p(\R) \mapsto L_\infty(\R), $
который определяется равенством
\begin{multline*} \tag{1.2.22}
(E_\kappa^p f)(x) = \sum_{\nu \in \Z} 2^{\kappa}
\biggl(\int_{\R} f(y) \overline {\tilde \phi(2^\kappa y -\nu)} dy\biggr)
\phi(2^\kappa x -\nu),\\ \text{ почти для всех } x \in \R, f \in L_p(\R),
\end{multline*}
где суммирование производится почти в каждой точке $ x \in \R $, является
непрерывным, и справедлива оценка
\begin{equation*} \tag{1.2.23}
\| E_\kappa^p f \|_{L_\infty(\R)} \le c_1 2^{\kappa /p} \| f \|_{L_p(\R)}.
\end{equation*}

Отметим, что если при $ 1 \le p, q \le \infty $ функция $ f \in L_p(\R) \cap L_q(\R), $
то $ E_\kappa^p f = E_\kappa^q f, \kappa \in \Z_+. $
Поэтому, если указание какого-либо индекса в обозначении по контексту не
существенно, то мы будем его опускать.

Доказательство.

Прежде всего, отметим, что в силу (1.2.19) почти для всех $ x \in \R $
справедливо неравенство
$$
| \phi(x)| \le \int_{(1/2) B^1} \Phi(x -u) du \le
\int_\R \Phi(x -u) du = \int_\R \Phi(u) du =
\|\Phi\|_{L_1(\R)},
$$
т.е. имеет место (1.2.20) при $ p = \infty, $
и
\begin{equation*}
\|\phi\|_{L_\infty(\R)} \le \|\Phi\|_{L_1(\R)}.
\end{equation*}

Кроме того, поскольку для $ u \in (1/2) B^1 $ неотрицательная функция $ \Phi(x -u) $
суммируема на $ \R, $ и $ \int_\R \Phi(x -u) dx = \|\Phi\|_{L_1(\R)} $
суммируема на $ (1/2) B^1, $ то по теореме Фубини функция $ \Phi(x -u) $
суммируема на $ \R \times ((1/2) B^1), $ а функция $ \int_{(1/2) B^1} \Phi(x
-u) du $ суммируема на $ \R, $ что в силу (1.2.19) влечет включение (1.2.20) при $ p =1. $

Теперь при $ 1 < p < \infty $ почти для всех $ x \in \R $ имеем
$$
| \phi(x)|^p = | \phi(x)| \cdot | \phi(x)|^{p -1} \le
| \phi(x)| \cdot \|\phi\|_{L_\infty(\R)}^{p -1} \in L_1(\R),
$$
а, значит, справедливо (1.2.20) при $ 1 < p < \infty. $

Далее, заметим, что при $ \kappa \in \Z_+, \nu \in \Z $ в силу (1.2.19) почти 
для всех $ x \in \R $ верно неравенство
\begin{multline*} \tag{1.2.24}
| \phi(2^\kappa x -\nu)| \le
\int_{(1/2) B^1} \Phi(2^\kappa x -\nu -u) du =
\int_{\nu +(1/2) B^1} \Phi(2^\kappa x -u) du \\
=2^\kappa \int_{2^{-\kappa}(\nu +(1/2) B^1)} \Phi(2^\kappa (x -u)) du,
\end{multline*}
а при $ 1 \le p \le \infty $ для $ f \in L_p(\R) $ вследствие неравенства 
Гельдера имеет место оценка
\begin{multline*} \tag{1.2.25}
\biggl| \int_{\R} f(y) \overline {\tilde \phi(2^\kappa y -\nu)} dy\biggr| \le
\| f \|_{L_p(\R)} \cdot
\| \tilde \phi(2^\kappa \cdot -\nu)\|_{L_{p^\prime}(\R)} \\
=2^{-\kappa /p^\prime} \| f \|_{L_p(\R)} \cdot
\| \tilde \phi\|_{L_{p^\prime}(\R)}
= c_2(\tilde \phi, p) 2^{-\kappa /p^\prime} \| f \|_{L_p(\R)}.
\end{multline*}
Из (1.2.24) и (1.2.25), пользуясь $ \sigma $-аддитивностью интеграла как
функции множеств, а также тем, что
$$
\mes(2^{-\kappa}(\nu +(1/2) B^1) \cap 2^{-\kappa}(\nu^\prime +(1/2) B^1)) =0,
\nu, \nu^\prime \in \Z: \nu \ne \nu^\prime, \kappa \in \Z_+,
$$
заключаем, что при $ 1 \le p \le \infty, \kappa \in \Z_+, f \in L_p(\R), \Nu \subset \Z $
почти для всех $ x \in \R $ справедливо соотношение
\begin{multline*}
\biggl| \sum_{ \nu \in \Nu} 2^\kappa (\int_\R f(y)
\overline {\tilde \phi(2^\kappa y -\nu)} dy) \phi(2^\kappa x -\nu) \biggr| \\
\le\sum_{\nu \in \Nu} 2^{\kappa}
\biggl| \int_{\R} f(y) \overline {\tilde \phi(2^\kappa y -\nu)} dy\biggr| \cdot
| \phi(2^\kappa x -\nu)| \\
\le 2^{\kappa} \sum_{\nu \in \Nu}
c_2 2^{-\kappa /p^\prime} \| f \|_{L_p(\R)} \cdot
2^\kappa \int_{2^{-\kappa}(\nu +(1/2) B^1)} \Phi(2^\kappa (x -u)) du \\
=c_2 2^{2\kappa -\kappa /p^\prime} \| f \|_{L_p(\R)} \cdot
\int_{\cup_{\nu \in \Nu} (2^{-\kappa}(\nu +(1/2) B^1))} \Phi(2^\kappa (x -u)) du\\
 \le c_2 2^{2\kappa -\kappa /p^\prime} \| f \|_{L_p(\R)} \cdot
\int_\R \Phi(2^\kappa (x -u)) du =
c_2 2^{2\kappa -\kappa /p^\prime} \| f \|_{L_p(\R)} \cdot
\int_\R \Phi(2^\kappa u) du \\
= c_2 2^{2\kappa -\kappa /p^\prime} \| f \|_{L_p(\R)} \cdot
\int_\R \Phi(u) 2^{-\kappa} du =
c_2 2^{\kappa -\kappa /p^\prime} \| f \|_{L_p(\R)} \cdot
\int_\R \Phi(u) du \\
= c_1(\Phi, \tilde \phi, p) 2^{\kappa /p} \| f \|_{L_p(\R)},
\end{multline*}

т.е. имеет место (1.2.21).
В частности, при $ \Nu = \Z $ из (1.2.21) следует, что при
$ 1 \le p \le \infty, \kappa \in \Z_+ $ для $ f \in L_p(\R) $ ряд в правой части
(1.2.22) сходится почти для всех $ x \in \R $ и
\begin{equation*}
\biggl| \sum_{\nu \in \Z} 2^{\kappa}
\biggl(\int_{\R} f(y) \overline {\tilde \phi(2^\kappa y -\nu)} dy\biggr) \cdot
\phi(2^\kappa x -\nu) \biggr| \le
c_1 2^{\kappa /p} \| f \|_{L_p(\R)} \text{ почти всюду на } \R,
\end{equation*}
откуда с учетом (1.2.22) вытекает (1.2.23). $ \square $

Предложение 1.2.3

Пусть $ X_0 $ -- замкнутое линейное подпространство в $ L_2(\R), $ такое,
что существуют функции $ \phi \in X_0, \tilde \phi \in X_0, $ для которых соблюдаются
условия (1.2.6), (1.2.7), (1.2.8), а
также выполняются условия предложения 1.2.2.
Тогда

1) при $ 1 \le p \le \infty, \kappa \in \Z_+ $ для $ f \in L_2(\R) \cap L_p(\R) $
почти для всех $ x \in \R $ имеет место равенство
\begin{equation*} \tag{1.2.26}
(E_\kappa^p f)(x) = (U_\kappa f)(x) (\text{см.} (1.2.22), (1.2.15));
\end{equation*}

2) при $ 1 \le p, q \le \infty, \kappa \in \Z_+ $ для $ f \in L_p(\R) \cap L_2(\R),
g \in L_q(\R) \cap L_2(\R) $ соблюдаются равенства
\begin{equation*} \tag{1.2.27}
\int_{\R} (E_\kappa^p f) \cdot \overline g dx = \int_{\R} f \cdot
\overline {(E_\kappa^q g) } dx
\end{equation*}
и
\begin{equation*} \tag{1.2.28}
\int_{\R} (\mathcal E_\kappa^p f) \cdot \overline g dx = \int_{\R} f \cdot
\overline {(\mathcal E_\kappa^q g) } dx,
\end{equation*}
где оператор $ \mathcal E_\kappa^p = \mathcal E_\kappa^{\phi, \tilde \phi, p}:
L_p(\R) \mapsto L_\infty(\R), $ определяется равенством
$$
\mathcal E_\kappa^p = E_\kappa^p -E_{\kappa -1}^p, \kappa \in \Z_+,
E_{-1}^p =0, 1 \le p \le \infty.
$$

Доказательство.

Чтобы убедиться в справедливости (1.2.26), заметим, что для каждой
функции $ f \in L_2(\R) $ в силу (1.2.18) последовательность функций
$$
\biggl\{ \sum_{\nu =-n}^n 2^{\kappa}
\biggl(\int_{\R} f(y) \overline {\tilde \phi(2^\kappa y -\nu)} dy\biggr) \cdot
\phi(2^\kappa \cdot -\nu): n \in \Z_+\biggr\}
$$
в $ L_2(\R) $ сходится к $ (U_\kappa f)(\cdot) $ при $ n \to \infty. $
Поэтому существует подпоследовательность
$$
\biggl\{ \sum_{\nu =-n_k}^{n_k} 2^{\kappa}
\biggl(\int_{\R} f(y) \overline {\tilde \phi(2^\kappa y -\nu)} dy\biggr) \cdot
\phi(2^\kappa x -\nu): n_k \in \Z_+, n_k < n_{k +1}, k \in \N\biggr\}
$$
которая почти для всех $ x \in \R $ сходится к $ (U_\kappa f)(x) $
при $ k \to \infty. $
Отсюда, учитывая, что для $ f \in L_p(\R) $ в виду (1.2.22) семейство функций
$$
\biggl\{ \sum_{\nu =-m}^n 2^{\kappa}
\biggl(\int_{\R} f(y) \overline {\tilde \phi(2^\kappa y -\nu)} dy\biggr) \cdot
\phi(2^\kappa x -\nu): m,n \in \Z_+\biggr\}
$$
почти для всех $ x \in \R $ сходится к $ (E_\kappa^p f)(x) $ при $ \min(m,n) \to \infty, $
заключаем, что для $ f \in L_2(\R) \cap L_p(\R) $ почти для всех $ x \in \R $
имеет место (1.2.26).

Для получения (1.2.27), принимая во внимание (1.2.26), (1.2.16), (1.2.17),
а также то обстоятельство, что $ \Im (E -U_\kappa) = \Ker U_\kappa, $
для $ f \in L_p(\R) \cap L_2(\R), g \in L_q(\R) \cap L_2(\R) $ имеем
\begin{multline*}
\int_{\R} (E_\kappa^p f) \cdot \overline g dx =
\int_{\R} (U_\kappa f) \cdot \overline g dx =
\int_{\R} (U_\kappa f) \cdot \overline {( U_\kappa g +g -U_\kappa g)} dx \\ =
\int_{\R} (U_\kappa f) \cdot \overline { U_\kappa g } dx +
\int_{\R} (U_\kappa f) \cdot \overline {( g -U_\kappa g)} dx =
\int_{\R} (U_\kappa f) \cdot \overline { (U_\kappa g )} dx,
\end{multline*}
и
\begin{multline*}
\int_{\R} f \cdot \overline {(E_\kappa^q g) } dx =
\int_{\R} f \cdot \overline {(U_\kappa g) } dx =
\int_{\R} (U_\kappa f +f -U_\kappa f) \cdot \overline {(U_\kappa g) } dx\\ =
\int_{\R} (U_\kappa f ) \cdot \overline {(U_\kappa g) } dx +
\int_{\R} (f -U_\kappa f) \cdot \overline {(U_\kappa g) } dx \\ =
\int_{\R} (U_\kappa f ) \cdot \overline {(U_\kappa g) } dx +
\overline{ \int_{\R} (U_\kappa g) \cdot \overline{ (f -U_\kappa f)} dx} \\ =
\int_{\R} (U_\kappa f ) \cdot \overline {(U_\kappa g) } dx.
\end{multline*}
Сопоставляя эти равенства, видим, что (1.2.27) выполняется для
$ f \in L_p(\R) \cap L_2(\R),
g \in L_q(\R) \cap L_2(\R). $

Для вывода (1.2.28), используя (1.2.27), получаем
\begin{multline*}
\int_{\R} (\mathcal E_\kappa^p f) \cdot \overline g dx =
\int_{\R} ((E_\kappa^p f)
-(E_{\kappa -1}^p f))
\cdot \overline g dx =
\int_{\R} (E_\kappa^p f) \cdot \overline g dx -
\int_{\R} (E_{\kappa -1}^p f) \cdot \overline g dx\\ =
\int_{\R} f \cdot \overline {(E_\kappa^q g) } dx -
\int_{\R} f \cdot \overline {(E_{\kappa -1}^q g) } dx\\ =
\int_{\R} f \cdot (\overline {(E_\kappa^q g) } -
\overline {(E_{\kappa -1}^q g) }) dx\\ =
\int_{\R} f \cdot \overline {(\mathcal E_\kappa^q g) } dx,
f \in L_p(\R) \cap L_2(\R), g \in L_q(\R) \cap L_2(\R). \square
\end{multline*}

Предложение 1.2.4

Пусть выполнены условия предложения 1.2.3 (без соблюдения условий
предложения 1.2.2), а также $ X_0 \subset X_1 = h_2(X_0) $
и $ \cup_{\kappa \in \Z_+} \span \{ \phi(2^\kappa \cdot -\nu),
\nu \in \Z \} $ плотно в $ L_2(\R). $ Тогда имеют место соотношения:

1) \begin{equation*} \tag{1.2.29}
\Im U_\kappa \subset \Im U_{\kappa +1}; \\
\Ker U_{\kappa +1} \subset \Ker U_\kappa, \kappa \in \Z_+;
\end{equation*}

2) для $ \kappa, \kappa^\prime \in \Z_+: \kappa^\prime \le \kappa, $ выполняются равенства
\begin{equation*} \tag{1.2.30}
U_{\kappa^\prime} U_\kappa = U_\kappa
U_{\kappa^\prime} =
U_{\kappa^\prime};
\end{equation*}

3) для $ \kappa, \kappa^\prime \in \Z_+ $ соблюдаются равенства
\begin{equation*} \tag{1.2.31}
\mathcal U_\kappa \mathcal U_{\kappa^\prime} =
\begin{cases} \mathcal U_\kappa,
  \text{при $ \kappa = \kappa^\prime $};\\
       0,  \text{при $ \kappa \ne \kappa^\prime $};
\end{cases}
\end{equation*}
где $ \mathcal U_\kappa = U_\kappa -U_{\kappa -1}, \kappa \in \Z_+, U_{-1} =0, $
и при $ \kappa, \kappa^\prime \in \Z_+: \kappa \ne \kappa^\prime, $
для $ f \in L_2(\R), g \in L_{2}(\R) $ имеет место равенство
\begin{equation*} \tag{1.2.32}
\int_{\R} (\mathcal U_\kappa f) \cdot \overline (\mathcal
U_{\kappa^\prime} g) dx = 0;
\end{equation*}

4) для $ \kappa \in \N $ соблюдаются равенства
$$
\Im \mathcal U_\kappa = \Im U_\kappa \cap \Ker U_{\kappa -1};
\Ker \mathcal U_\kappa = \Im U_{\kappa -1} +\Ker U_\kappa;
$$

5) в $ L_2(\R) $ справедливо равенство
\begin{equation*} \tag{1.2.33}
f = \sum_{ \kappa \in \Z_+} \mathcal U_\kappa f, f \in L_2(\R),
\end{equation*}
и для любой функции $ f \in L_2(\R) $ имеет место равенство
\begin{equation*} \tag{1.2.34}
\| f \|_{L_2(\R)} = \biggl(\sum_{ \kappa \in \Z_+} \| \mathcal
U_\kappa f \|_{L_2(\R)}^2\biggr)^{1/2}.
\end{equation*}

Доказательство.

Прежде всего, в условиях предложения с учетом (1.2.16), (1.2.14) имеем
\begin{multline*} 
\Im U_\kappa = X_\kappa = h_{2^\kappa} (X_0) \subset h_{2^\kappa} (X_1) =
h_{2^\kappa} (h_2 (X_0)) = (h_{2^\kappa} h_2) (X_0) \\ =
h_{2^{\kappa +1}} (X_0) = X_{\kappa +1} = \Im U_{\kappa +1}.
\end{multline*}

А ввиду (1.2.17), (1.2.14) и первого включения из (1.2.29) выводим
\begin{multline*}
\Ker U_{\kappa +1} =
\{ F \in L_2(\R): \int_\R F \overline G dx =0 \ \forall G \in \Im U_{\kappa +1} \} \\
\subset \{ F \in L_2(\R): \int_\R F \overline G dx =0 \ \forall G \in \Im U_\kappa \}
= \Ker U_\kappa,
\end{multline*}
что завершает вывод (1.2.29).

Убедимся в справедливости (1.2.30). Поскольку при $ \kappa^\prime,
\kappa \in \Z_+: \kappa^\prime \le \kappa, $
вследствие (1.2.29) имеет место включение
$ \Im U_{\kappa^\prime} \subset \Im U_\kappa, $
то для проектора $ U_\kappa $ вытекает, что
$$
U_\kappa (U_{\kappa^\prime} f) =
U_{\kappa^\prime} f,  f \in L_2(\R).
$$

Далее, в условиях п. 2 предложения для $ f \in L_2(\R) $ находим, что
$$
U_{\kappa^\prime} f = U_{\kappa^\prime} (U_\kappa
f) +U_{\kappa^\prime} (f -U_\kappa f) =
U_{\kappa^\prime} (U_\kappa f),
$$
ибо вследствие (1.2.29) справедливо соотношение
\begin{equation*}
\Im (E -U_\kappa ) = \Ker U_\kappa
\subset \Ker U_{\kappa^\prime},
\end{equation*}
и, значит, выполняется равенство
$$
U_{\kappa^\prime} (f -U_\kappa f) =0.
$$

Перейдем к проверке (1.2.31).
Пусть $ \kappa, \kappa^\prime \in \Z_+ $ и $ \kappa = \kappa^\prime. $
Тогда в силу (1.2.30) получаем
\begin{multline*}
(\mathcal U_\kappa)^2 = (U_\kappa)^2 - U_{\kappa
-1} U_\kappa -U_\kappa U_{\kappa -1}
+(U_{\kappa -1})^2 =\\
= U_\kappa -U_{\kappa -1} -U_{\kappa -1}
+U_{\kappa -1} = U_\kappa -U_{\kappa -1} =
\mathcal U_\kappa.
\end{multline*}

Пусть $ \kappa \ne \kappa^\prime. $ Предположим, что $
\kappa^\prime < \kappa. $ Тогда, снова используя (1.2.30), выводим
\begin{multline*}
\mathcal U_{\kappa^\prime} \mathcal u_\kappa =
U_{\kappa^\prime} U_\kappa -U_{\kappa^\prime
-1} U_\kappa -U_{\kappa^\prime} U_{\kappa
-1} +
U_{\kappa^\prime -1} U_{\kappa -1} =\\
U_{\kappa^\prime} -U_{\kappa^\prime -1} -
U_{\kappa^\prime} +U_{\kappa^\prime -1} =0.
\end{multline*}

Аналогично проверяется (1.2.31) при $
\kappa^\prime > \kappa. $

Для проверки (1.2.32), применяя равенство
$$
\int_{\R} (\mathcal U_\kappa F) \cdot \overline G dx = \int_{\R} F \cdot
\overline {(\mathcal U_\kappa G) } dx, F, G \in L_2(\R),
$$
справедливость которого видна из вывода (1.2.27), (1.2.28), с учетом (1.2.31),
для $ f \in L_2(\R), g \in L_{2}(\R) $ находим
\begin{equation*}
\int_{\R} (\mathcal U_\kappa f) \cdot \overline {(\mathcal
U_{\kappa^\prime} g)} dx =
\int_{\R} ( \mathcal U_{\kappa^\prime} \mathcal U_\kappa f) \cdot
\overline g dx =
\int_\R 0 \cdot \overline g dx =0, \kappa, \kappa^\prime \in \Z_+:
\kappa \ne \kappa^\prime.
\end{equation*}

Далее, сопоставляя (1.2.29) с (1.2.3), в соответствии
с (1.2.4) получаем равенства п. 4).

Для доказательства (1.2.33) с учетом предложения 1.1.1 покажем, что
для $ f \in L_2(\R) $ имеет место соотношение
\begin{equation*} \tag{1.2.35}
\| f -U_\kappa f \|_{L_2(\R)} \to 0 \text{ при } \kappa \to \infty.
\end{equation*}

Для $ f \in L_2(\R) $ и произвольного $ \epsilon >0, $ пользуясь плотностью
в $ L_2(\R) $ множества
$$
\cup_{\kappa \in \Z_+} \span \{ \phi(2^\kappa \cdot -\nu),
\nu \in \Z \},
$$
выберем $ \kappa_0 \in \Z_+ $ и $ g \in
\span \{ \phi(2^{\kappa_0} \cdot -\nu), \nu \in \Z \}, $ для которых
$$
\| f -g \|_{L_2(\R)} < \epsilon.
$$
Тогда при $ \kappa > \kappa_0 $ ввиду (1.2.16), (1.2.29) и того факта, что
$ \| U_\kappa F \|_{L_2(\R)} \le \| F \|_{L_2(\R)}, F \in L_2(\R), $
получаем, что
\begin{multline*}
\| f -U_\kappa f \|_{L_2(\R)} = \| f -g +g -U_\kappa f \|_{L_2(\R)} \\ =
\| f -g +U_\kappa g -U_\kappa f \|_{L_2(\R)} =
\| f -g +U_\kappa (g -f) \|_{L_2(\R)} \\ \le
\| f -g \|_{L_2(\R)} +
\| U_\kappa (g -f) \|_{L_2(\R)} \le
2 \| f -g \|_{L_2(\R)} < 2 \epsilon,
\end{multline*}
т.е. имеет место (1.2.35), а вместе с ним и (1.2.33).

Наконец, обращаясь к выводу (1.2.34), заметим, что для любого конечного
множества $ S \subset \Z_+ $ и любого набора функций $ \{ g_\kappa \in
\Im \mathcal U_\kappa, \kappa \in S \} $ верно равенство
\begin{equation*} \tag{1.2.36}
\| \sum_{ \kappa \in S} g_\kappa \|_{L_2(\R)} = \biggl(\sum_{ \kappa
\in S} \| g_\kappa \|_{L_2(\R)}^2\biggr)^{1/2}.
\end{equation*}

В самом деле, ввиду (1.2.32) выводим
\begin{multline*}
\| \sum_{ \kappa \in S} g_\kappa \|_{L_2(\R)}^2 = \int_{\R} \biggl|
\sum_{ \kappa \in S} g_\kappa\biggr|^2 dx =
\int_{\R} ( \sum_{ \kappa \in S} g_\kappa) ( \sum_{ \kappa^\prime \in S}
\overline { g_{\kappa^\prime}}) dx =\\
\int_{\R} \sum_{ \kappa \in S} \sum_{ \kappa^\prime \in S}
g_\kappa \overline {g_{\kappa^\prime}} dx =
\sum_{ \kappa \in S}
 \sum_{ \kappa^\prime \in S}
\int_{\R} g_\kappa \overline {g_{\kappa^\prime}} dx = \sum_{ \kappa \in S}
\int_{\R} | g_\kappa |^2 dx =
\sum_{ \kappa \in S} \| g_\kappa
\|_{L_2(\R)}^2,
\end{multline*}

откуда приходим к (1.2.36).

Для получения равенства (1.2.34), в условиях теоремы, благодаря (1.2.33)
и (1.2.36), имеем
\begin{multline*}
\| f \|_{L_2(\R)} = \biggl\| \lim_{ k \to \infty} \sum_{\kappa =0}^k
\mathcal U_\kappa f \biggr\|_{L_2(\R)}
= \lim_{ k \to \infty} \biggl\|\sum_{\kappa =0}^k \mathcal U_\kappa f \biggr\|_{L_2(\R)}\\
= \lim_{ k \to \infty} \biggl( \sum_{\kappa =0}^k
\| \mathcal U_\kappa f \|_{L_2(\R)}^2\biggr)^{1/2}
= \biggl( \lim_{ k \to \infty} \sum_{\kappa =0}^k
\| \mathcal U_\kappa f \|_{L_2(\R)}^2\biggr)^{1/2}\\ =
\biggl(\sum_{\kappa =0}^\infty
\| \mathcal U_\kappa f \|_{L_2(\R)}^2\biggr)^{1/2}. \square
\end{multline*}

Отметим, что функция $ \phi $ из предложения 1.2.4 является масштабирующей
для соответствующего кратно-масштабного анализа (см., например, [6], [7]).
\bigskip

1.3. В этом пункте приведем некоторые вспомогательные утверждения,
которые используются в дальнейшем.

Доказательство леммы 1.3.1, проведенное в [11], вполне аналогично доказательству
соответствующего утверждения из [12], [13].

Лемма 1.3.1

    Пусть $ d \in \N, 1 \le p < \infty. $ Тогда

1) при  $ j=1,\ldots,d$ для любого непрерывного линейного оператора $ T: L_p(\R)
\mapsto L_p(\R) $ существует единственный непрерывный линейный оператор
$ \mathcal T^j: L_p(\R^d) \mapsto L_p(\R^d), $ для которого для любой функции
$ f \in L_p(\R^d) $ почти для всех $ (x_1,\ldots,x_{j-1},x_{j+1},
\ldots,x_d) \in \R^{d-1} $ в $ L_p(\R) $ выполняется равенство
\begin{multline*} \tag{1.3.1}
(\mathcal T^j f)(x_1,\ldots, x_{j-1},\cdot,x_{j+1},
\ldots,x_d) \\
= (T(f(x_1,\ldots,x_{j-1},\cdot,x_{j+1},
\ldots,x_d)))(\cdot),
\end{multline*}

2) при этом, для каждого $ j=1,\ldots,d $ отображение $ V_j^{L_p}, $ которое
каждому оператору $ T \in \mathcal B(L_p(\R), L_p(\R)) $ ставит в соответствие
оператор $ V_j^{L_p}(T) = \mathcal T^j \in \mathcal B(L_p(\R^d), L_p(\R^d)), $
удовлетворяющий (1.3.1), является непрерывным гомоморфизмом банаховой алгебры
$ \mathcal B(L_p(\R), L_p(\R)) $ в банахову алгебру $ \mathcal B(L_p(\R^d),
L_p(\R^d)), $
и
\begin{equation*} \tag{1.3.2}
\| V_j^{L_p}(T) \| _{\mathcal B(L_p(\R^d), L_p(\R^d))} \le
\| T \|_{ \mathcal B(L_p(\R), L_p(\R))},
\end{equation*}

3) причём, для любых операторов $ S,T \in \mathcal B(L_p(\R), L_p(\R)) $
при любых $ i,j =1,\ldots,d: i\ne j, $ соблюдается равенство
\begin{equation*} \tag{1.3.3}
(V_i^{L_p}(S) V_j^{L_p}(T))f = (V_j^{L_p}(T) V_i^{L_p}(S))f, f \in L_p(\R^d).
\end{equation*}

Замечание.

Если при $ d \in \N, 1 \le p, q < \infty, $ оператор $ T \in
\mathcal B(L_p(\R), L_p(\R)) \cap \mathcal B(L_q(\R), L_q(\R)), $ то
при $ j=1,\ldots,d $ для $ f \in L_p(\R^d) \cap L_q(\R^d)$ справедливо равенство
\newline
$ (V_j^{L_p} T)f = (V_j^{L_q} T)f. $ Поэтому символы $ L_p, L_q $
в качестве индексов у $ V_j $ можно опускать.

Лемма 1.3.2

Пусть $ d \in \N, 1 \le p < \infty $ и $ T_j: L_p(\R) \mapsto
L_p(\R), j =1, \ldots, d $ -- набор непрерывных линейных операторов.
Тогда для любого множества $ \mathcal J \subset \{1, \ldots, d\} $
и любой функции $ f \in L_p(\R^d) $ вида
$ f(x_1, \ldots, x_d) = \prod_{j=1}^d f_j(x_j), $ где $ f_j \in
L_p(\R), j =1, \ldots, d, $ почти для всех $ (x_1, \ldots, x_d) \in \R^d $
выполняется равенство
\begin{multline*} \tag{1.3.4}
\biggl(\biggl(\prod_{j \in \mathcal J} V_j (T_j)\biggr) f\biggr)(x_1, \ldots, x_d) \\=
\biggl(\prod_{j \in \mathcal J} (T_j f_j)(x_j)\biggr) \times
\biggl(\prod_{j =1, \ldots, d: j \notin \mathcal J} f_j(x_j)\biggr).
\end{multline*}

Доказательство.

Доказательство (1.3.4) проведем по индукции относительно $ \card \mathcal J. $
При $ \card \mathcal J = 1, $ т.е. для  $ \mathcal J = \{j\}, j =1, \ldots, d, $
согласно (1.3.1) почти для всех
$ (x_1, \ldots, x_{j -1}, x_{j +1}, \ldots, x_d ) \in \R^{d -1} $
почти в каждой точке $ x_j \in \R $ имеем
\begin{multline*}
(V_j(T_j) f) (x_1, \ldots, x_{j -1}, x_j, x_{j +1}, \ldots, x_d ) \\=
(T_j f(x_1, \ldots, x_{j -1}, \cdot, x_{j +1}, \ldots, x_d))(x_j) \\=
\biggl(T_j \biggl(f_j(\cdot) \times \prod_{i =1, \ldots, d: i \ne j} f_i(x_i)\biggr)\biggr)(x_j) =
(T_j f_j)(x_j) \times \prod_{i =1, \ldots, d: i \ne j} f_i(x_i),
\end{multline*}
что совпадет с (1.3.4) при $ \card \mathcal J =1. $

Предположим, что (1.3.4) имеет место для любого множества $ \mathcal J, $
у которого $ \card \mathcal J \le m < d. $
Покажем, что тогда (.3.14) справедливо при $ \card \mathcal J = m +1. $
Представляя множество $ J \subset \{1, \ldots, d\}: \card J = m +1, $ в виде
$ J = \mathcal J \cup \{j\}, $
где $ j \notin \mathcal J, $ на основании предположения индукции
с учетом (1.3.3) выводим
\begin{multline*}
\biggl(\biggl(\prod_{i \in J} V_i (T_i)\biggr) f\biggr)(x_1, \ldots, x_d) =
\biggl(\biggl(\prod_{i \in \mathcal J \cup \{j\}} V_i(T_i)\biggr) f\biggr)(x_1, \ldots, x_d) \\=
\biggl(V_j(T_j)\biggl(\biggl(\prod_{i \in \mathcal J} V_i (T_i)\biggr) f\biggr)\biggr)(x_1, \ldots, x_d) \\=
\biggl(V_j(T_j)\biggl(\biggl(\prod_{i \in \mathcal J} (T_i f_i)(y_i)\biggr) \times
\biggl(\prod_{i =1, \ldots, d: i \notin \mathcal J} f_i(y_i)\biggr)\biggr)\biggr)(x_1, \ldots, x_d)\\ =
(T_j f_j)(x_j) \times \biggl(\prod_{i \in \mathcal J} (T_i f_i)(x_i)\biggr) \times
\biggl(\prod_{i =1, \ldots, d: i \notin (\mathcal J \cup \{j\})} f_i(x_i)\biggr) \\=
\biggl(\prod_{i \in (\mathcal J \cup \{j\})} (T_i f_i)(x_i)\biggr) \times
\biggl(\prod_{i =1, \ldots, d: i \notin (\mathcal J \cup \{j\})} f_i(x_i)\biggr) \\=
\biggl(\prod_{i \in JJ} (T_i f_i)(x_i)\biggr) \times
\biggl(\prod_{i =1, \ldots, d: i \notin J} f_i(x_i)\biggr). \square
\end{multline*}

\bigskip

1.4. В этом пункте приведены вспомогательные утверждения,
на которые опирается доказательство занимающей центральное место в работе
леммы 2.1.1. Все эти утверждения можно извлечь из [14, гл. I].

Будем обозначать $ |x| = (\sum_{j=1}^d |x_j|^2)^{1/2}, x \in \R^d.$

Для множества $ F \subset \R^d $ и $ x \in \R^d $
положим
$$
\rho(x) = \rho(x,F) = \inf_{y \in F} | x -y|.
$$

Лемма 1.4.1

Пусть $ d \in \N. $ Тогда существует константа $ c_1(d) >0 $
такая, что для любого замкнутого множества $ F \subset \R^d, $
для которого $ \mes (\R^d \setminus F) < \infty, $ функция $ M(x), $
определяемая для $ x \in F $ равенством
$$
M(x) = \int_{\R^d} \rho(u) |x -u|^{-(d+1)} du,
$$
суммируема на $ F $ и соблюдается  неравенство
\begin{equation*} \tag{1.4.1}
\int_F M(x) dx \le c_1 \mes (\R^d \setminus F).
\end{equation*}
(см. \S 2 гл. I из [14])

Предложение 1.4.2

Пусть $ d \in \N. $ Тогда существуют константы $ c_2(d) >0, c_3(d) >0,
c_4(d) >0, c_5(d) >0 $ такие, что
для любой (вещественной) функции $ f \in L_1(\R^d) $ при любом $ \alpha \in \R_+ $
можно построить замкнутое множесво $ F \subset \R^d $ и семейство открытых кубов
$ \{ Q_r, r \in \N \} $ со следующими свойствами:

1) почти для всех $ x \in F $  выполняется неравенство
\begin{equation*} \tag{1.4.2}
| f(x) | \le \alpha;
\end{equation*}

2) для $ W = \R^d \setminus F $ справедливо неравенство
\begin{equation*} \tag{1.4.3}
\mes W \le (c_2 / \alpha) \int_{\R^d} |f| dx;
\end{equation*}

3) \begin{equation*} \tag{1.4.4}
Q_r \cap Q_s = \emptyset, \text{ для}   r,s \in \N: r \ne s;
\end{equation*}

4) \begin{equation*} \tag{1.4.5}
W = \cup_{ r \in \N} \overline Q_r,
\end{equation*}
где  $ \overline Q_r $ -- замыкание куба $ Q_r, r \in \N; $

5) при $ r \in \N $ справедливы неравенства
\begin{equation*} \tag{1.4.6}
c_3 \diam Q_r < \inf_{ x \in Q_r} \rho(x, F) \le
c_4 \diam Q_r,
\end{equation*}
где  $ \diam Q = \sup_{x,y \in Q} |x -y|, Q \subset \R^d; $

а также

6) при  $ r \in \N $ имеет место оценка
\begin{equation*} \tag{1.4.7}
(1 / \mes Q_r) \int_{ Q_r} | f(x) | dx \le c_5 \alpha.
\end{equation*}
(см. \S 3 гл. I из [14] или [5])

При $ d \in \N $ через $ L(\R^d) $ обозначим пространство вещественных
измеримых по Лебегу функций в $ \R^d. $ Как обычно, при $ 1 \le p_1, p_2 \le
\infty $ под суммой $ L_{p_1}(\R^d) +L_{p_2}(\R^d) $ понимается
подпространство в $ L(\R^d), $ состоящее из всех функций $ f \in L(\R^d), $
для которых существуют функции $ f_1 \in L_{p_1}(\R^d) $ и
$ f_2 \in L_{p_2}(\R^d) $ такие, что $ f = f_1 +f_2.$

Напомним (см., например, [14]), что при $ 1 \le p_1 \le p \le p_2 < \infty $ справедливо включение
$ L_p(\R^d) \subset L_{p_1}(\R^d) +L_{p_2}(\R^d). $

Теорема 1.4.3

Пусть $ d \in \N, 1 < q < \infty, C_0 \in \R_+, C_1 \in \R_+ $ и
$ T: (L_1(\R^d) +L_q(\R^d)) \mapsto L(\R^d) $ --
отображение, удовлетворяющее следующим условиям:

1) для любых $ f, g \in (L_1(\R^d) +L_q(\R^d)) $ почти для всех
$ x \in \R^d $ выполняется неравенство
  \begin{equation*} \tag{1.4.8}
| (T(f +g))(x) | \le | (T f)(x) | +| (T g)(x) |;
\end{equation*}

2) для $ f \in L_1(\R^d) $ при $ \alpha >0 $ соблюдается неравенство
\begin{equation*} \tag{1.4.9}
\mes \{ x \in \R^d: | (T f)(x) | > \alpha \} \le
(C_0 / \alpha) \| f \|_{L_1(\R^d)};
\end{equation*}

3) для $ f \in L_q(\R^d) $ при $ \alpha >0 $ выполняется неравенство
\begin{equation*} \tag{1.4.10}
\mes \{ x \in \R^d: | (T f)(x) | > \alpha \} \le
((C_1 / \alpha) \| f \|_{L_q(\R^d)})^q.
\end{equation*}
Тогда при $ 1< p < q $ существует константа $ c_6(p, q, C_0, C_1)
>0 $ такая, что для любого отображения $ T: (L_1(\R^d) +L_q(\R^d))
\mapsto L(\R^d), $ подчиненного условиям (1.4.8) -- (1.4.10), для
$ f \in L_p(\R^d) $ имеет место неравенство
\begin{equation*} \tag{1.4.11}
\| T f \|_{L_p(\R^d)} \le c_6 \| f \|_{L_p(\R^d)}.
\end{equation*}
(см. \S 4 гл. I из [14])
\bigskip

\centerline{\S 2. Свойства ортопроекторов, порожденных КМА}

2.1. В этом пункте устанавливается справедливость утверждения, на котором
основан вывод аналога теоремы Литтлвуда-Пэли для ортопроекторов на
подпространства всплесков $ \{\mathcal E_\kappa, \kappa \in \Z_+^d\} $
(см. теорему 2.4.2), соответствующие кратно-масштабному анализу, порожденному
тензорным произведением гладких достаточно быстро стремящихся к нулю на
бесконечности функций.
При этом, при доказательстве теоремы 2.4.2 будем придерживаться того же
подхода, что в [3, п. 1.5.2] в случае теоремы Литтлвуда-Пэли для кратных рядов
Фурье (см. также [5]).
Убедимся, что имеет место

Лемма 2.1.1

Пусть выполнены условия предложений 1.2.2 и 1.2.4, функция
$ \tilde \phi \in C^1(\R), $ а также
соблюдаются условия:

1) для неотрицательной суммируемой на $ \R $ функции $ \Phi $ помимо выполнения
(1.2.19) имеет место включение
\begin{equation*} \tag{2.1.1}
\tau \times \int_{\eta \in \R: | \eta| \ge \tau} \Phi(\eta) d\eta \in L_1(\R_+);
\end{equation*}

2) существует неотрицательная суммируемая на $ \R $ функция
$ \tilde \Phi^\prime, $ для которой почти для всех $ x \in \R $ выполняется
неравенство
\begin{equation*} \tag{2.1.2}
\biggl| \frac{d \tilde \phi} {dx}(x)\biggr| \le \int_{(1/2) B^1} \tilde \Phi^\prime(x -u) du,
\end{equation*}
а функция
\begin{equation*} \tag{2.1.3}
t \times \int_{\eta \in \R: | \eta| \ge t}
\tilde \Phi^\prime(\eta) d\eta \in L_1(\R_+),
\end{equation*}
и $ 1 < p < \infty. $ Тогда существует константа
$ c_1(\phi, \tilde \phi, \Phi, \tilde \Phi^\prime, p) >0 $ такая, что при любом
$ k \in \Z_+ $ для любого набора чисел $ \sigma = \{ \sigma_\kappa \in
\{-1, 1\}: \kappa =0, \ldots, k \}, $
для $ f \in L_p(\R) $ справедливо неравенство
\begin{equation*} \tag{2.1.4}
\biggl\| \sum_{\kappa =0}^k \sigma_\kappa \cdot (\mathcal E_\kappa f)
\biggr\|_{L_p(\R)} \le c_1 \| f \|_{L_p(\R)},
\end{equation*}
и. в частности,
\begin{equation*} \tag{2.1.5}
\| E_k f\|_{L_p(\R)} \le c_1 \| f \|_{L_p(\R)}.
\end{equation*}

Например, если $ \phi_N^D $ -- масштабирующая функция Добеши
порядка $ N $ (см. п. 7.3 из [6] или [7]), то при достаточно
больших $ N \in \N $ пара функций $ \phi = \phi_N^D, \tilde \phi = \phi_N^D $
удовлетворяет условиям леммы 2.1.1.

Отметим также, что доказательство леммы 2.1.1 проводится по схеме, использованной
в [14] при доказательстве теоремы 1 из гл. II (см. также доказательство
аналогичного утверждения в [5]).

Доказательство.

Сначала установим справедливость (2.1.4) при $ 1 < p \le 2. $
Неравенство (2.1.4) достаточно получить для функций $ f, $
принимающих вещественные значения.
В самом деле, если (2.1.4) верно для вещественнозначных функий $ f, $
то для комплекснозначной функции $ f $  имеем
\begin{multline*}
\biggl\| \sum_{\kappa =0}^k \sigma_\kappa \cdot (\mathcal E_\kappa f)
\biggr\|_{L_p(\R)} =
\biggl\| \sum_{\kappa =0}^k \sigma_\kappa \cdot (\mathcal E_\kappa ( \Re f
+i \Im f)) \biggr\|_{L_p(\R)} \\ =
\biggl\| \sum_{\kappa =0}^k \sigma_\kappa \cdot (\mathcal E_\kappa ( \Re f))
+i \sum_{\kappa =0}^k \sigma_\kappa \cdot (\mathcal E_\kappa
( \Im f)) \biggr\|_{L_p(\R)}  \\ \le
\biggl\| \sum_{\kappa =0}^k \sigma_\kappa \cdot (\mathcal E_\kappa
( \Re f)) \biggr\|_{L_p(\R)}
+\biggl\| \sum_{\kappa =0}^k \sigma_\kappa \cdot (\mathcal E_\kappa
( \Im f)) \biggr\|_{L_p(\R)} \\ \le
c_1 \| \Re f \|_{L_p(\R)} +c_1 \| \Im f \|_{L_p(\R)}
\le 2 c_1 \| f \|_{L_p(\R)}.
\end{multline*}

Определим при $ k \in \Z_+, \sigma = \{ \sigma_\kappa \in \{-1, 1\}:
\kappa =0, \ldots, k \} $ отображение
$ T = T_{k,\sigma}: L_1(\R) +L_2(\R) \mapsto L(\R), $
(где $ L_1(\R), L_2(\R), L(\R) $ -- пространства
вещественнозначных функций), полагая для
$ f = f_1 +f_2: f_1 \in L_1(\R), f_2 \in L_2(\R) $ значение
$$
( T f)(x) = \biggl| \sum_{\kappa =0}^k \sigma_\kappa (\mathcal E_\kappa^1 ( f_1 ))(x) +
\sum_{\kappa =0}^k \sigma_\kappa (\mathcal E_\kappa^2 ( f_2 ))(x)\biggr|
\text{ почти для всех } x \in \R.
$$

Проверим корректность определения отображения $ T. $ Если $ f = f_1 +f_2 =
F_1 +F_2, $ где $ f_1, F_1 \in L_1(\R), f_2, F_2 \in L_2(\R), $ то
$$
f_1 -F_1 = F_2 -f_2 \in L_1(\R) \cap L_2(\R)
$$
и при $ \kappa =0, \ldots, k $ ввиду замечания после формулировки
предложения 1.2.2 справедливо равенство
$$
\mathcal E_\kappa^1 (f_1 -F_1) = \mathcal E_\kappa^2 (F_2 -f_2),
$$
и, значит,
$$
\sum_{\kappa =0}^k \sigma_\kappa \mathcal E_\kappa^1 (f_1 -F_1) =
\sum_{\kappa =0}^k \sigma_\kappa \mathcal E_\kappa^2 (F_2 -f_2),
$$
или
$$
\sum_{\kappa =0}^k \sigma_\kappa \mathcal E_\kappa^1 (f_1) -
\sum_{\kappa =0}^k \sigma_\kappa \mathcal E_\kappa^1 (F_1) =
\sum_{\kappa =0}^k \sigma_\kappa \mathcal E_\kappa^2 (F_2) -
\sum_{\kappa =0}^k \sigma_\kappa \mathcal E_\kappa^2 (f_2),
$$
поэтому
$$
\sum_{\kappa =0}^k \sigma_\kappa \mathcal E_\kappa^1 (f_1) +
\sum_{\kappa =0}^k \sigma_\kappa \mathcal E_\kappa^2 (f_2) =
\sum_{\kappa =0}^k \sigma_\kappa \mathcal E_\kappa^1 (F_1) +
\sum_{\kappa =0}^k \sigma_\kappa \mathcal E_\kappa^2 (F_2),
$$
и
$$
\biggl| \sum_{\kappa =0}^k \sigma_\kappa \mathcal E_\kappa^1 (f_1) +
\sum_{\kappa =0}^k \sigma_\kappa \mathcal E_\kappa^2 (f_2) \biggr| =
\biggl| \sum_{\kappa =0}^k \sigma_\kappa \mathcal E_\kappa^1 (F_1) +
\sum_{\kappa =0}^k \sigma_\kappa \mathcal E_\kappa^2 (F_2)\biggr|.
$$
Отметим, что с помощью (1.2.22) легко проверить, что если при $ 1 < p < 2 $
функция $ f \in L_p(\R) $ представляется в виде $ f = f_1 +f_2, $
где $ f_1 \in L_1(\R), f_2 \in L_2(\R), $ то при $ \kappa \in \Z_+ $ справедлво
равенство
$ E_\kappa^p f = E_\kappa^1 f_1 +E_\kappa^2 f_2. $ Имея в виду это обстоятельство,
не сложно убедиться в том, что при $ 1 < p < 2 $ для $ f \in L_p(\R) $ имеет место равенство
$$
( T f)(x) = \biggl| \sum_{\kappa =0}^k \sigma_\kappa (\mathcal E_\kappa^p f)(x)\biggr|
\text{ почти для всех } x \in \R.
$$

Теперь в рассматриваемой ситуации проверим соблюдение условий теоремы 1.4.3.
При $ k \in \Z_+, \sigma = \{ \sigma_\kappa \in \{-1, 1\}:
\kappa =0, \ldots, k \} $
для $ f = f_1 +f_2 \in (L_1(\R) +L_2(\R)), g = g_1 +g_2 \in (L_1(\R) +L_2(\R)) $
почти для всех $ x \in \R $ имеем
\begin{multline*} \tag{2.1.6}
| (T (f +g))(x) | =
| (T (f_1 +f_2 +g_1 +g_2))(x) | =
| (T (f_1 +g_1 +f_2 +g_2))(x) | \\=
\biggl| \sum_{\kappa =0}^k \sigma_\kappa (\mathcal E_\kappa^1 (f_1 +g_1))(x) +
\sum_{\kappa =0}^k \sigma_\kappa (\mathcal E_\kappa^2 (f_2 +g_2))(x) \biggr| \\=
\biggl| \sum_{\kappa =0}^k \sigma_\kappa (\mathcal E_\kappa^1 (f_1))(x) +
\sum_{\kappa =0}^k \sigma_\kappa (\mathcal E_\kappa^1 (g_1))(x) +
\sum_{\kappa =0}^k \sigma_\kappa (\mathcal E_\kappa^2 (f_2))(x) +
\sum_{\kappa =0}^k \sigma_\kappa (\mathcal E_\kappa^2 (g_2))(x) \biggr| \\=
\biggl| \sum_{\kappa =0}^k \sigma_\kappa (\mathcal E_\kappa^1 (f_1))(x) +
\sum_{\kappa =0}^k \sigma_\kappa (\mathcal E_\kappa^2 (f_2))(x) +
\sum_{\kappa =0}^k \sigma_\kappa (\mathcal E_\kappa^1 (g_1))(x) +
\sum_{\kappa =0}^k \sigma_\kappa (\mathcal E_\kappa^2 (g_2))(x) \biggr| \\ \le
\biggl| \sum_{\kappa =0}^k \sigma_\kappa (\mathcal E_\kappa^1 (f_1))(x) +
\sum_{\kappa =0}^k \sigma_\kappa (\mathcal E_\kappa^2 (f_2))(x) \biggr| +
\biggl| \sum_{\kappa =0}^k \sigma_\kappa (\mathcal E_\kappa^1 (g_1))(x) +
\sum_{\kappa =0}^k \sigma_\kappa (\mathcal E_\kappa^2 (g_2))(x) \biggr| \\=
| (T f)(x) | +| (T g)(x) |,
\end{multline*}
т.е. выполняется (1.4.8).

Далее, покажем, что при $ k \in \Z_+, \sigma = \{ \sigma_\kappa \in
\{-1, 1\}: \kappa =0, \ldots, k \}, $ для $ f \in L_2(\R) $ и $ \alpha >0 $
соблюдается неравенство
\begin{equation*} \tag{2.1.7}
\mes \{ x \in \R: | (T f)(x) | > \alpha \} \le
((1 / \alpha) \| f \|_{L_2(\R)})^2.
\end{equation*}

В самом деле, при $ k \in \Z_+, \sigma = \{ \sigma_\kappa \in \{-1, 1\}:
\kappa =0, \ldots, k \}, $ для $ f \in L_2(\R) $ и $ \alpha >0 $ в виду
(1.2.26), (1.2.32), (1.2.34) имеем
\begin{multline*}
\int_\R | (T f)(x)|^2 dx = \int_\R \biggl| \sum_{\kappa =0}^k \sigma_\kappa
(\mathcal E_\kappa^2(f))(x)\biggr|^2 dx \\ =
\int_\R \biggl| \sum_{\kappa =0}^k \sigma_\kappa
(\mathcal U_\kappa(f))(x)\biggr|^2 dx =
\int_\R \biggl(\sum_{\kappa =0}^k
\sigma_\kappa (\mathcal U_\kappa( f ))(x)\biggr)
\overline {\biggl(\sum_{\kappa^\prime =0}^k
\sigma_{\kappa^\prime} (\mathcal U_{\kappa^\prime}( f))(x)\biggr)} dx  \\ =
\int_\R \sum_{\kappa =0}^k \sum_{\kappa^\prime =0}^k
\sigma_\kappa (\mathcal U_\kappa( f ))(x)
\sigma_{\kappa^\prime} \overline {(\mathcal U_{\kappa^\prime}( f))(x)} dx  \\ =
\sum_{\kappa =0}^k \sum_{\kappa^\prime =0}^k
\sigma_\kappa \sigma_{\kappa^\prime}
\int_\R (\mathcal U_\kappa( f))(x)
\overline {(\mathcal U_{\kappa^\prime} ( f))(x)} dx \\
= \sum_{\kappa =0}^k \int_\R
| (\mathcal U_\kappa( f ))(x)|^2 dx
= \sum_{ \kappa =0}^k \| \mathcal U_\kappa f \|_{L_2(\R)}^2
\le \sum_{\kappa \in \Z_+} \| \mathcal U_\kappa f \|_{L_2(\R)}^2
= \| f \|_{L_2(\R)}^2,
\end{multline*}
откуда, как обычно, получаем
\begin{multline*}
\alpha^2 \mes \{ x \in \R: | (T f)(x) | > \alpha \} =
\int_{ \{x \in \R: | (T f)(x) | > \alpha \} } \alpha^2 dx \\
\le \int_{ \{x \in \R: | (T f)(x) | > \alpha \} } | (T f )(x)|^2 dx \le
\int_\R | (T f)(x)|^2 dx \le \| f \|_{L_2(\R)}^2,
\end{multline*}
и, значит, верно (2.1.7).

Теперь установим, что существует константа $ C_0 >0 $ такая,
что при $ k \in \Z_+, \sigma = \{ \sigma_\kappa \in \{-1, 1\}:
\kappa =0, \ldots, k \} $ для $ f \in L_1(\R) $ и $ \alpha >0 $
соблюдается неравенство
\begin{equation*} \tag{2.1.8}
\mes \{ x \in \R: | (T f)(x) | > \alpha \} \le
(C_0 / \alpha) \| f \|_{L_1(\R)}.
\end{equation*}

Пусть $ k \in \Z_+, \sigma = \{ \sigma_\kappa \in \{-1, 1\}:
\kappa =0, \ldots, k \}, f \in L_1(\R) $ и $ \alpha >0. $
Для функции $ f $ и числа $ \alpha $ построим замкнутое множество
$ F, $ множество $ W = \R \setminus F $ и семейство интервалов
$ \{ Q_r, r \in \N \}, $ для которых соблюдаются условия
(1.4.2) -- (1.4.7) при $ d =1. $

Определим функции $ g \in L_1(\R) \cap
L_2(\R) $ и $ h \in L_1(\R), $ полагая
$$
g(x) = f(x) \chi_F(x) +\sum_{r=1}^\infty (1/ \mes Q_r) \biggl(\int_{Q_r} f(y) dy\biggr)
\chi_{Q_r}(x), x \in \R,
$$
и $ h = f -g. $

Из (1.4.5), (1.4.4) с учетом того, что $ \mes (W \setminus (\cup_{
r \in \N } Q_r)) =0 $ (ибо $ (W \setminus (\cup_{ r \in \N } Q_r))
\subset \cup_{ r \in \N }( \overline Q_r \setminus Q_r)), $
вытекает, что почти для всех $ x \in \R $ имеет место равенство
$ \chi_W (x) = \sum_{r =1}^\infty \chi_{Q_r} (x). $ Поэтому почти
для всех $ x \in \R $ получаем
\begin{multline*}
h(x) = f(x) -g(x) = f(x) \chi_F (x) +f(x) \chi_W (x) -g(x) \\
=f(x) \chi_F (x) +f(x) \biggl(\sum_{r =1}^\infty \chi_{Q_r} (x)\biggr) -g(x)
\\ = f(x) (\sum_{r =1}^\infty \chi_{Q_r} (x)) -\sum_{r=1}^\infty (1/
\mes Q_r) \biggl(\int_{Q_r} f(y) dy\biggr) \chi_{Q_r}(x)\\
 = \sum_{r=1}^\infty
\biggl( f(x) -(1/ \mes Q_r) \int_{Q_r} f(y) dy\biggr) \chi_{Q_r}(x) =
\sum_{r=1}^\infty h_r (x),
\end{multline*}
где $ h_r (x) = \biggl( f(x)
-(1/ \mes Q_r) \int_{Q_r} f(y) dy\biggr) \chi_{Q_r}(x). $

 Учитывая (1.4.4), (1.4.5), на основании (1.4.2) и (1.4.7) заключаем,
что почти для всех $ x \in \R $ выполняется неравенство
$ | g(x) | \le c_2 \alpha, $ из которого вытекает оценка
\begin{multline*} \tag{2.1.9}
\| g \|_{L_2(\R)}^2 = \int_\R | g(x) |^2 dx \le \int_\R c_2 \alpha
| g(x) | dx  \\ =
c_2 \alpha \int_\R \biggl| f(x) \chi_F(x) +\sum_{r=1}^\infty
(1/ \mes Q_r) \biggl(\int_{Q_r} f(y) dy\biggr) \chi_{Q_r}(x) \biggr| dx \\ \le
c_2 \alpha \int_\R \biggl| f(x) | \chi_F(x) +\sum_{r=1}^\infty
(1/ \mes Q_r) \biggr| \int_{Q_r} f(y) dy | \chi_{Q_r}(x)  dx \\ =
c_2 \alpha \biggl(\int_\R | f(x) | \chi_F(x) dx +\sum_{r=1}^\infty \int_\R
(1/ \mes Q_r) \biggl| \int_{Q_r} f(y) dy \biggr| \chi_{Q_r}(x)  dx \biggr) \\ =
c_2 \alpha \biggl(\int_F | f(x) | dx +\sum_{r=1}^\infty
(1/ \mes Q_r) \biggl| \int_{Q_r} f(y) dy \biggr| \int_\R \chi_{Q_r}(x) dx \biggr)\\ \le
c_2 \alpha \biggl(\int_F | f(x) | dx +\sum_{r=1}^\infty
(1/ \mes Q_r) \biggl(\int_{Q_r} | f(y) | dy\biggr) \mes Q_r \biggr) \\ =
c_2 \alpha \biggl(\int_F | f(x) | dx +\sum_{r=1}^\infty
\int_{Q_r} | f(x) | dx \biggr) \\ =
c_2 \alpha \int_{F \cup (\cup_{r=1}^\infty Q_r)} | f(x) | dx  =
c_2 \alpha \int_\R | f(x) | dx = c_2 \alpha \| f \|_{L_1(\R)}.
\end{multline*}

Для получения (2.1.8), фиксируя множество $ A \subset \R:
\mes A =0 $ и для $ x \in \R \setminus A $ в виду (2.1.6)
имеет место неравенство
$$
| (T f)(x) | = | (T (g +h))(x) | \le | (T g)(x) | +| (T h)(x) |,
$$
видим, что
\begin{multline*}
(\{ x \in \R: | (T f)(x) | > \alpha \} \setminus A) \subset \{ x
\in \R: | (T g)(x)| +| (T h)(x) | > \alpha \} \\
\subset \{ x \in \R: | (T g)(x) | > \alpha /2 \} \cup \{ x \in
\R: | (T h)(x) |
> \alpha /2 \},
\end{multline*}
и, значит,
\begin{multline*} \tag{2.1.10}
\mes \{ x \in \R: | (T f)(x) | > \alpha \} = \mes (\{ x \in \R:
| (T f)(x) | > \alpha \} \setminus A ) \\
\le \mes \{ x \in \R: | (T g)(x) | > \alpha /2 \} + \mes \{ x
\in \R: | (T h)(x) | > \alpha /2 \}.
\end{multline*}

Из (2.1.7) и (2.1.9) выводим
\begin{multline*} \tag{2.1.11}
\mes \{ x \in \R: | (T g)(x) | > \alpha /2 \} \le ((2 / \alpha)
\| g \|_{L_2(\R)})^2 = (2 / \alpha)^2 \| g \|_{L_2(\R)}^2 \\
\le (2 / \alpha)^2 c_2 \alpha \| f \|_{L_1(\R)} = (c_3 / \alpha)
\| f \|_{L_1(\R)}.
\end{multline*}

Для оценки второго слагаемого в правой части (2.1.10) имеем
\begin{multline*}
\{ x \in \R: | (T h)(x) | > \alpha /2 \} \\
= \{ x \in F: | (T
h)(x) | > \alpha /2 \} \cup \{ x \in W: | (T h)(x) | > \alpha /2
\},
\end{multline*}
а, следовательно,
\begin{multline*} \tag{2.1.12}
\mes \{ x \in \R: | (T h)(x) | > \alpha /2 \} \\
= \mes \{ x \in
F: | (T h)(x) | > \alpha /2 \} + \mes \{ x \in W: | (T h)(x) |
> \alpha /2 \}.
\end{multline*}

Второе слагаемое в правой части (2.1.12) в силу (1.4.3) удовлетворяет
неравенству
\begin{equation*} \tag{2.1.13}
\mes \{ x \in W: | (T h)(x) | > \alpha /2 \} \le
\mes W \le (c_4 / \alpha) \int_\R |f| dx = (c_4 / \alpha) \| f \|_{L_1(\R)}.
\end{equation*}

С другой стороны
\begin{equation*} \tag{2.1.14}
\mes \{ x \in F: | (T h)(x) | > \alpha /2 \} \\
\le (2 / \alpha) \| T h \|_{L_1(F)}.
\end{equation*}

Для проведения оценки правой части (2.1.14) определим при $ m \in \N $
функцию $ h_m^\prime $ равенством
$$
h_m^\prime = h -\sum_{r=1}^m h_r
$$
и заметим, что вследствие (2.1.6) и с учетом предложения 1.2.2 при
$ m \in \N $ почти для всех $ x \in F $ справедливо неравенство
\begin{multline*} \tag{2.1.15}
| (T h)(x) | = | (T (\sum_{r=1}^m h_r +h_m^\prime))(x) | \le
\sum_{r=1}^m | (T h_r)(x) | +| (T h_m^\prime)(x) |\\ \le
\sum_{r=1}^m | (T h_r)(x) | +\| (T h_m^\prime)(\cdot) \|_{L_\infty(\R)}.
\end{multline*}

При $ r \in \N $ оценим сверху значения $ | (T h_r)(x) | $ для $ x \in F. $

Для этого при $ r \in \N $ представим множество $ F $ в виде объединения
$$
F = \cup_{\iota \in \Z} F_{r, \iota},
$$
попарно непересекающихся измеримых множеств
\begin{equation*} \tag{2.1.16}
F_{r, \iota} = \{x \in F: 2^{-\iota} \le \rho(x, Q_r) < 2^{-\iota +1} \},
\iota \in \Z.
\end{equation*}

Принимая во внимание (1.2.22), при $ r \in \N, \iota \in \Z $ почти для всех
$ x \in F_{r, \iota} $ справедливо неравенство
\begin{multline*} \tag{2.1.17}
| (T h_r)(x) | = \biggl| \sum_{\kappa =0}^k \sigma_\kappa
(\mathcal E_\kappa ( H_r))(x)\biggr|  =
\biggl| \sum_{\kappa =0}^k \sigma_\kappa
(( E_\kappa -E_{\kappa -1})( H_r))(x)\biggr| \\ \le
\sum_{\kappa =0}^k \biggl| \sigma_\kappa
(( E_\kappa(h_r))(x) -(E_{\kappa -1}( H_r))(x)) \biggr|  \le
\sum_{\kappa =0}^k ( | ( E_\kappa(h_r))(x) | +
| (E_{\kappa -1}( H_r))(x) |) \\ \le
2 \sum_{\kappa =0}^k | ( E_\kappa(h_r))(x) | =
2 \sum_{\kappa =0}^{k} \biggl| \sum_{ \nu_\kappa \in \Z} 2^\kappa
\times \biggl(\int_\R h_r(y)
\overline { \tilde \phi(2^\kappa y -\nu_\kappa) } dy \biggr) \times
\phi(2^\kappa x -\nu_\kappa) \biggr| \\
= 2 \sum_{\kappa =0}^{k} \biggl| \lim_{n \to \infty}
\sum_{ \nu_\kappa =-n}^n 2^\kappa
\times \biggl(\int_{Q_r} h_r(y)
\overline { \tilde \phi(2^\kappa y -\nu_\kappa) } dy \biggr) \times
\phi(2^\kappa x -\nu_\kappa) \biggr| \\
= 2 \sum_{\kappa =0}^{k} \lim_{n \to \infty}
\biggl| \sum_{ \nu_\kappa =-n}^n 2^\kappa
\times \biggl(\int_{Q_r} h_r(y)
\overline { \tilde \phi(2^\kappa y -\nu_\kappa) } dy \biggr) \times
\phi(2^\kappa x -\nu_\kappa) \biggr| \\
= 2 \lim_{n \to \infty} \sum_{\kappa =0}^{k}
\biggl| \sum_{ \nu_\kappa =-n}^n 2^\kappa
\times \biggl(\int_{Q_r} h_r(y)
\overline { \tilde \phi(2^\kappa y -\nu_\kappa) } dy \biggr) \times
\phi(2^\kappa x -\nu_\kappa) \biggr| \\
\le 2 \lim_{n \to \infty} \sum_{\kappa =0}^{k}
\sum_{ \nu_\kappa =-n}^n 2^\kappa
\times \biggl| \int_{Q_r} h_r(y)
\overline { \tilde \phi(2^\kappa y -\nu_\kappa) } dy \biggr|
\times | \phi(2^\kappa x -\nu_\kappa)|.
\end{multline*}

Проводя оценку правой части (2.1.17), при $ r \in \N, \iota \in \Z $ и
$ \Nu \subset \Z $ -- любое конечное множество, почти для всех
$ x \in F_{r, \iota} $ имеем
\begin{multline*} \tag{2.1.18}
\sum_{\kappa =0}^{k}
\sum_{ \nu_\kappa \in \Nu} 2^\kappa
\times \biggl| \int_{Q_r} h_r(y)
\overline { \tilde \phi(2^\kappa y -\nu_\kappa) } dy \biggr|
\times | \phi(2^\kappa x -\nu_\kappa)| \\
= \sum_{\kappa =0, \ldots, k: \kappa \le \iota} 2^\kappa \times
\sum_{ \nu_\kappa \in \Nu}
\biggl| \int_{Q_r} h_r(y)
\overline { \tilde \phi(2^\kappa y -\nu_\kappa) } dy \biggr|
\times | \phi(2^\kappa x -\nu_\kappa)| \\
+\sum_{\kappa =0, \ldots, k: \kappa > \iota} 2^\kappa \times
\sum_{ \nu_\kappa \in \Nu}
\biggl| \int_{Q_r} h_r(y)
\overline { \tilde \phi(2^\kappa y -\nu_\kappa) } dy \biggr|
\times | \phi(2^\kappa x -\nu_\kappa)|.
\end{multline*}

Для оценки правой части (2.1.18) заметим, что при $ r \in \N $
выполняется равенство
\begin{multline*} \tag{2.1.19}
\int_{Q_r} h_r (x) dx = \int_{Q_r} ( f(x) -(1/ \mes Q_r)
\int_{Q_r} f(y) dy) \chi_{Q_r}(x) dx\\ = \int_{Q_r} ( f(x) -(1/
\mes Q_r) \int_{Q_r} f(y) dy) dx \\
= \int_{Q_r} f(x) dx -(1/ \mes Q_r) \biggl(\int_{Q_r} f(y) dy\biggr)
\int_{Q_r} dx \\
= \int_{Q_r} f(x) dx -\int_{Q_r} f(y) dy =0.
\end{multline*}

Фиксируя для каждого $ r \in \N $ точку $ y_r \in Q_r, $
при $ r \in \N, \kappa =0, \ldots, k, \nu_\kappa \in \Z $
с учетом (2.1.19) выводим
\begin{multline*} \tag{2.1.20}
\biggl| \int_{ Q_r} h_r(y) \overline { \tilde \phi( 2^\kappa y -\nu_\kappa) } dy \biggr| \\ =
\biggl| \int_{ Q_r} h_r(y) (\overline { \tilde \phi( 2^\kappa y -\nu_\kappa) } -
\overline { \tilde \phi( 2^\kappa y_r -\nu_\kappa) }) +h_r(y)
\overline { \tilde \phi( 2^\kappa y_r -\nu_\kappa) } dy \biggr| \\ =
\biggl| \int_{ Q_r} h_r(y) (\overline { \tilde \phi( 2^\kappa y -\nu_\kappa) } -
\overline { \tilde \phi( 2^\kappa y_r -\nu_\kappa) }) dy +
\overline { \tilde \phi( 2^\kappa y_r -\nu_\kappa) } \times \int_{ Q_r} h_r(y) dy \biggr| \\ =
\biggl| \int_{ Q_r} h_r(y) (\overline { \tilde \phi( 2^\kappa y -\nu_\kappa) } -
\overline { \tilde \phi( 2^\kappa y_r -\nu_\kappa) }) dy \biggr| \\ \le
\int_{ Q_r} | h_r(y)| \times | \overline { \tilde \phi( 2^\kappa y -\nu_\kappa) } -
\overline { \tilde \phi( 2^\kappa y_r -\nu_\kappa) } | dy.
\end{multline*}

Далее, при $ r \in \N, \kappa =0, \ldots, k, \nu_\kappa \in \Z, $
пользуясь тем, что вследствие (2.1.2) функция $ \frac{d \tilde \phi} {dx} $
ограничена на $ \R $ (см. вывод (1.2.20) при $ p = \infty $), для $ y \in Q_r $
получаем
\begin{multline*} \tag{2.1.21}
| \overline { \tilde \phi( 2^\kappa y -\nu_\kappa) } -
\overline { \tilde \phi( 2^\kappa y_r -\nu_\kappa) } | =
| \tilde \phi( 2^\kappa y -\nu_\kappa) -\tilde \phi( 2^\kappa y_r -\nu_\kappa) | \\ =
\biggl| \int_0^1 \frac{d} {dt} ( \tilde \phi( 2^\kappa (y_r +t (y -y_r)) -
\nu_\kappa) ) dt \biggr| \\ =
\biggl| \int_0^1 2^\kappa (y -y_r) \frac{d \tilde \phi} {du} ( 2^\kappa (y_r +t (y -y_r))
-\nu_\kappa) dt \biggl| \\ =
2^\kappa | y -y_r | \times \biggl| \int_0^1 \frac{d \tilde \phi} {du} ( 2^\kappa (y_r +
t (y -y_r)) -\nu_\kappa) dt \biggr| \\ \le
2^\kappa | y -y_r | \times \int_0^1 \biggl| \frac{d \tilde \phi} {du} ( 2^\kappa (y_r +
t (y -y_r)) -\nu_\kappa) \biggr| dt \\ \le
2^\kappa | y -y_r | \times \biggl\| \frac{d \tilde \phi} {du} \biggr\|_{L_\infty(\R)} \le
2^\kappa ( \diam Q_r ) \biggl\| \frac{d \tilde \phi} {du} \biggr\|_{L_\infty(\R)} =
c_5(\tilde \phi) 2^\kappa ( \diam Q_r ).
\end{multline*}

Соединяя (2.1.20) и (2.1.21), находим, что при $ r \in \N,
\kappa =0, \ldots, k, \nu_\kappa \in \Z $ выполняется неравенство
\begin{multline*} \tag{2.1.22}
\biggl| \int_{ Q_r} h_r(y) \overline { \tilde \phi( 2^\kappa y -\nu_\kappa) } dy \biggr| \\ \le
\int_{ Q_r} c_5 2^\kappa ( \diam Q_r ) | h_r(y) | dy = c_5
2^\kappa ( \diam Q_r ) \int_{ Q_r} | h_r(y) | dy.
\end{multline*}

Используя неравенство (2.1.22) при оценке первой суммы в правой части
(2.1.18), получаем, что
при $ r \in \N, \iota \in \Z, \Nu \subset \Z:
\card \Nu < \infty, $ почти для всех $ x \in F_{r, \iota} $ соблюдается
неравенство
\begin{multline*} \tag{2.1.23}
\sum_{\kappa =0, \ldots, k: \kappa \le \iota} 2^\kappa \times
\sum_{ \nu_\kappa \in \Nu}
\biggl| \int_{Q_r} h_r(y)
\overline { \tilde \phi(2^\kappa y -\nu_\kappa) } dy \biggr|
\times | \phi(2^\kappa x -\nu_\kappa)| \\
\le \sum_{\kappa =0, \ldots, k: \kappa \le \iota} 2^\kappa \times
\sum_{ \nu_\kappa \in \Nu}
c_5 2^\kappa ( \diam Q_r ) \biggl(\int_{ Q_r} | h_r(y) | dy\biggr) \\
\times | \phi(2^\kappa x -\nu_\kappa)| \\
= c_5 ( \diam Q_r ) \biggl(\int_{ Q_r} | h_r(y) | dy\biggr) \times
\sum_{\kappa =0, \ldots, k: \kappa \le \iota} 2^{2\kappa} \times
\sum_{ \nu_\kappa \in \Nu} | \phi(2^\kappa x -\nu_\kappa)|.
\end{multline*}
Из (1.2.24) и аддитивности интеграла как функции множеств вытекает, что при
$ \kappa =0, \ldots, k $ для $ \Nu \subset \Z: \card \Nu < \infty, $ почти 
для всех $ x \in \R $ справедливо неравенство
\begin{multline*} \tag{2.1.24}
\sum_{ \nu_\kappa \in \Nu} | \phi(2^\kappa x -\nu_\kappa)| \le
\sum_{ \nu_\kappa \in \Nu}
2^\kappa \int_{2^{-\kappa}(\nu_\kappa +(1/2) B^1)}
\Phi(2^\kappa (x -u)) du \\
= 2^\kappa \int_{\cup_{ \nu_\kappa \in \Nu}
2^{-\kappa} (\nu_\kappa +(1/2) B^1)}
\Phi(2^\kappa (x -u)) du \\
\le 2^\kappa \int_{\R} \Phi(2^\kappa (x -u)) du
= 2^\kappa \int_{\R} \Phi(t) 2^{-\kappa} dt
= \|\Phi\|_{L_1(\R)}.
\end{multline*}

Подставляя (2.1.24) в (2.1.23), заключаем, что
при $ r \in \N, \iota \in \Z, \Nu \subset \Z: \card \Nu < \infty, $
почти для всех $ x \in F_{r, \iota} $ верно неравенство
\begin{multline*} \tag{2.1.25}
\sum_{\kappa =0, \ldots, k: \kappa \le \iota} 2^\kappa \times
\sum_{ \nu_\kappa \in \Nu}
\biggl| \int_{Q_r} h_r(y)
\overline { \tilde \phi(2^\kappa y -\nu_\kappa) } dy \biggr|
\times | \phi(2^\kappa x -\nu_\kappa)| \\
\le c_5 ( \diam Q_r ) \biggl(\int_{ Q_r} | h_r(y) | dy\biggr) \times
\sum_{\kappa =0, \ldots, k: \kappa \le \iota} 2^{2\kappa} \|\Phi\|_{L_1(\R)} \\ =
c_6 ( \diam Q_r ) \biggl(\int_{ Q_r} | h_r(y) | dy\biggr) \times
\sum_{\kappa =0, \ldots, k: \kappa \le \iota} 2^{2\kappa} \\ \le
c_{7} ( \diam Q_r ) \biggl(\int_{ Q_r} | h_r(y) | dy\biggr) \times 2^{2 \iota}.
\end{multline*}

Далее, заметим, что при $ r \in \N $ в силу (1.4.7) справедливо неравенство
\begin{multline*} \tag{2.1.26}
\int_{Q_r} | h_r (y) | dy = \int_{Q_r} \biggl| ( f(y) -(1/ \mes Q_r)
\int_{Q_r} f(z) dz) \chi_{Q_r}(y) \biggr| dy \\ =
\int_{Q_r} \biggl| f(y) -(1/ \mes Q_r)
\int_{Q_r} f(z) dz \biggr| dy \\ \le
\int_{Q_r} | f(y) | dy +(1/ \mes Q_r)
\biggl| \int_{Q_r} f(z) dz \biggr| \int_{Q_r} dy \\ \le
\int_{Q_r} | f(y) | dy +\int_{Q_r} | f(z) | dz = 2 \int_{Q_r} | f(y) | dy
\le c_{8} \alpha \mes Q_r.
\end{multline*}

Кроме того, благодаря (1.4.6), при $ r \in \N, $ для $ y \in Q_r $
справедливо неравенство
\begin{equation*} \tag{2.1.27}
\diam Q_r \le c_{9} \rho(y, F),
\end{equation*}
а также при $ r \in \N $ для $ x \in F $ и $ y \in Q_r $
соблюдается соотношение
\begin{multline*} \tag{2.1.28}
\rho(x, Q_r) = \inf_{\eta \in Q_r} | x -\eta| \le | x -y| \\
\le (\inf_{\eta \in Q_r}| x -\eta|) +(\diam Q_r) = \rho(x, Q_r) +(\diam Q_r) \\
\le \rho(x, Q_r)
+c_{10} \inf_{ \eta \in Q_r} \rho(\eta, F)  = \rho(x, Q_r)
+c_{10} \inf_{ \eta \in Q_r} \inf_{\xi \in F} | \eta -\xi| \\ =
\rho(x, Q_r) +c_{10} \inf_{\xi \in F} \inf_{ \eta \in Q_r} | \eta -\xi| =
\rho(x, Q_r) +c_{10} \inf_{\xi \in F} \rho(\xi, Q_r)\\ \le c_{11} \rho(x, Q_r).
\end{multline*}

Сопоставляя (2.1.28) с (2.1.16), видим, что при $ r \in \N,
\iota \in \Z $ для $ x \in F_{r, \iota}, y \in Q_r $ имеет место соотношение
\begin{equation*}
2^{-\iota} \le \rho(x, Q_r) \le | x -y| \le c_{11} \rho(x, Q_r) < c_{12} 2^{-\iota}
\end{equation*}
или
\begin{equation*} \tag{2.1.29}
c_{13} 2^\iota < | x -y|^{-1} \le 2^\iota.
\end{equation*}

Подставляя в неравенство (2.1.25) оценки (2.1.26), (2.1.27) и (2.1.29),
приходим к выводу, что при $ r \in \N, \iota \in \Z,
\Nu \subset \Z: \card \Nu < \infty, $ почти для всех $ x \in F_{r, \iota} $ выполняется неравенство
\begin{multline*} \tag{2.1.30}
\sum_{\kappa =0, \ldots, k: \kappa \le \iota} 2^\kappa \times
\sum_{ \nu_\kappa \in \Nu}
\biggl| \int_{Q_r} h_r(y) \overline { \tilde \phi(2^\kappa y -\nu_\kappa) } dy \biggr|
\times | \phi(2^\kappa x -\nu_\kappa)| \\
\le c_{14} ( \diam Q_r ) 2^{2 \iota} \alpha \mes Q_r =
c_{14} \alpha ( \diam Q_r ) 2^{2 \iota} \int_{ Q_r } dy \\
= c_{14} \alpha \int_{ Q_r } ( \diam Q_r ) 2^{2 \iota} dy \le
c_{15} \alpha \int_{ Q_r } \rho(y, F) | x -y |^{-2} dy.
\end{multline*}

Переходя к оценке второй суммы в правой части (2.1.18), разобьем ее на две
части, а, именно, при $ r \in \N, \iota \in \Z $ и
$ \Nu \subset \Z: \card \Nu < \infty, $ почти для всех $ x \in F_{r, \iota} $ 
выводим
\begin{multline*} \tag{2.1.31}
\sum_{\kappa =0, \ldots, k: \kappa > \iota} 2^\kappa \times
\sum_{ \nu_\kappa \in \Nu}
\biggl| \int_{Q_r} h_r(y)
\overline { \tilde \phi(2^\kappa y -\nu_\kappa) } dy \biggr|
\times | \phi(2^\kappa x -\nu_\kappa)| \\
= \sum_{\kappa =0, \ldots, k: \kappa > \iota} 2^\kappa \times
\biggl(\sum_{ \nu_\kappa \in \Nu: | x -2^{-\kappa} \nu_\kappa| < (1/2) \rho(x, Q_r)}
\biggl| \int_{Q_r} h_r(y)
\overline { \tilde \phi(2^\kappa y -\nu_\kappa) } dy \biggr|
\times | \phi(2^\kappa x -\nu_\kappa)| \\
+\sum_{ \nu_\kappa \in \Nu: | x -2^{-\kappa} \nu_\kappa| \ge (1/2) \rho(x, Q_r)}
\biggl| \int_{Q_r} h_r(y)
\overline { \tilde \phi(2^\kappa y -\nu_\kappa) } dy \biggr|
\times | \phi(2^\kappa x -\nu_\kappa)|\biggr) \\
= \sum_{\kappa =0, \ldots, k: \kappa > \iota} 2^\kappa \times
\sum_{ \nu_\kappa \in \Nu: | x -2^{-\kappa} \nu_\kappa| < (1/2) \rho(x, Q_r)}
\biggl| \int_{Q_r} h_r(y)
\overline { \tilde \phi(2^\kappa y -\nu_\kappa) } dy \biggr|
\times | \phi(2^\kappa x -\nu_\kappa)| \\
+\sum_{\kappa =0, \ldots, k: \kappa > \iota} 2^\kappa \times
\sum_{ \nu_\kappa \in \Nu: | x -2^{-\kappa} \nu_\kappa| \ge (1/2) \rho(x, Q_r)}
\biggl| \int_{Q_r} h_r(y)
\overline { \tilde \phi(2^\kappa y -\nu_\kappa) } dy \biggr|
\times | \phi(2^\kappa x -\nu_\kappa)|.
\end{multline*}

Оценивая первую сумму в правой части (2.1.31), при $ r \in \N, \iota \in \Z $
и $ \Nu \subset \Z: \card \Nu < \infty, $ почти для всех
$ x \in F_{r, \iota}, $ в силу (2.1.20), (2.1.21), (1.2.20) при $ p = \infty, $
имеем
\begin{multline*} \tag{2.1.32}
\sum_{\kappa =0, \ldots, k: \kappa > \iota} 2^\kappa \times
\sum_{ \nu_\kappa \in \Nu: | x -2^{-\kappa} \nu_\kappa| < (1/2) \rho(x, Q_r)}
\biggl| \int_{Q_r} h_r(y)
\overline { \tilde \phi(2^\kappa y -\nu_\kappa) } dy \biggr|
\times | \phi(2^\kappa x -\nu_\kappa)| \\
\le \sum_{\kappa =0, \ldots, k: \kappa > \iota} 2^\kappa \times
\sum_{ \nu_\kappa \in \Nu: | x -2^{-\kappa} \nu_\kappa| < (1/2) \rho(x, Q_r)}
\biggl(\int_{ Q_r} | h_r(y)| \times
2^\kappa \times | y -y_r | \\
\times \int_0^1 \biggl| \frac{d \tilde \phi} {du} ( 2^\kappa (y_r +
t (y -y_r)) -\nu_\kappa) \biggr| dt dy\biggr) \times | \phi(2^\kappa x -\nu_\kappa)| \\
\le \sum_{\kappa =0, \ldots, k: \kappa > \iota} 2^{2 \kappa} \times
\sum_{ \nu_\kappa \in \Nu: | x -2^{-\kappa} \nu_\kappa| < (1/2) \rho(x, Q_r)}
\biggl(\int_{ Q_r} | h_r(y)| \times
(\diam Q_r) \\
\times \int_0^1 \biggl| \frac{d \tilde \phi} {du} ( 2^\kappa (y_r +
t (y -y_r)) -\nu_\kappa) \biggr| dt dy\biggr) \times \| \phi \|_{L_\infty(\R)} \\
= c_{16}(\phi) (\diam Q_r) \times
\sum_{\kappa =0, \ldots, k: \kappa > \iota} 2^{2 \kappa} \times
\sum_{ \nu_\kappa \in \Nu: | x -2^{-\kappa} \nu_\kappa| < (1/2) \rho(x, Q_r)}
\biggl(\int_{ Q_r} | h_r(y)| \\ \times
\int_0^1 \biggl| \frac{d \tilde \phi} {du} ( 2^\kappa (y_r +
t (y -y_r)) -\nu_\kappa) \biggr| dt dy\biggr) \\
= c_{16} (\diam Q_r) \times
\sum_{\kappa =0, \ldots, k: \kappa > \iota} 2^{2 \kappa} \times
\biggl(\int_{ Q_r} | h_r(y)| \\ \times
\int_0^1 \sum_{ \nu_\kappa \in \Nu: | x -2^{-\kappa} \nu_\kappa| < (1/2) \rho(x, Q_r)}
\biggl| \frac{d \tilde \phi} {du} ( 2^\kappa (y_r +t (y -y_r)) -\nu_\kappa) \biggr| dt dy\biggr).
\end{multline*}

Принимая во внимание, что при $ r \in \N, \kappa \in \Z_+, \nu_\kappa \in \Z, $
для $ y \in Q_r \setminus \{ y_r \} $ почти для всех $ t \in (0, 1) $ 
вследствие (2.1.2) выполняется неравенство (ср. с (1.2.24))
\begin{equation*}
\biggl| \frac{d \tilde \phi} {du} ( 2^\kappa (y_r +t (y -y_r)) -\nu_\kappa) \biggr| \le
2^\kappa \int_{2^{-\kappa}(\nu_\kappa +(1/2) B^1)}
\tilde \Phi^\prime ( 2^\kappa (y_r +t (y -y_r) -z)) dz,
\end{equation*}
а также замечая, что при $ r \in \N, \iota \in \Z, \kappa =0, \ldots, k: \iota < \kappa, $
для $ x \in F_{r, \iota}, y \in Q_r, t \in (0, 1), $ и $ \nu_\kappa \in \Z:
| x -2^{-\kappa} \nu_\kappa| < (1/2) \rho(x, Q_r), z \in
2^{-\kappa}(\nu_\kappa +(1/2) B^1), $ с учетом (2.1.16) справедливо неравенство
\begin{multline*}
| y_r +t(y -y_r) -z | =
| y_r +t(y -y_r) -x +x -2^{-\kappa} \nu_\kappa +2^{-\kappa} \nu_\kappa -z | \\
\ge | y_r +t(y -y_r) -x | -| x -2^{-\kappa} \nu_\kappa| -| 2^{-\kappa} \nu_\kappa -z | \\ \ge
\inf_{\eta \in Q_r} | \eta -x| -(1/2) \rho(x, Q_r) -2^{-(\kappa +1)}
\ge \rho(x, Q_r) -(1/2) \rho(x, Q_r) -2^{-(\iota +2)}\\
= (1/2) \rho(x, Q_r) -(1/4) 2^{-\iota} \ge (1/2) 2^{-\iota} -(1/4) 2^{-\iota}
= 2^{-\iota -2},
\end{multline*}
заключаем, что при $ r \in \N, \iota \in \Z, \kappa =0, \ldots, k:
\iota < \kappa, \Nu \subset \Z: \card \Nu < \infty, $
для всех $ x \in F_{r, \iota}, y \in Q_r \setminus \{ y_r \} $ имеет 
место неравенство
\begin{multline*} \tag{2.1.33}
\int_0^1 \sum_{ \nu_\kappa \in \Nu: | x -2^{-\kappa} \nu_\kappa| < (1/2) \rho(x, Q_r)}
\biggl| \frac{d \tilde \phi} {du} ( 2^\kappa (y_r +t (y -y_r)) -\nu_\kappa) \biggr| dt \\ \le
\int_0^1 \sum_{ \nu_\kappa \in \Nu: | x -2^{-\kappa} \nu_\kappa| < (1/2) \rho(x, Q_r)}
2^\kappa \times \int_{2^{-\kappa}(\nu_\kappa +(1/2) B^1)}
\tilde \Phi^\prime ( 2^\kappa (y_r +t (y -y_r) -z)) dz dt \\
= 2^\kappa \times \int_0^1
\int_{\cup_{ \nu_\kappa \in \Nu: | x -2^{-\kappa} \nu_\kappa| < (1/2) \rho(x, Q_r)}
2^{-\kappa}(\nu_\kappa +(1/2) B^1)}
\tilde \Phi^\prime ( 2^\kappa (y_r +t (y -y_r) -z)) dz dt \\
\le 2^\kappa \times \int_0^1
\int_{z \in \R: | y_r +t(y -y_r) -z | \ge 2^{-\iota -2}}
\tilde \Phi^\prime ( 2^\kappa (y_r +t (y -y_r) -z)) dz dt \\
= 2^\kappa \times \int_0^1
\int_{\xi \in \R: | \xi| \ge 2^{-\iota -2}}
\tilde \Phi^\prime ( 2^\kappa \xi) d\xi dt
= 2^\kappa \times \int_{\xi \in \R: | \xi| \ge 2^{-\iota -2}}
\tilde \Phi^\prime ( 2^\kappa \xi) d\xi \\
= 2^\kappa \times \int_{\eta \in \R: | \eta| \ge 2^{\kappa -\iota -2}}
\tilde \Phi^\prime ( \eta) 2^{-\kappa} d\eta
= \int_{\eta \in \R: | \eta| \ge 2^{\kappa -\iota -2}}
\tilde \Phi^\prime ( \eta) d\eta.
\end{multline*}

Применяя (2.1.33) для оценки правой части (2.1.32) и учитывая
(2.1.3), получаем, что
при $ r \in \N, \iota \in \Z $ и
$ \Nu \subset \Z: \card \Nu < \infty, $ почти для всех $ x \in F_{r, \iota} $
справедливо неравенство
\begin{multline*} \tag{2.1.34}
\sum_{\kappa =0, \ldots, k: \kappa > \iota} 2^{2 \kappa}  \times
\biggl(\int_{ Q_r} | h_r(y)| \\ \times
\int_0^1 \sum_{ \nu_\kappa \in \Nu: | x -2^{-\kappa} \nu_\kappa| < (1/2) \rho(x, Q_r)}
\biggl| \frac{d \tilde \phi} {du} ( 2^\kappa (y_r +t (y -y_r)) -\nu_\kappa) \biggr| dt dy\biggr) \\
\le \sum_{\kappa =0, \ldots, k: \kappa > \iota} 2^{2 \kappa} \times
\int_{ Q_r} | h_r(y)| \times
\int_{\eta \in \R: | \eta| \ge 2^{\kappa -\iota -2}}
\tilde \Phi^\prime ( \eta) d\eta dy \\
= \sum_{\kappa =0, \ldots, k: \kappa > \iota} 2^{2 \kappa} \times
\biggl(\int_{ Q_r} | h_r(y)| dy\biggr) \times
\int_{\eta \in \R: | \eta| \ge 2^{\kappa -\iota -2}}
\tilde \Phi^\prime ( \eta) d\eta \\
= \biggl(\int_{ Q_r} | h_r(y)| dy\biggr) \times
\sum_{\kappa =0, \ldots, k: \kappa > \iota} 2^{2 \kappa} \times
\int_{\eta \in \R: | \eta| \ge 2^{\kappa -\iota -2}}
\tilde \Phi^\prime ( \eta) d\eta \\
= \biggl(\int_{ Q_r} | h_r(y)| dy\biggr) \times
\sum_{\kappa =0, \ldots, k: \kappa > \iota} \int_{2^\kappa +2^\kappa I} 2^{ \kappa} \times
\int_{\eta \in \R: | \eta| \ge 2^{\kappa +1} 2^{-\iota -3}}
\tilde \Phi^\prime ( \eta) d\eta ds \\
\le \biggl(\int_{ Q_r} | h_r(y)| dy\biggr) \times
\sum_{\kappa =0, \ldots, k: \kappa > \iota} \int_{2^\kappa +2^\kappa I} s \times
\int_{\eta \in \R: | \eta| \ge s 2^{-\iota -3}}
\tilde \Phi^\prime ( \eta) d\eta ds \\
= \biggl(\int_{ Q_r} | h_r(y)| dy\biggr) \times
\int_{\cup_{\kappa =0, \ldots, k: \kappa > \iota} (2^\kappa +2^\kappa I)}
s \times \int_{\eta \in \R: | \eta| \ge s 2^{-\iota -3}}
\tilde \Phi^\prime ( \eta) d\eta ds \\
\le \biggl(\int_{ Q_r} | h_r(y)| dy\biggr) \times
\int_{2^{\iota +1}}^\infty s \times
\int_{\eta \in \R: | \eta| \ge s 2^{-\iota -3}}
\tilde \Phi^\prime ( \eta) d\eta ds \\
= \biggl(\int_{ Q_r} | h_r(y)| dy\biggr) \times
\int_{1/4}^\infty 2^{\iota +3} t \times
\int_{\eta \in \R: | \eta| \ge t}
\tilde \Phi^\prime ( \eta) d\eta 2^{\iota +3} dt \\
= c_{17}(\tilde \Phi^\prime) 2^{2 \iota} \times
\biggl(\int_{ Q_r} | h_r(y)| dy\biggr).
\end{multline*}
Объединяя (2.1.34) с (2.1.32), а затем используя (2.1.26),
(2.1.27), (2.1.29), так же, как при выводе (2.1.30), находим, что
при $ r \in \N, \iota \in \Z $ и
$ \Nu \subset \Z: \card \Nu < \infty, $ почти для всех $ x \in F_{r, \iota} $
соблюдается неравенство
\begin{multline*} \tag{2.1.35}
\sum_{\kappa =0, \ldots, k: \kappa > \iota} 2^\kappa \times
\sum_{ \nu_\kappa \in \Nu: | x -2^{-\kappa} \nu_\kappa| < (1/2) \rho(x, Q_r)}
\biggl| \int_{Q_r} h_r(y)
\overline { \tilde \phi(2^\kappa y -\nu_\kappa) } dy \biggr|
\times | \phi(2^\kappa x -\nu_\kappa)| \\
\le c_{18} (\diam Q_r) 2^{2 \iota} \times
\biggl(\int_{ Q_r} | h_r(y)| dy\biggr) \le
c_{19} \alpha \int_{ Q_r } \rho(y, F) | x -y |^{-2} dy.
\end{multline*}

Перейдем к оценке второй суммы в правой части (2.1.31).
Используя (2.1.22), (1.2.24) и замечая, что при $ r \in \N, \iota \in \Z,
\kappa =0, \ldots, k: \iota < \kappa, $
для $ x \in F_{r, \iota}, \nu_\kappa \in \Z: | x -2^{-\kappa} \nu_\kappa| \ge
(1/2) \rho(x, Q_r), z \in 2^{-\kappa}(\nu_\kappa +(1/2) B^1), $
в виду (2.1.16) выполняется неравенство
\begin{multline*} \tag{2.1.36}
| x -z| = | x -2^{-\kappa} \nu_\kappa +2^{-\kappa} \nu_\kappa -z| \ge
| x -2^{-\kappa} \nu_\kappa| -| 2^{-\kappa} \nu_\kappa -z| \\ \ge
(1/2) \rho(x, Q_r) -(1/2) 2^{-\kappa} \ge (1/2) 2^{-\iota} -(1/2) 2^{-(\iota +1)}
= 2^{-\iota -2},
\end{multline*}
с учетом (2.1.1) приходим к выводу, что при $ r \in \N, \iota \in \Z, \Nu \subset \Z: \card \Nu < \infty, $
почти для всех $ x \in F_{r, \iota} $ выполняется неравенство
\begin{multline*} \tag{2.1.37}
\sum_{\kappa =0, \ldots, k: \kappa > \iota} 2^\kappa \times
\sum_{ \nu_\kappa \in \Nu: | x -2^{-\kappa} \nu_\kappa| \ge (1/2) \rho(x, Q_r)}
\biggl| \int_{Q_r} h_r(y)
\overline { \tilde \phi(2^\kappa y -\nu_\kappa) } dy \biggr|
\times | \phi(2^\kappa x -\nu_\kappa)| \\ \le
\sum_{\kappa =0, \ldots, k: \kappa > \iota} 2^\kappa \times
\sum_{ \nu_\kappa \in \Nu: | x -2^{-\kappa} \nu_\kappa| \ge (1/2) \rho(x, Q_r)}
c_5 2^\kappa ( \diam Q_r ) \biggl(\int_{ Q_r} | h_r(y) | dy\biggr)
\times | \phi(2^\kappa x -\nu_\kappa)| \\ =
c_5 ( \diam Q_r ) \biggl(\int_{ Q_r} | h_r(y) | dy\biggr)
\sum_{\kappa =0, \ldots, k: \kappa > \iota} 2^{2 \kappa} \times
\sum_{ \nu_\kappa \in \Nu: | x -2^{-\kappa} \nu_\kappa| \ge (1/2) \rho(x, Q_r)}
| \phi(2^\kappa x -\nu_\kappa)| \\ \le
c_5 ( \diam Q_r ) \biggl(\int_{ Q_r} | h_r(y) | dy\biggr)
\sum_{\kappa =0, \ldots, k: \kappa > \iota} 2^{2 \kappa} \\ \times
\sum_{ \nu_\kappa \in \Nu: | x -2^{-\kappa} \nu_\kappa| \ge (1/2) \rho(x, Q_r)}
2^\kappa \int_{2^{-\kappa}(\nu_\kappa +(1/2) B^1)}
\Phi(2^\kappa (x -z)) dz \\ =
c_{5} ( \diam Q_r ) \biggl(\int_{ Q_r} | h_r(y) | dy\biggr)
\sum_{\kappa =0, \ldots, k: \kappa > \iota} 2^{2 \kappa} \\ \times
\int_{\cup_{ \nu_\kappa \in \Nu: | x -2^{-\kappa} \nu_\kappa| \ge
(1/2) \rho(x, Q_r)} 2^{-\kappa} (\nu_\kappa +(1/2) B^1)}
2^\kappa \times \Phi(2^\kappa (x -z)) dz \\ \le
c_{5} ( \diam Q_r ) \biggl(\int_{ Q_r} | h_r(y) | dy\biggr)
\sum_{\kappa =0, \ldots, k: \kappa > \iota} 2^{2 \kappa} \times
\int_{z \in \R: | x -z| \ge 2^{-\iota -2}}
2^\kappa \times \Phi(2^\kappa (x -z)) dz \\
= c_{5} ( \diam Q_r ) \biggl(\int_{ Q_r} | h_r(y) | dy\biggr)
\sum_{\kappa =0, \ldots, k: \kappa > \iota} 2^{2 \kappa} \times
\int_{\xi \in \R: | \xi| \ge 2^{-\iota -2}}
2^\kappa \times \Phi(2^\kappa \xi ) d\xi  \\
= c_{5} ( \diam Q_r ) \biggl(\int_{ Q_r} | h_r(y) | dy\biggr)
\sum_{\kappa =0, \ldots, k: \kappa > \iota} 2^{2 \kappa} \times
\int_{\eta \in \R: | \eta| \ge 2^{\kappa -\iota -2}}
2^\kappa \times \Phi(2^\kappa 2^{-\kappa} \eta ) 2^{-\kappa} d\eta  \\
= c_{5} ( \diam Q_r ) \biggl(\int_{ Q_r} | h_r(y) | dy\biggr)
\sum_{\kappa =0, \ldots, k: \kappa > \iota} 2^{2 \kappa} \times
\int_{\eta \in \R: | \eta| \ge 2^{\kappa +1} 2^{-\iota -3}}
\Phi(\eta) d\eta \\
= c_{5} ( \diam Q_r ) \biggl(\int_{ Q_r} | h_r(y) | dy\biggr)
\sum_{\kappa =0, \ldots, k: \kappa > \iota} \int_{2^\kappa +2^\kappa I} 2^{\kappa} \times
\int_{\eta \in \R: | \eta| \ge 2^{\kappa +1} 2^{-\iota -3}}
\Phi(\eta) d\eta ds \\
\le c_{5} ( \diam Q_r ) \biggl(\int_{ Q_r} | h_r(y) | dy\biggr)
\sum_{\kappa =0, \ldots, k: \kappa > \iota} \int_{2^\kappa +2^\kappa I} s \times
\int_{\eta \in \R: | \eta| \ge s 2^{-\iota -3}}
\Phi(\eta) d\eta ds \\
= c_{5} ( \diam Q_r ) \biggl(\int_{ Q_r} | h_r(y) | dy\biggr)
\int_{\cup_{\kappa =0, \ldots, k: \kappa > \iota} (2^\kappa +2^\kappa I)}
s \times \int_{\eta \in \R: | \eta| \ge s 2^{-\iota -3}} \Phi(\eta) d\eta ds \\
\le c_{5} ( \diam Q_r ) \biggl(\int_{ Q_r} | h_r(y) | dy\biggr)
\int_{2^{\iota +1}}^\infty s \times
\int_{\eta \in \R: | \eta| \ge s 2^{-\iota -3}} \Phi(\eta) d\eta ds \\
=  c_{5} ( \diam Q_r ) \biggl(\int_{ Q_r} | h_r(y) | dy\biggr)
\int_{1/4}^\infty 2^{\iota +3} \tau \times
\int_{\eta \in \R: | \eta| \ge \tau} \Phi(\eta) d\eta 2^{\iota +3} d\tau \\
= c_{20} ( \diam Q_r ) 2^{2 \iota}
\biggl(\int_{ Q_r} | h_r(y) | dy\biggr)
\int_{1/4}^\infty \tau \times
\int_{\eta \in \R: | \eta| \ge \tau} \Phi(\eta) d\eta d\tau \\
= c_{21} ( \diam Q_r ) 2^{2 \iota} \biggl(\int_{ Q_r} | h_r(y) | dy\biggr).
\end{multline*}

Применяя (2.1.26), (2.1.27), (2.1.29) так же, как при выводе (2.1.30) или
(2.1.35), из (2.1.37) получаем, что при $ r \in \N, \iota \in \Z $ и
$ \Nu \subset \Z: \card \Nu < \infty, $ почти для всех $ x \in F_{r, \iota} $
справедливо неравенство
\begin{multline*}
\sum_{\kappa =0, \ldots, k: \kappa > \iota} 2^\kappa \times
\sum_{ \nu_\kappa \in \Nu: | x -2^{-\kappa} \nu_\kappa| \ge (1/2) \rho(x, Q_r)}
\biggl| \int_{Q_r} h_r(y)
\overline { \tilde \phi(2^\kappa y -\nu_\kappa) } dy \biggr|
\times | \phi(2^\kappa x -\nu_\kappa)| \\ \le
c_{22} \alpha \int_{ Q_r } \rho(y, F) | x -y |^{-2} dy,
\end{multline*}
которое в соединении с (2.1.35) и (2.1.31) дает оценку
\begin{multline*} \tag{2.1.38}
\sum_{\kappa =0, \ldots, k: \kappa > \iota} 2^\kappa \times
\sum_{ \nu_\kappa \in \Nu}
\biggl| \int_{Q_r} h_r(y)
\overline { \tilde \phi(2^\kappa y -\nu_\kappa) } dy \biggr|
\times | \phi(2^\kappa x -\nu_\kappa)| \\
\le c_{23} \alpha \int_{ Q_r } \rho(y, F) | x -y |^{-2} dy
\end{multline*}
почти для всех $ x \in F_{r, \iota}, r \in \N, \iota \in \Z, \Nu \subset \Z: \card \Nu < \infty. $
Из (2.1.18), (2.1.30), (2.1.38) вытекает, что существует константа
$ c_{24}(\phi, \tilde \phi, \Phi, \tilde \Phi^\prime) > 0 $ такая, что
при $ r \in \N, \iota \in \Z, \Nu \subset \Z: \card \Nu < \infty $
почти для всех $ x \in F_{r, \iota} $ имеет место неравенство

\begin{multline*}
\sum_{\kappa =0}^{k}
\sum_{ \nu_\kappa \in \Nu} 2^\kappa
\times \biggl| \int_{Q_r} h_r(y)
\overline { \tilde \phi(2^\kappa y -\nu_\kappa) } dy \biggr|
\times | \phi(2^\kappa x -\nu_\kappa)| \\
\le c_{24} \alpha \int_{ Q_r } \rho(y, F) | x -y |^{-2} dy,
\end{multline*}
а, следовательно, при $ r \in \N, $ почти для всех $ x \in F $ при любом 
$ n \in \Z_+ $ справедливо неравенство
\begin{multline*} \tag{2.1.39}
\sum_{\kappa =0}^{k}
\sum_{ \nu_\kappa =-n}^n 2^\kappa
\times \biggl| \int_{Q_r} h_r(y)
\overline { \tilde \phi(2^\kappa y -\nu_\kappa) } dy \biggr|
\times | \phi(2^\kappa x -\nu_\kappa)| \\
\le c_{24} \alpha \int_{ Q_r } \rho(y, F) | x -y |^{-2} dy.
\end{multline*}
причем, в каждой точке $ x \in F $ левая часть (2.1.39) представляет
собой монотонно возрастающую последовательность.
Поэтому, переходя к пределу при $ n \to \infty $ в неравенстве (2.1.39) и
учитывая (2.1.17), приходим к неравенству
\begin{equation*} \tag{2.1.40}
| (T h_r)(x)| \le c_{25} \alpha \int_{ Q_r } \rho(y, F) | x -y |^{-2} dy,
\text{ почти для всех} x \in F, r \in \N.
\end{equation*}

Соединяя (2.1.15) и (2.1.40), с учетом (1.4.4), (1.4.5) находим, что
почти для всех $ x \in F $ при любом $ m \in \N $ выполняется неравенство
\begin{multline*} \tag{2.1.41}
| (T h)(x)| \le \sum_{r=1}^m c_{25} \alpha \int_{ Q_r } \rho(y, F)
| x -y|^{-2} dy +\| (T h_m^\prime) \|_{L_\infty(\R)} \\ =
c_{25} \alpha \int_{ \cup_{r=1}^m Q_r } \rho(y, F) | x -y|^{-2} dy
+\| (T h_m^\prime) \|_{L_\infty(\R)} \\ \le
c_{25} \alpha \int_W \rho(y, F) | x -y|^{-2} dy
+\| (T h_m^\prime) \|_{L_\infty(\R)}.
\end{multline*}

Принимая во внимание, что при $ m \in \N $ в силу
оценок (1.2.23) и (2.1.26), условия (1.4.4) и включения
$ f \in L_1(\R) $ имеет место соотношение
\begin{multline*}
\| (T h_m^\prime) \|_{L_\infty(\R)} \\
= \| \sum_{\kappa =0}^k \sigma_\kappa \cdot \mathcal E_\kappa^1
( h_m^\prime ) \|_{L_\infty(\R)} \le
\sum_{\kappa =0}^k \| \mathcal E_\kappa^1( h_m^\prime ) \|_{L_\infty)\R)} \\
\le \sum_{\kappa =0}^k (\| E_\kappa^1(h_m^\prime) \|_{L_\infty(\R)} +
\| E_{\kappa -1}^1(h_m^\prime) \|_{L_\infty(\R)}) \\
\le \sum_{\kappa =0}^k c_{26} 2^\kappa \| h_m^\prime \|_{L_1(\R)} \le c_{27} 2^k \| h_m^\prime \|_{L_1(\R)} \\ =
c_{27} 2^k \| \sum_{r=m+1}^\infty h_r \|_{L_1(\R)} \le
c_{27} 2^k \sum_{r=m+1}^\infty \| h_r \|_{L_1(\R)} \\ =
c_{27} 2^k \sum_{r=m+1}^\infty \int_\R | h_r (y) | dy =
c_{27} 2^k \sum_{r=m+1}^\infty \int_{ Q_r} | h_r (y) | dy \\ \le
c_{27} 2^k \sum_{r=m+1}^\infty 2 \int_{ Q_r} | f(y) | dy \to 0
\text{ при } m \to \infty,
\end{multline*}
из (2.1.41) следует, что почти для всех $ x \in F $ соблюдается неравенство
\begin{equation*}
| (T h(x)| \le c_{25} \alpha \int_W \rho(y, F) | x -y|^{-2} dy.
\end{equation*}
Отсюда, применяя (1.4.1), (1.4.3), приходим к неравенству
\begin{multline*} \tag{2.1.42}
\| T h \|_{L_1(F)} = \int_F | (T h)(x)| dx \le \int_F c_{25} \alpha \int_W
\rho(y, F) | x -y |^{-2} dy dx \\
= c_{25} \alpha \int_F \int_W \rho(y, F) | x -y |^{-2} dy dx \le
c_{26} \alpha \mes (\R \setminus F) \\
= c_{26} \alpha \mes W \le c_{27} \int_{\R} | f(x)| dx =
c_{27} \| f \|_{L_1(\R)}.
\end{multline*}

Подставляя (2.1.42) в (2.1.14), выводим неравенство
$$
\mes \{ x \in F: | (T h)(x)| > \alpha /2 \} \le
(c_{28} / \alpha) \| f \|_{L_1(\R)},
$$
которое в соединении с (2.1.12), (2.1.13) дает оценку
\begin{equation*} \tag{2.1.43}
\mes \{ x \in \R: | (T h)(x) | > \alpha /2 \} \le
(c_{29} / \alpha) \| f \|_{L_1(\R)}.
\end{equation*}

Объединяя (2.1.10), (2.1.11) и (2.1.43), приходим к (2.1.8).
Сопоставляя (2.1.6) -- (2.1.8) с (1.4.8) -- (1.4.10), заключаем, что
существует константа $ c_{30}(\phi, \tilde \phi, \Phi, \tilde \Phi^\prime, p) >0 $
такая, что при любых $ k \in \Z_+, \sigma = \{ \sigma_\kappa \in \{-1, 1\}:
\kappa =0, \ldots, k \}, 1 < p < 2 $ для вещественнозначной $ f \in L_p(\R) $
согласно (1.4.11) верно неравенство
$$
\| T f \|_{L_p(\R)} \le c_{30} \| f \|_{L_p(\R)},
$$
из которого с учетом замечаний перед выводом (2.1.6) следует (2.1.4)
при $ 1 < p < 2 $ для комплекснозначных функций $ f \in L_p(\R). $
Справедливость (2.1.4) при $ p =2 $ установлена при выводе (2.1.7).

В частности, при $ k \in \Z_+ $ для $ f \in L_p(\R), 1 < p \le 2, $
в виду (1.1.1), (2.1.4) имеем
\begin{equation*}
\| E_k^p f \|_{L_p(\R)} = \| \sum_{\kappa =0}^k (\mathcal E_\kappa^p f) \|_{L_p(\R)}
\le c_1 \| f \|_{L_p(\R)},
\end{equation*}
что совпадает с (2.1.5).

Теперь проверим соблюдение (2.1.4)
при $ 2 < p < \infty. $
Учитывая сказанное, рассмотрим при $ k \in \Z_+, 2 < p < \infty,
p^\prime = p/(p-1) \in (1,2), $ непрерывный линейный оператор
$ E_k^{p^\prime} = E_k: L_{p^\prime}(\R) \mapsto L_{p^\prime}(\R) $ и
сопряженный к нему непрерывный линейный оператор $ (E_k^{p^\prime})^*:
(L_{p^\prime})^* = L_p(\R) \mapsto l_p(\R), $ который обладает тем свойством,
что для $ f \in L_p(\R), g \in L_{p^\prime}(\R) $ выполняется равенство
$$
\int_\R ((E_k^{p^\prime})^* f) \cdot g dx = \int_\R f \cdot (E_k^{p^\prime} g) dx.
$$
Сопоставляя это равенство с (1.2.27), видим, что для $ f, g \in C_0^\infty(\R) $
соблюдается равенство
\begin{multline*}
\int_\R (E_k^p f) \cdot g dx = \int_\R f \cdot \overline{ E_k^{p^\prime}
\overline g} dx = \overline{\int_\R \overline f \cdot (E_k^{p^\prime} \overline g) dx}\\ =
\overline{\int_\R ((E_k^{p^\prime})^* \overline f) \cdot \overline g dx} =
\int_\R \overline{((E_k^{p^\prime})^* \overline f)} \cdot g dx,
\end{multline*}
из которого следует, что почти для всех $ x \in \R $ справедливо равенство
\begin{equation*} \tag{2.1.44}
(E_k^p f)(x) = (S \circ (E_k^{p^\prime})^* (S f))(x), f \in C_0^\infty(\R),
k \in \Z_+, 2 < p < \infty,
\end{equation*}
где отображение
$S: L_{p}(\R) \ni f \mapsto S f = \overline f \in
L_{p}(\R), $ является изометрией $ L_{p}(\R) $ на себя.
Теперь при $ k \in \Z_+, 2 < p < \infty $ для $ f \in L_p(\R) $ выберем последовательность
$ \{ f_n \in C_0^\infty(\R), n \in \N \}, $ сходящуюся к $ f $ в $ L_p(\R). $
Тогда, учитывая сказанное, последовательность
$ \{ S \circ (E_k^{p^\prime})^* (S f_n), n \in \N\} $ сходится к
$ (S \circ (E_k^{p^\prime})^* (S f)) $ в $ L_p(\R), $ и, благодаря, (2.1.44),
(1.2.23), последовательность $ \{(E_k^p f_n) =
(S \circ (E_k^{p^\prime})^* (S f_n)), n \in \N\} $ сходится к $ (E_k^p f) $
в $ L_\infty(\R),$ что влечет равенство (2.1.44) для $ f \in L_p(\R). $

Принимая во внимание, что при $ k \in \Z_+, \sigma = \{ \sigma_\kappa \in \{-1, 1\}:
\kappa =0, \ldots, k \}, 2 < p < \infty, p^\prime = p/(p-1) \in (1,2), $
в силу (2.1.44) справедливо равенство
$$
\sum_{\kappa =0}^k \sigma_\kappa \cdot \mathcal E_\kappa^p =
S \biggl(\sum_{\kappa =0}^k \sigma_\kappa \cdot \mathcal E_\kappa^{p^\prime}\biggr)^* S,
$$
а, значит,
\begin{multline*}
\biggl\| \biggl(\sum_{\kappa =0}^k \sigma_\kappa \cdot \mathcal E_\kappa^{p}\biggr)
\biggr\|_{\mathcal B(L_p(\R), L_p(\R))} =
\biggl\| S \biggl(\sum_{\kappa =0}^k \sigma_\kappa \cdot \mathcal
E_\kappa^{p^\prime}\biggr)^* S \biggr\|_{ \mathcal
B(L_{p}(\R), L_{p}(\R))} \\ =
\biggl\| \biggl(\sum_{\kappa =0}^k \sigma_\kappa \cdot \mathcal
E_\kappa^{p^\prime}\biggr)^* \biggr\|_{ \mathcal
B(L_{p}(\R), L_{p}(\R))} =
\biggl\| \biggl(\sum_{\kappa =0}^k \sigma_\kappa \cdot \mathcal
E_\kappa^{p^\prime}\biggr) \biggr\|_{ \mathcal
B(L_{p^\prime}(\R), L_{p^\prime}(\R))}.
\end{multline*}
Отсюда, учитывая выполнение (2.1.4) при $ p^\prime $ вместо $ p, $
заключаем, что (2.1.4) соблюдается при $ 2 < p < \infty, $ и, в частности,
выполняется (2.1.5). $ \square $

Используя технику проведения оценок, применявшуюся при доказательстве 
теоремы 11.1.1 из [7] (см. также доказательство теоремы 1 из [15] и 
теоремы 9 из [16, гл. 7]) можно показать, что имеет место 

Лемма 2.1.2 

Пусть выполнены условия предложения 1.2.4 и пусть существует
четная ограниченная убывающая на $ [0, \infty) $ функция $ \mu $, 
удовлетворяющая условию
\begin{equation*} \tag{2.1.45}
\mu(x) \ln(1 +x) \in L_1([0, \infty)),
\end{equation*}
и такая, что
$$	
| \phi(x)|, | \tilde \phi(x)| \le \mu(x)
$$
для всех $ x \in \R.$
Тогда при $ 1 < p < \infty $ существует константа
$ c_{31}(\phi, \tilde \phi, \mu, p) >0 $ такая, что при любом 
$ k \in \Z_+ $ для любого набора чисел $ \sigma = \{ \sigma_\kappa \in 
\{-1, 1\}: \kappa =0, \ldots, k \}, $
для $ f \in L_p(\R) $ справедливо неравенство
\begin{equation*} \tag{2.1.46}
\biggl\| \sum_{\kappa =0}^k \sigma_\kappa \cdot (\mathcal E_\kappa f)
\biggr\|_{L_p(\R)} \le c_{31} \| f \|_{L_p(\R)}.
\end{equation*}

Отметим, что соблюдение условий леммы 2.1.1 не влечет соблюдение всех условий 
леммы 2.1.2, т.е. лемма 2.1.1 не является следствием леммы 2.1.2.

Для доказательства леммы 2.1.2 используются вспомогательные 
утверждения, часть из которых взята без изменений из [7], а часть является 
естественной модификацией соответствующих утверждений из [7].

Лемма 2.1.3

Пусть $ \mu $ -- четная ограниченная суммируемая и убывающая на $ [0, \infty) $ 
функция. Тогда для $ x \in \R $ выполняется неравенство
\begin{equation*} \tag{2.1.47}
\sum_{\nu \in \Z} \mu(x -\nu) \le 2 \mu(0) +\int_\R \mu(t) dt,
\end{equation*}
и существует константа $ c_{32}(\mu) > 0 $ такая, что для $ x, y \in \R $ верно 
неравенство
\begin{equation*} \tag{2.1.48}
\sum_{\nu \in \Z} \mu(x -\nu) \mu(y -\nu) \le c_{32} \mu(| x -y| /2).
\end{equation*}

Доказательство.

Поскольку $ x \in \R $ можно представить в виде $ x = \xi +\eta, $ где 
$ \xi \in \Z, \eta \in [0, 1), $ то (2.1.46) достаточно доаказать при 
$ x \in [0, 1), $ что мы и предполагаем. Тогда, учитывая свойства $ \mu, $ имеем
\begin{multline*}
\sum_{\nu \in \Z} \mu(x -\nu) = \mu(x) +\mu(x -1) +\sum_{\nu =-1}^{-\infty} 
\mu(x -\nu) +\sum_{\nu =2}^\infty \mu(x -\nu) \\ \le
2 \mu(0) +\sum_{\nu =-1}^{-\infty} \int_{\nu}^{\nu +1} \mu(x -u) du +
\sum_{\nu =2}^\infty \int_{\nu -1}^\nu \mu(x -u) du \\ =
2 \mu(0) +\int_{-\infty}^0 \mu(x -u) du +\int_{1}^\infty \mu(x -u) du\\ \le 
2 \mu(0) +\int_\R \mu(x -u) du = 2 \mu(0) +\int_\R \mu(t) dt.
\end{multline*}

Для получения (2.1.48) заметим, что для $ x, y \in \R, \nu \in \Z $ соблюдается 
неравенство
$$
| x -y| = | x -\nu -(y -\nu)| \le | x -\nu| +| -(y -\nu)| = | x -\nu| +| y -\nu|,
$$
и, следовательно, либо $ | x -\nu| \ge | x -y| /2, $ либо $ | y -\nu| \ge | x -y| /2. $
Учитывая это обстоятельство, свойства функции $ \mu $ и (2.1.47), получаем
\begin{multline*}
\sum_{\nu \in \Z} \mu(x -\nu) \mu(y -\nu) =  
\sum_{\nu \in \Z} \mu(| x -\nu|) \mu(| y -\nu|) \\ = 
\sum_{\nu \in \Z: | x -\nu| \ge | x -y| /2} \mu(| x -\nu|) \mu(| y -\nu|) +  
\sum_{\nu \in \Z: | x -\nu| < | x -y| /2} \mu(| x -\nu|) \mu(| y -\nu|) \\ \le 
\sum_{\nu \in \Z: | x -\nu| \ge | x -y| /2} \mu(| x -\nu|) \mu(| y -\nu|) +  
\sum_{\nu \in \Z: | y -\nu| \ge | x -y| /2} \mu(| x -\nu|) \mu(| y -\nu|) \\ \le 
\sum_{\nu \in \Z: | x -\nu| \ge | x -y| /2} \mu(| x -y| /2) \mu(| y -\nu|) +  
\sum_{\nu \in \Z: | y -\nu| \ge | x -y| /2} \mu(| x -\nu|) \mu(| x -y| /2) \\ =
\mu(| x -y| /2) \biggl(\sum_{\nu \in \Z: | x -\nu| \ge | x -y| /2} \mu(y -\nu) +  
\sum_{\nu \in \Z: | y -\nu| \ge | x -y| /2} \mu(x -\nu)\biggr) \\ \le
\mu(| x -y| /2) \biggl(\sum_{\nu \in \Z} \mu(y -\nu) +  
\sum_{\nu \in \Z} \mu(x -\nu)\biggr) \le c_{32} \mu(| x -y| /2). \square
\end{multline*} 

Следствие 2.1.4

Пусть функция $ \mu $ удовлетворяет условиям леммы 2.1.3 и для функций 
$ \phi, \tilde \phi $ почти для всех $ x \in \R $ выполняется неравенство
\begin{equation*} \tag{2.1.49}
| \phi(x)|, | \tilde \phi(x)| \le \mu(x). 
\end{equation*}
Тогда при $ 1 \le p \le \infty $ оператор $ E_\kappa = E_\kappa^p, $ определяемый равенством
\begin{multline*} \tag{2.1.50}
(E_\kappa f)(x) := \sum_{\nu \in \Z} 2^\kappa \cdot (\int_\R f(y) 
\overline{\tilde \phi(2^\kappa y -\nu)} dy) \phi(2^\kappa x -\nu)\\ =
\int_\R f(y) \sum_{\nu \in \Z} 2^\kappa \cdot \phi(2^\kappa x -\nu) 
\overline{\tilde \phi(2^\kappa y -\nu)} dy, \\
\text{ почти для всех } x \in \R,
f \in L_p(\R), 1 \le p \le \infty,
\end{multline*}
является непрерывным оператором из $ L_p(\R) $ в $ L_p(\R), $ причем,
\begin{equation*} \tag{2.1.51}
\| E_\kappa \|_{\mathcal B(L_p(\R), L_p(\R))} \le c_{33}(\mu, p).  
\end{equation*}
При этом, если $ | \phi(x)| \le \mu(x) $ для всех $ x \in \R, $ то равенство 
(2.1.50) имеет место для всех $ x \in \R. $

Доказательство.

Пусть
$$
| \phi(2^\kappa x -\nu)| \le \mu(2^\kappa x -\nu)
$$
для всех $ \nu \in \Z.$
Тогда по лемме 2.1.3 частичные суммы ряда
$$
\sum_{\nu \in \Z} | \phi(2^\kappa x -\nu) \overline{\tilde \phi(2^\kappa y -\nu)}|
$$
почти для всех $ y $ мажорируются функцией $ c(\mu) \mu(2^{\kappa -1} (x -y))$,
которая как функция от аргумента $ y $ принадлежит $ L_{p^\prime}(\R) $ при 
любом $ 1 \le p^\prime \le \infty. $ 
Значит в силу неравенства Гельдера и теоремы Лебега имеет место равенство 
(2.150) и благодаря неравенству Гельдера и теореме Фубини выполняется 
неравенство
\begin{multline*}
\biggl\| E_\kappa f \|_{L_p(\R)} \le c \| 2^\kappa \int_\R | f(y)| 
\mu(2^{\kappa -1} (\cdot -y)) dy \biggr\|_{L_p(\R)} \\ =
c \biggl\| 2^\kappa \int_\R | f(\cdot -y)| \mu(2^{\kappa -1} y) dy \biggr\|_{L_p(\R)} \\ =
c 2^\kappa \biggl(\int_\R \biggl| \int_\R | f(x -y)| \mu(2^{\kappa -1} y) dy\biggr|^p dx\biggr)^{1 /p} \\ = 
c 2^\kappa \biggl(\int_\R \biggl| \int_\R | f(x -y)| (\mu(2^{\kappa -1} y))^{1 /p} 
(\mu(2^{\kappa -1} y))^{1 /p^\prime} dy\biggr|^p dx\biggr)^{1 /p} \\ \le 
c 2^\kappa \biggl(\int_\R \biggl(\int_\R | f(x -y)|^p \mu(2^{\kappa -1} y) dy\biggr) \cdot
\biggl(\int_\R \mu(2^{\kappa -1} y) dy\biggr)^{1 /p^\prime} dx\biggr)^{1 /p} \\ = 
c 2^\kappa \biggl(\int_\R \mu(2^{\kappa -1} y) dy\biggr)^{1 /p^\prime} 
\biggl(\int_\R \biggl(\int_\R | f(x -y)|^p \mu(2^{\kappa -1} y) dy\biggr) dx\biggr)^{1 /p} \\ = 
c 2^\kappa \biggl(\int_\R \mu(2^{\kappa -1} y) dy\biggr)^{1 /p^\prime} \cdot
\biggl(\int_\R \int_\R | f(x -y)|^p \mu(2^{\kappa -1} y) dx dy\biggr)^{1 /p} \\ = 
c 2^\kappa \biggl(\int_\R \mu(2^{\kappa -1} y) dy\biggr)^{1 /p^\prime} \cdot
\biggl(\int_\R \mu(2^{\kappa -1} y) \biggl(\int_\R | f(x -y)|^p dx\biggr) dy\biggr)^{1 /p} \\ = 
c 2^\kappa \biggl(\int_\R \mu(2^{\kappa -1} y) dy\biggr)\biggl(\int_\R | f(x)|^p dx\biggr)^{1 /p}  = 
2 c \biggl(\int_\R \mu(y) dy\biggr) \| f \|_{L_p(\R)}. \square
\end{multline*}

Будем рассматривать интервалы:
$$
Q_{\kappa, \nu} := (2^{-\kappa} \nu +2^{-\kappa} I), \kappa, \nu \in \Z, 
$$
и отрезки:
$$
\overline Q_{\kappa, \nu} := (2^{-\kappa} \nu +2^{-\kappa} \overline I), \kappa, \nu \in \Z. 
$$

Лемма 2.1.5

Пусть функция $ \mu $ удовлетворяет условиям леммы 2.1.2 и для функций 
$ \phi, \tilde \phi $ для всех $ x \in \R $ выполняется неравенство (2.1.49).
Тогда существует константа $ c_{34}(\mu) >0 $ такая, что при $ \mathbf \kappa, 
\mathbf \nu \in \Z $ для $ f \in L_1(\R): \supp f \subset 
\overline Q_{\mathbf \kappa, \mathbf \nu}, $ для $ x \in \R: | x -2^{-\mathbf \kappa} 
\mathbf \nu| \ge 2 2^{-\mathbf \kappa} $ справедливо неравенство
\begin{multline*} \tag{2.1.52}
\sum_{\kappa \ge \mathbf \kappa} \sum_{ \nu \in \Z} 2^\kappa \biggl(\int_\R | f(y)|  
\cdot | \tilde \phi(2^\kappa y -\nu)| dy\biggr) | \phi(2^\kappa x -\nu)| \\ \le
c_{34} \| f \|_{L_1(\R)} 2^{\mathbf \kappa} \gamma((2^{\mathbf \kappa} x -\mathbf \nu) /4),
\end{multline*}
где $ \gamma(x) $ -- функция, определяемая равенством
\begin{equation*}
\gamma(x) := \sum_{\kappa \ge 0} 2^\kappa \mu(2^\kappa x), x \in \R,
\end{equation*} 
является четной, убывает на $ [0, \infty), $ суммируема и ограничена на 
каждом интервале $ (\\delta, \infty), \delta > 0. $

Доказательство.

Учитывая (2.1.49), (2.1.48),  на основании теоремы Лебега о предельном 
переходе под знаком интеграла при $ x \in \R $ имеем
\begin{multline*} \tag{2.1.53}
\sum_{\kappa \ge \mathbf \kappa} \sum_{ \nu \in \Z} 2^\kappa \biggl(\int_\R | f(y)|  
\cdot | \tilde \phi(2^\kappa y -\nu)| dy\biggr) | \phi(2^\kappa x -\nu)| \\ =
\sum_{\kappa \ge \mathbf \kappa} 2^\kappa \cdot 
\int_{Q_{\mathbf \kappa, \mathbf \nu}} | f(y)| \cdot 
\sum_{ \nu \in \Z} | \tilde \phi(2^\kappa y -\nu)| \cdot 
| \phi(2^\kappa x -\nu)| dy \\ \le
\sum_{\kappa \ge \mathbf \kappa} 2^\kappa \cdot 
\int_{Q_{\mathbf \kappa, \mathbf \nu}} | f(y)| \cdot 
\sum_{ \nu \in \Z} \mu(2^\kappa y -\nu) \cdot \mu(2^\kappa x -\nu) dy \\ \le
c_{32} \sum_{\kappa \ge \mathbf \kappa} 2^\kappa \cdot 
\int_{Q_{\mathbf \kappa, \mathbf \nu}} | f(y)| \cdot \mu(2^\kappa | x -y| /2) dy.
\end{multline*}

Принимая во внимание, что для $ x \in \R: | x -2^{-\mathbf \kappa} \mathbf \nu| 
\ge 2 2^{-\mathbf \kappa}, y \in Q_{\mathbf \kappa, \mathbf \nu} $ соблюдается 
неравенство
\begin{multline*}
| x -y| = | x -2^{-\mathbf \kappa} \mathbf \nu +2^{-\mathbf \kappa} \mathbf \nu -y| \ge
| x -2^{-\mathbf \kappa} \mathbf \nu | -| 2^{-\mathbf \kappa} \mathbf \nu -y| \\ \ge
| x -2^{-\mathbf \kappa} \mathbf \nu | -2^{-\mathbf \kappa} \ge
| x -2^{-\mathbf \kappa} \mathbf \nu | -| x -2^{-\mathbf \kappa} \mathbf \nu | /2  =
| x -2^{-\mathbf \kappa} \mathbf \nu | /2,
\end{multline*}
и учитывая монотонность и четность $ \mu, $ выводим
\begin{multline*} \tag{2.1.54}
\sum_{\kappa \ge \mathbf \kappa} 2^\kappa \cdot 
\int_{Q_{\mathbf \kappa, \mathbf \nu}} | f(y)| \cdot \mu(2^\kappa | x -y| /2) dy
\\ \le
\sum_{\kappa \ge 0} 2^{\mathbf \kappa} 2^\kappa \cdot 
\int_{Q_{\mathbf \kappa, \mathbf \nu}} | f(y)| \cdot \mu(2^\kappa 
| 2^{\mathbf \kappa} x -\mathbf \nu| /4) dy \\ = 
2^{\mathbf \kappa} \sum_{\kappa \ge 0} 2^\kappa \cdot 
\mu(2^\kappa | 2^{\mathbf \kappa} x -\mathbf \nu| /4) \cdot
\int_{Q_{\mathbf \kappa, \mathbf \nu}} | f(y)| dy \\ = 
2^{\mathbf \kappa} \biggl(\sum_{\kappa \ge 0} 2^\kappa \cdot 
\mu(2^\kappa ( 2^{\mathbf \kappa} x -\mathbf \nu) /4)\biggr) \cdot \int_\R | f(y)| dy \\ = 
2^{\mathbf \kappa} \gamma(( 2^{\mathbf \kappa} x -\mathbf \nu) /4) \cdot 
\| f \|_{L_1(\R)}. 
\end{multline*}

Соединяя (2.1.53) и (2.1.54), приходим к (2.1.52).

Осталось проверить свойства функции $ \gamma. $
Четность и монотонность функции $ \gamma $ вытекают из соответствующих свойств 
функции $ \mu. $
Пусть $ x > \delta. $ Из монотонности $ \mu $ следует
\begin{multline*} 
\sum_{\kappa \ge 0} 2^\kappa \mu(2^\kappa x) = 
2 \sum_{\kappa \ge 0} \mu(2^\kappa x) \int_{2^{\kappa -1}}^{2^\kappa} dt \\ \le
2 \sum_{\kappa \ge 0} \int_{2^{\kappa -1}}^{2^\kappa} \mu(t x) dt = 
2 \int_{1/2}^\infty \mu(t x) dt  = 
(2 /x) \int_{x /2}^\infty \mu(t) dt\\ \le
(1 /x) \| \mu \|_{L_1(\R)} \le (1 / \delta) \| \mu \|_{L_1(\R)}.
\end{multline*}

Теперь используя предыдущее неравенство и (2.1.45), считая $ \delta < 0, $ имеем
\begin{multline*}
\int_\delta^\infty \gamma(x) dx \le
\int_\delta^\infty (2 /x) \int_{x /2}^\infty \mu(t) dt dx \\ =
2 \int_{\delta /2}^\infty (1 /x) \int_{x }^\infty \mu(t) dt dx \le 
2 \int_{\delta /2}^\infty (1 +2 /\delta) (1 /(x +1)) \int_{x }^\infty \mu(t) dt  dx \\ \le
c(\delta) \int_{0}^\infty (1 /(x +1)) \int_{x }^\infty \mu(t) dt dx =
c(\delta) \int_{0}^\infty \int_0^t (1 /(x +1)) \mu(t) dx dt \\ =
c(\delta) \int_{0}^\infty \ln(t +1) \mu(t) dt. \square 
\end{multline*}

Для функции $ f: \R \mapsto \C $ будем обозначать 
\begin{equation*}
f^{\kappa, \nu} := f \cdot \chi_{Q_{\kappa, \nu}}, \kappa, \nu \in \Z.
\end{equation*}

Следствие 2.1.6

В условиях леммы 2.1.2 существует константа $ c_{35}(\mu) > 0 $ такая, что
при любом $ k \in \N $ для любого набора чисел $ \sigma = \{ \sigma_\kappa \in 
\{-1, 1\}: \kappa =1, \ldots, k \}, $
для $ f \in L_1(\R) \cap L_2(\R), \mathbf \kappa, \mathbf \nu \in \Z $ и 
$ x \in \R: | x -2^{-\mathbf \kappa} 
\mathbf \nu| \ge 2 2^{-\mathbf \kappa}, $ справедливо неравенство 
\begin{equation*} \tag{2.1.55}
\biggl| \biggl(\biggl(\sum_{\kappa =1}^k \sigma_\kappa \mathcal E_\kappa\biggr)(f^{\mathbf \kappa, 
\mathbf \nu} -E_{\mathbf \kappa} f^{\mathbf \kappa, \mathbf \nu})\biggr)(x)\biggr| \le
c_{35} \| f^{\mathbf \kappa, \mathbf \nu}\|_{L_1(\R)} 2^{\mathbf \kappa} 
\gamma((2^{\mathbf \kappa} x -\mathbf \nu) /4).
\end{equation*}

Доказательство.

Прежде всего, замечая, что в силу (1.2.26), (1.2.30) при 
$ \kappa, \mathbf \kappa, \mathbf \nu \in \Z $ для $ f \in L_1(\R) \cap L_2(\R) $ 
имеют место равенства
\begin{equation*} 
\mathcal E_\kappa (f^{\mathbf \kappa, \mathbf \nu} -E_{\mathbf \kappa} f^{\mathbf \kappa,
\mathbf \nu}) = \begin{cases}
0, \text{при} \kappa \le \mathbf \kappa; \\
\mathcal E_\kappa f^{\mathbf \kappa, \mathbf \nu}, \text{при} \kappa > \mathbf \kappa,
\end{cases}
\end{equation*}
для $ x \in \R: | x -2^{-\mathbf \kappa} 
\mathbf \nu| \ge 2 2^{-\mathbf \kappa}, $ выводим
\begin{multline*} \tag{2.1.56}
\biggl| \biggl(\biggl(\sum_{\kappa =1}^k \sigma_\kappa \mathcal E_\kappa\biggr)(f^{\mathbf \kappa, 
\mathbf \nu} -E_{\mathbf \kappa} f^{\mathbf \kappa, \mathbf \nu})\biggr)(x)\biggr| \\ =
\biggl| \biggl(\sum_{\kappa =1}^k \sigma_\kappa \mathcal E_\kappa(f^{\mathbf \kappa, 
\mathbf \nu} -E_{\mathbf \kappa} f^{\mathbf \kappa, \mathbf \nu})\biggr)(x)\biggr|  =
\biggl| \sum_{\kappa =1, \ldots, k: \kappa > \mathbf \kappa} \sigma_\kappa 
(\mathcal E_\kappa f^{\mathbf \kappa, \mathbf \nu})(x)\biggr| \\ \le
\sum_{\kappa =1, \ldots, k: \kappa > \mathbf \kappa} | \sigma_\kappa| \cdot 
| (\mathcal E_\kappa f^{\mathbf \kappa, \mathbf \nu})(x)|  =
\sum_{\kappa =1, \ldots, k: \kappa > \mathbf \kappa} 
| (E_\kappa f^{\mathbf \kappa, \mathbf \nu})(x) - 
( E_{\kappa -1} f^{\mathbf \kappa, \mathbf \nu})(x)| \\ \le
\sum_{\kappa =1, \ldots, k: \kappa > \mathbf \kappa} 
(| (E_\kappa f^{\mathbf \kappa, \mathbf \nu})(x) | + 
| ( E_{\kappa -1} f^{\mathbf \kappa, \mathbf \nu})(x)|)  \le
2 \sum_{\kappa =1, \ldots, k: \kappa \ge \mathbf \kappa} 
| (E_\kappa f^{\mathbf \kappa, \mathbf \nu})(x) | \\ \le
2 \sum_{\kappa \ge \mathbf \kappa} 
\biggl| \sum_{ \nu \in \Z} 2^\kappa \biggl(\int_\R f^{\mathbf \kappa, \mathbf \nu}(y) \cdot 
\overline{ \tilde \phi(2^\kappa y -\nu)} dy\biggr) \phi(2^\kappa x -\nu)\biggr| \\ \le
2 \sum_{\kappa \ge \mathbf \kappa} 
\sum_{ \nu \in \Z} 2^\kappa \biggl(\int_\R | f^{\mathbf \kappa, \mathbf \nu}(y) \cdot 
\overline{ \tilde \phi(2^\kappa y -\nu)}| dy\biggr) | \phi(2^\kappa x -\nu)| \\ =
2 \sum_{\kappa \ge \mathbf \kappa} 
\sum_{ \nu \in \Z} 2^\kappa \biggl(\int_\R | f^{\mathbf \kappa, \mathbf \nu}(y)| \cdot 
| \tilde \phi(2^\kappa y -\nu)| dy\biggr) | \phi(2^\kappa x -\nu)|.
\end{multline*} 

Соединяя (2.1.56) с (2.1.52), приходим к (2.1.55) для $ x \in \R: 
| x -2^{-\mathbf \kappa} \mathbf \nu| \ge 2 2^{-\mathbf \kappa}. \square $ 
  
Лемма 2.1.7

Пусть выполнены условия леммы 2.1.2. Тогда существует константа $ c_{36}(\mu) > 0 $ 
такая, что при любом $ k \in \N $ для любого набора чисел $ \sigma = 
\{ \sigma_\kappa \in \{-1, 1\}: \kappa =1, \ldots, k \}, $
для $ f \in L_1(\R) \cap L_2(\R), \mathbf \kappa \in \Z $ справедливо неравенство 
\begin{equation*} \tag{2.1.57}
\biggl\| \sum_{\kappa =1}^k \sigma_\kappa \mathcal E_\kappa 
(E_{\mathbf \kappa} f) \biggr\|_{L_\infty(\R)} \le c_{36} 2^{\mathbf \kappa} \| f \|_{L_1(\R)}.
\end{equation*} 

Доказательство.

Принимая во внимание, что в силу (1.2.26), (1.2.30) для $ f \in L_1(\R) 
\cap L_2(\R) $ при $ \kappa, \mathbf \kappa \in \Z $ соблюдаются равенства
\begin{equation*}
\mathcal E_\kappa (E_{\mathbf \kappa} f) =
\begin{cases}
0, \text{при} \kappa > \mathbf \kappa; \\
\mathcal E_\kappa f, \text{при} \kappa \le \mathbf \kappa,
\end{cases}
\end{equation*}
с учетом (2.1.47) для $ x \in \R $ имеем
\begin{multline*}
\biggl| \biggl(\sum_{\kappa =1}^k \sigma_\kappa \mathcal E_\kappa 
(E_{\mathbf \kappa} f)\biggr)(x)\biggr|  =
\biggl| \sum_{\kappa =1, \ldots, k: \kappa \le \mathbf \kappa} \sigma_\kappa 
(\mathcal E_\kappa f)(x)\biggr| \\ \le
\sum_{\kappa =1, \ldots, k: \kappa \le \mathbf \kappa} | \sigma_\kappa| \cdot 
| (\mathcal E_\kappa f)(x)|  =
\sum_{\kappa =1, \ldots, k: \kappa \le \mathbf \kappa} 
| (E_\kappa f)(x) -(E_{\kappa -1} f)(x)| \\ \le
2 \cdot \sum_{\kappa =0, \ldots, k: \kappa \le \mathbf \kappa} 
| (E_\kappa f)(x)|  \le
2 \cdot \sum_{\kappa \in \Z: \kappa \le \mathbf \kappa} 
| (E_\kappa f)(x)| \\ \le
2 \cdot \sum_{\kappa \in \Z: \kappa \le \mathbf \kappa} 
\biggl| \sum_{\nu \in \Z} 2^\kappa \biggl(\int_\R f(y) 
\overline{\tilde \phi(2^\kappa y -\nu)} dy\biggr) \phi(2^\kappa x -\nu)\biggr| \\ \le
2 \sum_{\kappa \le \mathbf \kappa} 
\sum_{ \nu \in \Z} 2^\kappa \biggl(\int_\R | f(y) \cdot 
\overline{ \tilde \phi(2^\kappa y -\nu)}| dy\biggr) | \phi(2^\kappa x -\nu)| \\ =
2 \sum_{\kappa \le \mathbf \kappa} 
\sum_{ \nu \in \Z} 2^\kappa \biggl(\int_\R | f(y)| \cdot 
| \tilde \phi(2^\kappa y -\nu)| dy\biggr) | \phi(2^\kappa x -\nu)| \\ \le
2 \sum_{\kappa \le \mathbf \kappa} 2^\kappa 
\sum_{ \nu \in \Z} \biggl(\int_\R | f(y)| \cdot 
\| \tilde \phi \|_{L_\infty(\R)} dy\biggr) | \phi(2^\kappa x -\nu)| \\ \le
2 \sum_{\kappa \le \mathbf \kappa} 2^\kappa 
\| f \|_{L_1(\R)} \cdot \| \mu \|_{L_\infty(\R)} 
\sum_{ \nu \in \Z} \mu(2^\kappa x -\nu) \\ \le
c_{37} \| f \|_{L_1(\R)} \sum_{\kappa \le \mathbf \kappa} 2^\kappa  \le
c_{36} 2^{\mathbf \kappa} \| f \|_{L_1(\R)}, 
\end{multline*}
что влечет (2.1.57). $ \square $

Лемма 2.1.8

Пусть соблюдаются условия леммы 2.1.2, $ \mathbf \kappa, \mathbf \nu \in \Z $ и 
$ f \in L_1(\R): \supp f \subset \overline Q_{\mathbf \kappa, \mathbf \nu}. $ 
Тогда при $ x \in \R $ имеет место неравенство
\begin{equation*} \tag{2.1.58}
| (E_{\mathbf \kappa} f)(x)| \le 2^{\mathbf \kappa} \| f \|_{L_1(\R)} \cdot
\beta(2^{\mathbf \kappa} x -\mathbf \nu),
\end{equation*}
где $ \beta(x) $ -- четная суммируемая убывающая на $ [0, \infty) $ функция, 
обладающая тем свойством, что при $ l \in \Z_+, x \in \R: | x| \ge 1. $ 
выполняется неравенство
\begin{equation*} \tag{2.1.59}
\beta(2^l x) \le 2^{1 -l} \beta(x).
\end{equation*}

Доказательство.

В виду (2.1.52) при $ x \in \R: | 2^{\mathbf \kappa} x -\mathbf \nu| \ge 2, $ имеем
\begin{multline*}
| (E_{\mathbf \kappa} f)(x)| =
\biggl| \sum_{\nu \in \Z} 2^{\mathbf \kappa} \biggl(\int_\R f(y) 
\overline{\tilde \phi(2^{\mathbf \kappa} y -\nu)} dy\biggr) \phi(2^{\mathbf \kappa} x -\nu)\biggr| \\ \le
\sum_{ \nu \in \Z} 2^{\mathbf \kappa} \biggl(\int_\R | f(y) \cdot 
\overline{ \tilde \phi(2^{\mathbf \kappa} y -\nu)}| dy\biggr) | \phi(2^{\mathbf \kappa} x -\nu)| \\ =
\sum_{ \nu \in \Z} 2^{\mathbf \kappa} \biggl(\int_\R | f(y)| \cdot 
| \tilde \phi(2^{\mathbf \kappa} y -\nu)| dy\biggr) | \phi2^{\mathbf \kappa} x -\nu)| \\ \le
c_{34} \| f \|_{L_1(\R)} 2^{\mathbf \kappa} \gamma((2^{\mathbf \kappa} x -\mathbf \nu) /4).
\end{multline*}
И как видно из доказательства (2.1.57) , при $ x \in \R: | 2^{\mathbf \kappa} x 
-\mathbf \nu| \le 2, $ выполняется неравенство
\begin{multline*}
| (E_{\mathbf \kappa} f)(x)| \le
\sum_{ \nu \in \Z} 2^{\mathbf \kappa} \biggl(\int_\R | f(y)| \cdot 
| \tilde \phi(2^{\mathbf \kappa} y -\nu)| dy\biggr) | \phi(2^{\mathbf \kappa} x -\nu)| \\ \le
2^{\mathbf \kappa} \sum_{ \nu \in \Z} \biggl(\int_\R | f(y)| \cdot 
\| \tilde \phi \|_{L_\infty(\R)} dy\biggr) | \phi(2^{\mathbf \kappa} x -\nu)| \\ \le
2^{\mathbf \kappa} \| f \|_{L_1(\R)} \cdot \| \mu \|_{L_\infty(\R)} 
\sum_{ \nu \in \Z} \mu(2^{\mathbf \kappa} x -\nu)  \le
c_{38} 2^{\mathbf \kappa} \| f \|_{L_1(\R)}. 
\end{multline*}

Полагая
\begin{equation*}
\beta(x) = \begin{cases} c_{39} \gamma(x /4), \text{ при } | x| \ge 2; \\
\beta(2) = c_{39} \gamma(1 /2), \text{ при } | x| < 2,
\end{cases}
\end{equation*}
с константой $ c_{39} = \max(c_{34}, c_{38} / \gamma(1/2)), $
убеждаемся в справедливости (2.1.58).

Для получения соотношения (2.1.59) при $ | x| \ge 2 $ имеем
\begin{multline*}
\beta(2 x) = c_{39} \gamma(2 x /4) = c_{39} \sum_{\kappa \ge 0} 2^\kappa 
\mu(2^\kappa 2 x /4)  =
c_{39} \sum_{\kappa \ge 0} 2^\kappa \mu(2^{\kappa +1} x /4) \\ =
(c_{39} /2) \sum_{\kappa \ge 0} 2^{\kappa +1} \mu(2^{\kappa +1} x /4)  =
(c_{39} /2) \sum_{\kappa \ge 1} 2^\kappa  \mu(2^\kappa x /4) \\ \le
(c_{39} /2) \sum_{\kappa \ge 0} 2^\kappa  \mu(2^\kappa x /4)  =
(c_{39} /2) \gamma(x /4) = (1 /2) \beta(x),
\end{multline*}
откуда по индукции относительно $ l $ следует, что при
$ l \in \Z_+, | x| \ge 2 $ выполняется неравенство
$$
\beta(2^l x) \le 2^{-l} \beta(x).
$$

Отсюда при $ x \in \R: 1 \le | x| < 2, l \in \N, $ с учетом четности и 
монотонности функции $ \beta $ выводим
\begin{multline*}
\beta(2^l x) = \beta(2^{l -1} 2 x) \le 2^{-(l -1)} \beta(2 x)\\ = 2^{-(l -1)} 
\beta(2| x|) \le 2^{-(l-1)} \beta(| x|) = 2^{-(l -1)} \beta(x).
\end{multline*}
Сопоставляя полученные оценки, приходим к (2.1.59). $ \square $

Лемма 2.1.9

Пусть выполнены условия леммы 2.1.2. Тогда существует кнстанта $ c_{40}(\mu) >0 $ 
такая, что для любой функции $ f \in L_1(\R) \cap L_2(\R), \alpha \in \R_+$ 
и множества $ S = \{(\kappa_r, \nu_r) \in \Z \times \Z, r \in \N\}, $ для 
которых двоичные интервалы $ \{ Q_{\kappa, \nu}, (\kappa, \nu) \in S \} $ 
попарно не пересекаются и для любого $ (\kappa, \nu) \in S $ соблюдается 
неравенство
\begin{equation*}
\| f^{\kappa, \nu} \|_{L_1(\R)} \le 2^{1 -\kappa} \alpha.
\end{equation*}
справедливо неравенство
\begin{equation*} \tag{2.1.60}
\int_\R \biggl| \sum_{(\kappa, \nu) \in S} E_\kappa f^{\kappa, \nu}\biggr|^2 dx \le 
c_{40} \alpha \| f \|_{L_1(\R)}.
\end{equation*}   

Доказательство.

Сначала установим справедливость (2.1.60) в случае, когда $ S $ -- 
конечное множество.

Перепишем левую часть (2.1.60) в виде
$$
\int_\R \biggl| \sum_{(\kappa, \nu) \in S} E_\kappa f^{\kappa, \nu}\biggr|^2 dx  =
\int_\R \sum_{(\kappa, \nu) \in S} \sum_{(\kappa^\prime, \nu^\prime) \in S} 
(E_\kappa f^{\kappa, \nu}) \cdot \overline{(E_{\kappa^\prime} 
f^{\kappa^\prime, \nu^\prime})} dx := J.
$$
Используя (2.1.58), и заменяя переменную в интеграле, имеем
\begin{multline*} \tag{2.1.61}
J \le 2 \sum_{(\kappa^\prime, \nu^\prime) \in S} \sum_{(\kappa, \nu) \in S: 
\kappa \ge \kappa^\prime} \int_\R | (E_{\kappa^\prime} f^{\kappa^\prime, 
\nu^\prime})(x) (E_\kappa f^{\kappa, \nu})(x)| dx \\ \le
2 \sum_{(\kappa^\prime, \nu^\prime) \in S} 2^{\kappa^\prime} 
\| f^{\kappa^\prime, \nu^\prime} \|_{L_1(\R)} \sum_{(\kappa, \nu) \in S: 
\kappa \ge \kappa^\prime} 2^\kappa \| f^{\kappa, \nu} \|_{L_1(\R)} 
\int_\R \beta(2^{\kappa^\prime} x -\nu^\prime) \beta(2^\kappa x -\nu) dx \\ \le
4 \alpha \sum_{(\kappa^\prime, \nu^\prime) \in S} 2^{\kappa^\prime} 
\| f^{\kappa^\prime, \nu^\prime} \|_{L_1(\R)} \sum_{(\kappa, \nu) \in S:
\kappa \ge \kappa^\prime} \int_\R \beta(2^{\kappa^\prime} x -\nu^\prime) 
\beta(2^\kappa x -\nu) dx \\ =
4 \alpha \sum_{(\kappa^\prime, \nu^\prime) \in S} \| f^{\kappa^\prime, \nu^\prime} \|_{L_1(\R)} 
\sum_{(\kappa, \nu) \in S^\prime: \kappa \ge 0}
\int_\R \beta(x) \beta(2^\kappa x -\nu) dx,
\end{multline*}     

где множество $ S^\prime = S^\prime(\kappa^\prime, \nu^\prime) $ состоит из 
пар $ (\kappa -\kappa^\prime, \nu -2^{\kappa -\kappa^\prime} \nu^\prime): (\kappa, \nu) \in S. $

Ясно, что интервалы $ Q_{\kappa, \nu}: (\kappa, \nu) \in S^\prime, \kappa \ge 0, $ 
попарно дизъюнктны.

Докажем равномерную ограниченность внутренней суммы в правой части 
(2.1.61). Обозначим через $ S_{n, \kappa}, n \in \Z, \kappa \in \Z_+, $ 
множество пар $ (\kappa, \nu) \in S^\prime, \kappa \ge 0, $ для которых 
двоичные интервалы $ Q_{\kappa, \nu} $ лежат в отрезке $ [n, n +1], $ 
и пусть $ C_{n, \kappa} $ -- число элементов в $ S_{n, \kappa}. $ Тогда
\begin{multline*}
\sum_{(\kappa, \nu) \in S^\prime: \kappa \ge 0} \int_\R \beta(x) \beta(2^\kappa x -\nu) dx  =
\sum_{n \in \Z} \sum_{\kappa =0}^\infty \sum_{(\kappa, \nu) \in S_{n, \kappa}} 
\int_\R \beta(x) \beta(2^\kappa x -\nu) dx \\ =
\sum_{n \in \Z} \int_\R \sum_{ \kappa \in \Z_+} 
\sum_{(\kappa, \nu) \in S_{n, \kappa}} \beta(x) \beta(2^\kappa x -\nu) dx.
\end{multline*}

Разобьём последний интеграл по схеме
$$
\int_\R = \int_{n -1}^{n +2} +\int_{-\infty}^{n -1} +\int_{n +2}^\infty :=
J_{n, 1} +J_{n, 2} +J_{n, 3}. 
$$
Из неравенства
\begin{multline*}
J_{n, 1} \le \left( \max_{[n -1], [n +2]} \beta\right) \sum_{\kappa \in \Z_+} 
\sum_{(\kappa, \nu) \in S_{n, \kappa}} \int_\R \beta(2^\kappa x -\nu) dx \\ =
\biggl( \max_{[n -1], [n +2]} \beta\biggr) 
\sum_{\kappa \in \Z_+} 2^{-\kappa} C_{n,\kappa} \int_\R \beta(x) dx  \le
\biggl(\max_{[n -1, n +2]} \beta\biggr) \int_\R \beta(x) dx,
\end{multline*}
принимая во внимание монотонность и суммируемость функции $ \beta, $ получаем
равномерную ограниченность суммы $ \sum_{n \in \Z} J_{n, 1}. $

Для оценки интегралов $ J_{n, 2}, J_{n, 3} $ заметим, что из включения 
$ Q_{\kappa, \nu} \subset [n, n +1] $ следует неравенство
$$
n \le 2^{-\kappa} \nu < n +1.
$$
Поэтому, если $ x \le n -1, $  то 
\begin{multline*}
| 2^\kappa x -\nu| = 2^\kappa | x -2^{-\kappa} \nu| = 2^\kappa | x -(n +1) 
+(n +1) -2^{-\kappa} \nu| \\
\ge 2^\kappa (| x -(n +1) | -| (n +1) -2^{-\kappa} 
\nu| ) \\ \ge
2^\kappa | x -(n +1)| (1 -| (n +1) -2^{-\kappa} \nu| /| x -(n +1)|) \\ \ge
2^\kappa | x -(n +1)| (1 -1/2) = 2^{\kappa -1} | x -(n +1)|,
\end{multline*} 
а, следовательно, с учетом четности и монотонности функции $ \beta $
\begin{equation*} \tag{2.1.62}
\beta(2^\kappa x -\nu) \le \beta(2^{\kappa -1} (x -(n +1))).        
\end{equation*}

Если же $ x \ge n +2, $ то
\begin{multline*}
| 2^\kappa x -\nu| = 2^\kappa | x -2^{-\kappa} \nu| = 2^\kappa | x -n  
+n -2^{-\kappa} \nu| \ge 2^\kappa (| x -n| -| n -2^{-\kappa} \nu|) \\ =
2^\kappa | x -n| (1 -| n -2^{-\kappa} \nu| /| x -n|)  \ge
2^\kappa | x -n| (1 -1/2) = 2^{\kappa -1} | x -n|,
\end{multline*}
и, значит,
\begin{equation*} \tag{2.1.63}
\beta(2^\kappa x -\nu) \le \beta(2^{\kappa -1} (x -n)).        
\end{equation*}

Применяя (2.1.62) и (2.1.59), имеем
\begin{multline*}
J_{n, 2} \le \sum_{\kappa \in \Z_+} \sum_{(\kappa, \nu) \in S_{n, \kappa}} 
\int_{-\infty}^{n -1} \beta(x) \beta(2^{\kappa -1} (n +1-x)) dx \\ =
\sum_{\kappa \in \Z_+} \sum_{(\kappa, \nu) \in S_{n, \kappa}} \int_2^\infty 
\beta(n +1 -x) \beta(2^{\kappa -1} x) dx \\ \le
\sum_{\kappa \in \Z_+} 2^{2 -\kappa} C_{n, \kappa} \int_2^\infty 
\beta(n +1-x) \beta(x) dx \le 2^2 \int_2^\infty \beta(n +1 -x) \beta(x) dx.
\end{multline*}
Отсюда, учитывая монотонность и суммируемость функции $ \beta, $ 
а также (2.1.47), получаем
\begin{equation*}
\sum_{n \in \Z} J_{n, 2} \le 2^2 \sum_{n \in \Z} \int_2^\infty \beta(n +1 -x) \beta(x) dx \le c_{41}.
\end{equation*}
И аналогично, используя (2.1.63) и (2.1.59), находим
\begin{multline*}
\sum_{n \in \Z} J_{n, 3} \le \sum_{n \in \Z} \int_{n +2}^\infty \sum_{ \kappa 
\in \Z_+} \sum_{(\kappa, \nu) \in S_{n, \kappa}} \beta(x) 
\beta(2^{\kappa -1} (x -n)) dx \\ =
\sum_{n \in \Z} \sum_{\kappa \in \Z_+} \sum_{(\kappa, \nu) \in S_{n, \kappa}}
\int_2^\infty \beta(n +x) \beta(2^{\kappa -1} x) dx \\ \le
2^2 \sum_{n \in \Z} \int_2^\infty \beta(n +x) \beta(x) dx \le c_{42}.
\end{multline*}

Сопоставляя полученные оценки с (2.1.61), получаем (2.1.60) в случае 
конечного множества $ S. $

Пусть теперь множество $ S $  -- бесконечно, т.е. $ S = \{(\kappa_r, 
\nu_r) \in \Z \times \Z: r \in \N\}. $ Для каждого $ m \in \N $ положим 
$ S^m := \{(\kappa_r, \nu_r): r =1, \ldots, m\}. $  Учитывая, что в силу 
(2.1.51) при $ (\kappa, \nu) \in S $ выполняется неравенство
\begin{equation*}
\| E_\kappa f^{\kappa, \nu} \|_{L_1(\R)} \le c_{33} \| f^{\kappa, \nu} \|_{L_1(\R)},
\end{equation*} 
и ряд 
\begin{equation*}
\sum_{r =1}^\infty \| f^{\kappa_r, \nu_r} \|_{L_1(\R)} =
\sum_{r =1}^\infty \int_{Q_{\kappa_r, \nu_r}} | f(x)| dx 
\end{equation*}
сходится, получаем, что ряд
\begin{equation*}
\sum_{r =1}^\infty E_{\kappa_r} f^{\kappa_r, \nu_r} 
\end{equation*}
сходится в $ L_1(\R) $ к некоторой функции $ g^\prime. $ 
Тогда выберем подпоследовательноть
\begin{equation*}
\biggl\{ \sum_{r =1}^{m_n} E_{\kappa_r} f^{\kappa_r, \nu_r}: m_n \in N, m_n < m_{n +1}\biggr\},
\end{equation*}
сходящуюся почти всюду на $ R $ к функции $ g^\prime, $
и применим к каждому элементу этой подпоследовательности неравенство
(2.1.60), в результате чего получим 
\begin{equation*}
\int_\R \biggl| \sum_{r =1}^{m_n} (E_{\kappa_r} f^{\kappa_r, \nu_r})(x)\biggr|^2 dx \le
c_{40} \alpha \| f \|_{L_1(\R)}.
\end{equation*}
Переходя в этом неравенстве к пределу на основании теоремы Фату, убеждаемся 
в справедливости (2.1.60) и для бесконечных множеств $ S. \square $

Лемма 2.1.10

Для любой функции $ f \in L_1(\R) $ и любого $ \alpha > 0 $ существует 
множество попарно дизъюнктных двоичных интервалов $ \{ Q_{\kappa, \nu}: 
(\kappa, \nu) \in S \} $ такое, что
\begin{equation*} \tag{2.1.64}
\alpha < (1 / \mes Q_{\kappa, \nu}) \int_{Q_{\kappa, \nu}} | f(x)| dx \le 
2 \alpha,
\end{equation*}
при $ (\kappa, \nu) \in S, $ а почти для всех $ x \in \R \setminus 
\cup_{(\kappa, \nu) \in S} Q_{\kappa, \nu} $ выполняется неравенство
\begin{equation*} \tag{2.1.65}
| f(x)| \le \alpha.
\end{equation*}

Доказательство леммы 2.1.10 приведено в [7].

Лемма 2.1.11 

В условиях леммы 2.1.2 существует константа $ c_{43}(\mu) > 0$ такая, что при 
любом $ k \in \N $ для любого набора чисел $ \sigma = \{ \sigma_\kappa \in 
\{-1, 1\}: \kappa =1, \ldots, k \}, $
для $ f \in L_1(\R), \alpha > 0 $
имеет место неравенство
\begin{equation*} \tag{2.1.66}
\mes \biggl\{x \in \R: \biggl| \sum_{\kappa =1}^k \sigma_\kappa (\mathcal E_\kappa f)(x)\biggr| 
> \alpha\biggr\} \le \frac{C_{43}}{\alpha} \| f \|_{L_1(\R)}.	
\end{equation*}

Доказательство. 

Обозначим $ T f = T^{k, \sigma} f := | \sum_{ \kappa =1}^k \sigma_\kappa 
(\mathcal E_\kappa f)|. $

Сначала установим справедливость (2.1.66) для $ f \in L_1(\R) \cap L_2(\R). $
 
Для $ f \in L_1(\R) \cap L_2(\R), \alpha >0 $ рассмотрим множество двоичных 
интервалов $ \{Q_{\kappa, \nu}: (\kappa, \nu) \in S\} $ из леммы 2.1.10. 
Из (2.1.64) следует, что
\begin{equation*} \tag{2.1.67}
\mes(\cup_{(\kappa, \nu) \in S} Q_{\kappa, \nu}) = 
\sum_{(\kappa, \nu) \in S} \mes(Q_{\kappa, \nu}) \le \frac{1}{\alpha} 
\| f \|_{L_1(\R)}.
\end{equation*}

Положим $ F = \R \setminus \cup_{(\kappa, \nu) \in S} Q_{\kappa, \nu}. $ 
Представим $ f $ в виде суммы функций
\begin{equation*} 
g := f \cdot \chi_F +\sum_{(\kappa, \nu) \in S} E_\kappa(f^{\kappa, \nu})  
\end{equation*} 
и
\begin{equation*} 
h := f -g = \sum_{(\kappa, \nu) \in S} (f^{\kappa, \nu} -E_\kappa(f^{\kappa, \nu})) = 
\sum_{r =1}^\infty h_r,  
\end{equation*}
где 
\begin{equation*}
h_r := (f^{\kappa_r, \nu_r} -E_{\kappa_r}(f^{\kappa_r, \nu_r})), r \in \N.
\end{equation*}
 
Ясно, что
\begin{multline*} \tag{2.1.68}
\mes \{x \in \R: | ( T f)(x)| > \alpha \}
 \le \mes \{x \in \R: | (T g)(x)| > \alpha /2 \} \\+
\mes \{x \in \R: | (T h)(x)| > \alpha /2 \}. 
\end{multline*}

Из (2.1.65) следует
\begin{equation*}
\int_F | f(x)|^2 dx \le \alpha \int_F | f(x)| dx \le \alpha \| f \|_{L_1(\R)}.
\end{equation*}
Отсюда и из леммы 2.1.9 (см. также (2.1.64)) получаем, что
\begin{multline*} 
\| g \|_{L_2(\R)}^2 = \int_\R \biggl| f \cdot \chi_F +\sum_{(\kappa, \nu) \in S} 
E_\kappa(f^{\kappa, \nu})\biggr|^2 dx \\ \le
2 \int_\R \biggl(| f \cdot \chi_F|^2 +\biggl| \sum_{(\kappa, \nu) \in S} 
E_\kappa(f^{\kappa, \nu})\biggr|^2\biggr) dx \\ =  
2 \biggl(\int_F | f|^2 dx +\int_\R \biggl| \sum_{(\kappa, \nu) \in S} 
E_\kappa(f^{\kappa, \nu})\biggr|^2 dx\biggr)  \le  
c_{44} \alpha \| f \|_{L_1(\R)}. 
\end{multline*} 

Имея в виду выкладку после неравенства (2.1.7), выводим
\begin{multline*} \tag{2.1.69}
\mes \{x \in \R: | (T g)(x)| > \alpha /2 \} = \mes \{x \in \R: | (T g)(x)|^2 
> \alpha^2 /4 \} \\
\le \frac{4}{\alpha^2} \| T g \|_{L_2(\R)}^2 \le 
\frac{4}{\alpha^2} \| g \|_{L_2(\R)}^2 \le \frac{c_{45}}{\alpha^2} \alpha 
\| f \|_{L_1(\R)} = \frac{c_{45}}{\alpha} \| f \|_{L_1(\R)}. 
\end{multline*}

Далее, положим $ G = \cup_{(\kappa, \nu) \in S} Q_{\kappa, \nu}^\prime, $ 
где 
$$
Q_{\kappa, \nu}^\prime := \{ x \in \R: | x -2^{-\kappa} \nu| < 2 2^{-\kappa} \}.
$$

Из (2.1.67) следует, что
\begin{equation*} \tag{2.1.70}
\mes G \le \sum_{(\kappa, \nu) \in S} \mes(Q_{\kappa, \nu}^\prime) =
\sum_{(\kappa, \nu) \in S} 4 \mes(Q_{\kappa, \nu}) \le \frac{4}{\alpha} \| f \|_{L_1(\R)}.   
\end{equation*}

С другой стороны,
\begin{equation*} \tag{2.1.71}
\mes \{ x \in \R \setminus G: | (T h)(x) | > \alpha /2 \} 
\le (2 / \alpha) \| T h \|_{L_1(\R \setminus G)}.
\end{equation*}

Для проведения оценки правой части (2.1.71) определим при $ m \in \N $
функцию $ h_m^\prime $ равенством
$$
h_m^\prime = h -\sum_{r=1}^m h_r
$$
и заметим, что вследствие (2.1.6) при $ m \in \N $ почти для всех
$ x \in (\R \setminus G) $ справедливо неравенство
$$
| (T h)(x) | = \biggl| \biggl(T \biggl(\sum_{r=1}^m h_r +h_m^\prime\biggr)\biggr)(x) \biggr| \le
\sum_{r=1}^m | (T h_r)(x) | +| (T h_m^\prime)(x) |,
$$
которое влечет оценку
\begin{multline*} 
\| T h \|_{L_1(\R \setminus G)} = \int_{\R \setminus G} | (T h)(x) | dx
\le \int_{\R \setminus G} \biggl(\sum_{r=1}^m | (T h_r)(x) |\biggr) +| (T
h_m^\prime)(x) | dx \\
= \sum_{r=1}^m \int_{\R \setminus G} | (T h_r)(x) | dx +\int_{\R \setminus G} 
| (T h_m^\prime)(x) | dx \\ =
\sum_{r =1}^m \int_{\R \setminus G}| T (f^{\kappa_r, \nu_r} 
-E_{\kappa_r}(f^{\kappa_r, \nu_r}))| dx +\int_{\R \setminus G} 
| (T h_m^\prime)(x) | dx \\ \le
\sum_{r =1}^m \int_{\R \setminus Q_{\kappa_r, \nu_r}^\prime}
\biggl| \biggl(\sum_{\kappa =1}^k \sigma_\kappa \mathcal E_\kappa\biggr) (f^{\kappa_r,
\nu_r} -E_{\kappa_r}(f^{\kappa_r, \nu_r}))\biggr| dx +\int_{\R} 
| (T h_m^\prime)(x) | dx.
\end{multline*}

Применив (2.1.55), получим

\begin{multline*} \tag{2.1.72}
\int_{\R \setminus G}| T h| dx\\ \le 
\sum_{r =1}^m \int_{\R \setminus Q_{\kappa_r, \nu_r}^\prime}
\biggl| \biggl(\sum_{\kappa =1}^k \sigma_\kappa \mathcal E_\kappa\biggr) (f^{\kappa_r,
\nu_r} -E_{\kappa_r}(f^{\kappa_r, \nu_r}))\biggl| dx +\int_{\R} 
| (T h_m^\prime)(x) | dx \\ \le
\sum_{r =1}^m \int_{x \in \R: | x -2^{-\kappa_r} 
\nu_r| \ge 2 2^{-\kappa_r}} c_{35} \| f^{\kappa_r, \nu_r}\|_{L_1(\R)} \\
2^{\kappa_r} \gamma((2^{\kappa_r} x -\nu_r) /4) dx +\int_{\R} 
| (T h_m^\prime)(x) | dx  =
c_{35} \sum_{r =1}^m \| f^{\kappa_r, \nu_r}\|_{L_1(\R)}\\ \times
\int_{x \in \R: | x -2^{-\kappa_r} \nu_r| \ge 2 2^{-\kappa_r}} 
2^{\kappa_r} \gamma((2^{\kappa_r} x -\nu_r) /4) dx +\int_{\R} 
| (T h_m^\prime)(x) | dx \\ =
c_{35} \sum_{r =1}^m \| f^{\kappa_r, \nu_r}\|_{L_1(\R)} 
\int_{x \in \R: | x | \ge 2} \gamma(x /4) dx +\int_{\R} 
| (T h_m^\prime)(x) | dx \\ =
c_{46} \sum_{r =1}^m \| f^{\kappa_r, \nu_r}\|_{L_1(\R)} +\int_{\R} 
| (T h_m^\prime)(x) | dx \\
\le c_{46} \| f \|_{L_1(\R)} +\int_{\R} | (T h_m^\prime)(x) | dx.
\end{multline*}

Принимая во внимание, что при $ m \in \N $ в силу
оценки (2.1.51), дизъюнкности интервалов $ Q_{\kappa, \nu}, (\kappa, \nu) \in S, $
и включения $ f \in L_1(\R) $ имеет место соотношение
\begin{multline*}
\int_R | (T h_m^\prime)(x) | dx \le
\sum_{\kappa =1}^k \| \mathcal E_\kappa(h_m^\prime) \|_{L_1(\R)}
\le c_{47} k \| h_m^\prime \|_{L_1(\R)} \\ =
c_{47} k \| \sum_{r=m+1}^\infty h_r \|_{L_1(\R)} \le
c_{47} k \sum_{r=m+1}^\infty \| h_r \|_{L_1(\R)} \\ =
c_{47} k \sum_{r=m+1}^\infty \| f^{\kappa_r, \nu_r} -E_{\kappa_r} 
f^{\kappa_r, \nu_r} \|_{L_1(\R)} \\ \le
c_{47} k \sum_{r=m+1}^\infty (\| f^{\kappa_r, \nu_r} \|_{L_1(\R)} 
+\| E_{\kappa_r} f^{\kappa_r, \nu_r} \|_{L_1(\R)}) \\ \le
c_{48} k \sum_{r=m+1}^\infty \| f^{\kappa_r, \nu_r} \|_{L_1(\R)} =
c_{48} k \sum_{r=m+1}^\infty \int_{ Q_{\kappa_r, \nu_r}} | f(x)| dx \to 0
\text{ при } m \to \infty,
\end{multline*}
из (2.1.72) следует, что выполняется неравенство
\begin{equation*} \tag{2.1.73}
\| T h \|_{L_1(\R \setminus G)} \le c_{46} \| f \|_{L_1(\R)}.
\end{equation*}

Используя соотношения (2.1.70), ( 2.1.71), (2.1.73), имеем
\begin{multline*} \tag{2.1.74}
\mes \{x \in \R: | (T h)(x)| > \alpha /2 \}\\ 
\le (\mes  G) +\mes \{x \in \R 
\setminus G: | (T h)(x)| > \alpha /2 \}\\
\le \frac{4}{\alpha} \| f \|_{L_1(\R)} 
+\frac{c_{49}}{\alpha} \| f \|_{L_1(\R)} = \frac{c_{50}}{\alpha} \| f \|_{L_1(\R)}.  
\end{multline*}

Для доказательства (2.1.66) в случае, когда $ f \in L_1(\R) \cap L_2(\R), $ 
осталось сопоставить (2.1.68) с (2.1.69) и (2.1.74). 

Если $ f \in L_1(\R), \alpha >0, $ то  выбирая с учетом (2.1.51) $ f^\prime \in 
L_1(\R) \cap L_2(\R) $ так, чтобы соблюдались неравенства
$$
\| f -f^\prime \|_{L_1(\R)} < \| f \|_{L_1(\R)}, 
\biggl\| \sum_{\kappa =1}^k \sigma_\kappa \mathcal E_\kappa (f -f^\prime) \biggr\|_{L_1(\R)} <
\| f \|_{L_1(\R)}, 
$$
благодаря (2.1.66), примененному к $ f^\prime, $ и неравенству Чебышева, получаем
\begin{multline*}
\mes \biggl\{x \in \R: \biggl| \sum_{\kappa =1}^k \sigma_\kappa (\mathcal E_\kappa f)(x)\biggr| 
> \alpha\biggr\} \\ \le
\mes \biggl\{x \in \R: \biggl| \sum_{\kappa =1}^k \sigma_\kappa (\mathcal E_\kappa f^\prime )(x)\biggr| 
> \alpha /2\biggr\}\\ +
\mes \biggl\{x \in \R: \biggl| \sum_{\kappa =1}^k \sigma_\kappa (\mathcal E_\kappa (f 
-f^\prime))(x)\biggr| > \alpha /2\biggr\} \\ \le
(c_{51} /\alpha) \| f^\prime \|_{L_1(\R)} +(2 /\alpha) \biggl\| \sum_{\kappa =1}^k 
\sigma_\kappa (\mathcal E_\kappa (f -f^\prime))\biggr\|_{L_1(\R)} \\ 
\le (c_{51} /\alpha) (\| f \|_{L_1(\R)} +\| f -f^\prime \|_{L_1(\R)} +(2 /\alpha)
\| f \|_{L_1(\R)} \le (c_{43} /\alpha) \| f \|_{L_1(\R)}. \square 
\end{multline*}

Для доказательства леммы 2.1.2 достаточно с очевидными изменениями 
повторить доказательство леммы 2.1.1 и вместо (2.1.8) использовать (2.1.66).

Дальнейшие рассмотрения проводятся на основе леммы 2.1.1. Однако все 
построения можно провести, опираясь на лемму 2.1.2.

\bigskip

2.2. В этом пункте описываются некоторые свойства проекторов $ E_\kappa^p,
\mathcal E_\kappa^p, \kappa \in \Z_+, 1 < p < \infty, $ установленные в п. 1.2.
при $ p =2. $

Предложение 2.2.1

Пусть выполнены условия леммы 2.1.1 и существует неотрицательная суммируемая
на $ \R $ функция $ \tilde Phi $ такая, что
почти для всех $ x \in \R $ соблюдается неравенство
\begin{equation*} \tag{2.2.1}
| \tilde \phi(x)| \le \int_{(1/2) B^1} \tilde \Phi(x -u) du.
\end{equation*}
Тогда при $ \kappa \in \Z_+ $ справедливы соотношения:

1) при $ 1 < P < \infty $ для $ f \in L_p(\R),
g \in L_{p^\prime}(\R) $ имеют место
равенства
\begin{equation*} \tag{2.2.2}
\int_{\R} (E_\kappa^p f) \cdot \overline g dx = \int_{\R} f \cdot
\overline {(E_\kappa^{p^\prime} g) } dx
\end{equation*}
и
\begin{equation*} \tag{2.2.3}
\int_{\R} (\mathcal E_\kappa^p f) \cdot \overline g dx = \int_{\R} f \cdot
\overline {(\mathcal E_\kappa^{p^\prime} g) } dx;
\end{equation*}

2) при $ 1 < p < \infty $ оператор $ E_\kappa^p $ является оператором
проектирования в $ L_p(\R) $, для которого
\begin{multline*} \tag{2.2.4}
\Im E_\kappa^p = \close_{L_p(\R)} (\span \{ \phi(2^\kappa \cdot -\nu),
\nu \in \Z \}) \\ =
h_{2^\kappa} (\close_{L_p(\R)} (\span
\{ \phi(\cdot -\nu), \nu \in \Z \}));
\end{multline*}
\begin{multline*} \tag{2.2.5}
\Ker E_\kappa^p = \{ f \in L_p(\R): \int_\R f(x)
\overline {\tilde \phi(2^\kappa x -\nu)} dx =0 \ \forall \nu \in \Z \} \\ =
\{ f \in L_p(\R): \int_\R f(x)
\overline {\phi(2^\kappa x -\nu)} dx =0 \ \forall \nu \in \Z \} \\ =
\{ f \in L_p(\R): \int_\R f(x)
\overline {g(x)} dx =0 \ \forall g \in
\span \{\phi(2^\kappa \cdot -\nu),  \nu \in \Z\} \}.
\end{multline*}

Доказательство.

Для вывода (2.2.2) для $ f \in L_p(\R), g \in L_{p^\prime}(\R) $
достаточно выбрать последовательность $ \{ f_n \in L_p(\R) \cap L_2(\R),
n \in \N \} $, сходящуюся к $ f $ в $ L_p(\R), $
и последовательность $ \{g_n \in L_{p^\prime}(\R) \cap L_2(\R), n \in \N \}, $
сходящуюся к $ g $ в $ L_{p^\prime}(\R), $ применить к ним
доказанное ранее равенство (1.2.27), получив равенство
$$
\int_{\R} (E_\kappa^p f_n) \cdot \overline {g_n} dx = \int_{\R} f_n \cdot
\overline {(E_\kappa^{p^\prime} g_n) } dx, n \in \N,
$$

а затем, учитывая (2.1.5) и неравенство Гельдера, перейти к пределу в этом
равенстве при $ n \to \infty. $

Для получения (2.2.3), используя (2.2.2), имеем
\begin{multline*}
\int_{\R} (\mathcal E_\kappa^p f) \cdot \overline g dx =
\int_{\R} ((E_\kappa^p f)
-(E_{\kappa -1}^p f))
\cdot \overline g dx =
\int_{\R} (E_\kappa^p f) \cdot \overline g dx -
\int_{\R} (E_{\kappa -1}^p f) \cdot \overline g dx \\ =
\int_{\R} f \cdot \overline {(E_\kappa^{p^\prime} g) } dx -
\int_{\R} f \cdot \overline {(E_{\kappa -1}^{p^\prime} g) } dx =
\int_{\R} f \cdot (\overline {(E_\kappa^{p^\prime} g) } -
\overline {(E_{\kappa -1}^{p^\prime} g) }) dx =
\int_{\R} f \cdot \overline {(\mathcal E_\kappa^{p^\prime} g) } dx.
\end{multline*}

Проверим, что в условиях п. 2) $ E_\kappa^p $ является проектором.
Для этого, учитывая, что для $ f \in L_2(\R) \cap L_p(\R), $ благодаря (2.1.5), имеет место
соотношение $ E_\kappa^p f = E_\kappa^2 f \in L_2(\R) \cap L_p(\R), $ используя
(1.2.26) и тот факт, что $ U_\kappa $ -- оператор проектирования, выводим
\begin{equation*}
E_\kappa^p E_\kappa^p f =
E_\kappa^p U_\kappa f =
U_\kappa U_\kappa f = U_\kappa f = E_\kappa^p f.
\end{equation*}
Отсюда в силу непрерывности оператора $ E_\kappa^p: L_p(\R) \mapsto L_p(\R) $
(см. (2.1.5)) и плотности $ L_2(\R) \cap L_p(\R) $ в $ L_p(\R) $
при $ 1 < p < \infty $ вытекает, что $ E_\kappa^p E_\kappa^p f = E_\kappa^p f $
для $ f \in L_p(\R), $ т.е. $ E_\kappa^p $ -- проектор в $ L_p(\R). $

Для проверки справедливости (2.2.4) сначала установим включение
\begin{equation*} \tag{2.2.6}
\close_{L_p(\R)} (\span \{ \phi(2^\kappa \cdot -\nu),
\nu \in \Z \}) \subset \Im E_\kappa^p.
\end{equation*}

С этой целью, беря для $ f \in \close_{L_p(\R)} (\span \{ \phi(2^\kappa \cdot -\nu),
\nu \in \Z \}) $ послдовательность $ \{ f_n \in
\span \{ \phi(2^\kappa \cdot -\nu), \nu \in \Z \}, n \in \N \}, $
сходящуюся к $ f $ в $ L_p(\R), $ и используя непрерывность оператора
$ E_\kappa^p $ в $ L_p(\R), $ а также (1.2.26), (1.2.16) и тот факт,
что $ U_\kappa $ суть проектор, получаем,
$$
E_\kappa^p f = \lim_{n \to \infty} E_\kappa^p f_n =
\lim_{n \to \infty} U_\kappa f_n =
\lim_{n \to \infty} f_n = f,
$$
следовательно, имеет место (2.2.6).

Теперь установим обратное включение
\begin{equation*} \tag{2.2.7}
\Im E_\kappa^p \subset \close_{L_p(\R)}
(\span \{ \phi(2^\kappa \cdot -\nu), \nu \in \Z \}).
\end{equation*}
При $ p = 2 $ (2.2.7) с учетом (1.2.26) следует из (1.2.16).
Выбирая для $ f \in L_p(\R), 1 < p < \infty, $
последовательность функций $ \{ f_n \in C_0^\infty(\R), n \in \N \}, $
(имеющих компактные носители), сходящуюся к $ f $ в $ L_p(\R), $
вследствие непрерывности оператора $ E_\kappa^p: L_p(\R) \mapsto L_p(\R) $
получаем, что $ E_\kappa^p f_n \to E_\kappa^p f $ в $ L_p(\R) $
при $ n \to \infty, $ т.е.
\begin{equation*}
\Im E_\kappa^p \subset \close_{L_p(\R)}
\{ E_\kappa^p f, f \in C_0^\infty(\R)\}.
\end{equation*}
Таким образом, чтобы доказать (2.2.7), достаточно убедиться в справедливости
соотношения
\begin{equation*} \tag{2.2.8}
\{ E_\kappa^p f, f \in C_0^\infty(\R)\}
\subset \close_{L_p(\R)}
(\span \{ \phi(2^\kappa \cdot -\nu), \nu \in \Z \}).
\end{equation*}
Для этого, фиксируя для $ f \in C_0^\infty(\R) $ компакт $ K = \supp f, $
при $ \kappa \in \Z_+, 1 < p < \infty, $ для  $ R \in \R: R \ge 1, $ рассмотрим
линейную комбинацию (конечную в виду компактности $ K $ )
\begin{multline*}
\sum_{\nu_\kappa \in \Z: \rho(2^{-\kappa} \nu_\kappa, K ) < R} 2^{\kappa}
\biggl(\int_{\R} f(y) \overline {\tilde \phi(2^\kappa y -\nu_\kappa)} dy\biggr) \times
\phi(2^\kappa x -\nu_\kappa)\\ \in \span \{ \phi(2^\kappa \cdot -\nu), \nu \in \Z \}.
\end{multline*}
Тогда в силу (1.2.22), (1.2.20) имеем
\begin{multline*} \tag{2.2.9}
\| E_\kappa^p f -\sum_{\nu_\kappa \in \Z: \rho(2^{-\kappa} \nu_\kappa, K ) < R}
2^{\kappa} \biggl(\int_{\R} f(y) \overline {\tilde \phi(2^\kappa y -\nu_\kappa)} dy\biggr) \times
\phi(2^\kappa \cdot -\nu_\kappa) \|_{L_p(\R)}\\ =
\biggl\| \sum_{\nu_\kappa \in \Z: \rho(2^{-\kappa} \nu_\kappa, K ) \ge R}
2^{\kappa} \biggl(\int_{\R} f(y) \overline {\tilde \phi(2^\kappa y -\nu_\kappa)} dy\biggr) \times
\phi(2^\kappa \cdot -\nu_\kappa) \biggr\|_{L_p(\R)} \\ \le
\sum_{\nu_\kappa \in \Z: \rho(2^{-\kappa} \nu_\kappa, K ) \ge R}
2^{\kappa} \times \biggl| \int_{\R} f(y) \overline {\tilde \phi(2^\kappa y -\nu_\kappa)} dy\biggr| \times
\| \phi(2^\kappa \cdot -\nu_\kappa) \|_{L_p(\R)}.
\end{multline*}
Замечая, что при $ \kappa \in \Z_+, \nu_\kappa \in \Z $ имеет место равенство
\begin{multline*}
\biggl\| \phi(2^\kappa \cdot -\nu_\kappa)\biggr \|_{L_p(\R)} =
\biggl(\int_\R | \phi(2^\kappa y -\nu_\kappa) |^p dy\biggr)^{1/p} \\ =
\biggl(\int_\R | \phi(2^\kappa (2^{-\kappa} \nu_\kappa +2^{-\kappa} z) -\nu_\kappa) |^p
2^{-\kappa} dz\biggr)^{1/p} \\ =
2^{-\kappa /p} \biggl(\int_\R | \phi(z)|^p dz\biggr)^{1/p} = 2^{-\kappa /p} \| \phi \|_{L_p(\R)},
\end{multline*}
а, благодаря (2.2.1) (см. также (1.2.24)), выполняется неравенство
\begin{multline*}
\biggl| \int_{\R} f(y) \overline {\tilde \phi(2^\kappa y -\nu_\kappa)} dy\biggr| =
\biggl| \int_{K} f(y) \overline {\tilde \phi(2^\kappa y -\nu_\kappa)} dy\biggr| \\ \le
\| f \|_{L_\infty(\R)} \times \int_{K} \biggl| \tilde \phi(2^\kappa y -\nu_\kappa)\biggr| dy \\ \le
\| f \|_{L_\infty(\R)} \times \int_{K}
2^\kappa \int_{2^{-\kappa}(\nu_\kappa +(1/2) B^1)} \tilde \Phi(2^\kappa (y -u)) du dy\\ =
2^\kappa \times \| f \|_{L_\infty(\R)}
\int_{K} \int_{2^{-\kappa}(\nu_\kappa +(1/2) B^1)}
\tilde \Phi(2^\kappa (y -u)) du dy,
\end{multline*}
и учитывая, что при $ \kappa \in \Z_+, R \ge 1,  \nu_\kappa \in \Z:
\rho(2^{-\kappa} \nu_\kappa , K) \ge R, y \in K,
u \in 2^{-\kappa}(\nu_\kappa +(1/2) B^1), $ соблюдается неравенство
\begin{multline*}
| y -u| = | y -2^{-\kappa} \nu_\kappa +2^{-\kappa} \nu_\kappa -u| \ge
| y -2^{-\kappa} \nu_\kappa| -| 2^{-\kappa} \nu_\kappa -u|\\ \ge
\rho(2^{-\kappa} \nu_\kappa, K) -(1/2) 2^{-\kappa} \ge R -1/2 \ge R/2,
\end{multline*}

из (2.2.9) в силу $\sigma$-аддитивности интеграла как функции множеств, теоремы
Лебега о предельном переходе под знаком интеграла получаем соотношение
\begin{multline*}
\biggl\| E_\kappa^p f -\sum_{\nu_\kappa \in \Z: \rho(2^{-\kappa} \nu_\kappa, K ) < R}
2^{\kappa} \biggl(\int_{\R} f(y) \overline {\tilde \phi(2^\kappa y -\nu_\kappa)} dy\biggr) \times
\phi(2^\kappa \cdot -\nu_\kappa) \biggr\|_{L_p(\R)} \\ \le
\sum_{\nu_\kappa \in \Z: \rho(2^{-\kappa} \nu_\kappa, K ) \ge R}
2^{\kappa} \times 2^\kappa \times \| f \|_{L_\infty(\R)} \\
\times\biggl(\int_{K} \int_{2^{-\kappa}(\nu_\kappa +(1/2) B^1)}
\tilde \Phi(2^\kappa (y -u)) du dy\biggr) \times
2^{-\kappa /p} \| \phi \|_{L_p(\R)} \\ =
c_1(\phi, p) \| f \|_{L_\infty(\R)} 2^{ \kappa(2 -1/p)}
\sum_{\nu_\kappa \in \Z: \rho(2^{-\kappa} \nu_\kappa, K ) \ge R}
\int_{K} \int_{2^{-\kappa}(\nu_\kappa +(1/2) B^1)}
\tilde \Phi(2^\kappa (y -u)) du dy \\ =
c_1 \| f \|_{L_\infty(\R)} 2^{ \kappa(2 -1/p)} \times
\int_{K} \biggl(\sum_{\nu_\kappa \in \Z: \rho(2^{-\kappa} \nu_\kappa, K ) \ge R}
\int_{2^{-\kappa}(\nu_\kappa +(1/2) B^1)}
\tilde \Phi(2^\kappa (y -u)) du\biggr) dy \\ =
c_1 \| f \|_{L_\infty(\R)} 2^{ \kappa(2 -1/p)} \times
\int_{K} \int_{\cup_{\nu_\kappa \in \Z: \rho(2^{-\kappa} \nu_\kappa, K ) \ge R}
(2^{-\kappa}(\nu_\kappa +(1/2) B^1))}
\tilde \Phi(2^\kappa (y -u)) du dy \\ \le
c_1 \| f \|_{L_\infty(\R)} 2^{ \kappa(2 -1/p)} \times
\int_{K} \int_{u \in \R: | y -u| \ge R/2}
\tilde \Phi(2^\kappa (y -u)) du dy \\ =
c_1 \| f \|_{L_\infty(\R)} 2^{ \kappa(2 -1/p)} \times
\int_{K} \int_{z \in \R: | z| \ge R/2}
\tilde \Phi(2^\kappa z) dz dy \\ =
c_1 \| f \|_{L_\infty(\R)} (\mes K)
2^{ \kappa(2 -1/p)} \times \int_{z \in \R: | z| \ge R/2} \tilde \Phi(2^\kappa z) dz\\  =
c_1 \| f \|_{L_\infty(\R)} (\mes K)
2^{ \kappa(2 -1/p)} \times \int_{\zeta \in \R: | \zeta| \ge 2^\kappa R/2}
\tilde \Phi(\zeta) 2^{-\kappa} d\zeta\\ =
c_1 \| f \|_{L_\infty(\R)} (\mes K)
2^{ \kappa(1 -1/p)} \times \int_{\zeta \in \R: | \zeta| \ge 2^{(\kappa -1)} R}
\tilde \Phi(\zeta) d\zeta \to 0\\
 \text{ при } R \to \infty,
\kappa \in \Z_+, 1 < p < \infty.
\end{multline*}

Тем самым, установлена справедливость (2.2.8), а вместе с ним и (2.2.7), что завершает доказательство (2.2.4).
Перейдем к выводу (2.2.5). Сначала получим первое равенство в (2.2.5) при $ \kappa =0. $
Пусть $ f \in \Ker E_0^p. $ Тогда при $ \mu \in \Z $ ввиду теоремы Лебега о
предельном переходе под знаком интеграла с учетом включения $ \tilde \phi \in L_1(\R) $
и соотношений (1.2.22), (1.2.21), (1.2.8) соблюдаются равенства
\begin{multline*}
0 = \int_\R (E_0^p f)(x) \overline {\tilde \phi(x -\mu)} dx \\ =
\int_\R \biggl(\sum_{\nu \in \Z} \biggl(\int_\R f(y) \overline {\tilde \phi(y -\nu)} dy\biggr)
\phi(x -\nu)\biggr) \overline {\tilde \phi(x -\mu)} dx \\ =
\sum_{\nu \in \Z} \int_\R \biggl(\int_\R f(y) \overline {\tilde \phi(y -\nu)} dy\biggr)
\phi(x -\nu) \overline {\tilde \phi(x -\mu)} dx \\ =
\sum_{\nu \in \Z} \biggl(\int_\R f(y) \overline {\tilde \phi(y -\nu)} dy\biggr) \times \int_\R
\phi(x -\nu) \overline {\tilde \phi(x -\mu)} dx =
\int_\R f(y) \overline {\tilde \phi(y -\mu)} dy.
\end{multline*}
Если же для $ f \in L_p(\R) $ выполняются равенства
$ \int_\R f(y) \overline {\tilde \phi(y -\nu)} dy =0 , \nu \in \Z, $
то согласно (1.2.22) $ E_0^p f =0. $
Таким образом, первое равенство в (2.2.5) при $ \kappa = 0 $ установлено.
Принимая во внимание, что $ E_\kappa^p =
h_{2^\kappa} E_0^p (h_{2^\kappa})^{-1}, $ на основании (1.2.2), (1.2.5)
и первого равенства в (2.2.5) при $ \kappa =0, $ получаем
\begin{multline*}
\Ker E_\kappa^p =
h_{2^\kappa} \biggl(\biggl\{ f \in L_p(\R): \int_\R f(x)
\overline {\tilde \phi( x -\nu)} dx =0 \ \forall \nu \in \Z \biggr\}\biggr) \\ =
\biggl\{ h_{2^\kappa} f| f \in L_p(\R): \int_\R f(x)
\overline {\tilde \phi( x -\nu)} dx =0 \ \forall \nu \in \Z \biggr\} \\ =
\biggl\{ h_{2^\kappa} f| f \in L_p(\R): \int_\R (h_{2^\kappa} f)(x)
\overline {(h_{2^\kappa} \tilde \phi(\cdot -\nu))(x)} dx =0 \ \forall \nu \in \Z \biggr\} \\ =
\biggl\{ f \in L_p(\R): \int_\R f(x)
\overline { \tilde \phi(2^\kappa x -\nu)} dx =0 \ \forall \nu \in \Z \biggr\}.
\end{multline*}
Для проверки справедливости остальных равенств в (2.2.5) заметим, что
для $ f \in L_p(\R), 1 < p < \infty, $ при $ \kappa \in \Z_+ $ вследствие
(2.2.4), (2.2.2) имеют место утверждения:
\begin{multline*}
\int_\R f(x) \overline{\phi(2^\kappa x -\nu)} dx = 0 \ \forall \nu \in \Z \\ \iff
\int_\R f(x) \overline{g(x)} dx =0 \ \forall g \in
\span \{\phi(2^\kappa \cdot -\nu), \nu \in \Z \} \\ \iff
\int_\R f(x) \overline{g(x)} dx =0 \ \forall g \in
\close_{L_{p^\prime}(\R)} (\span \{\phi(2^\kappa \cdot -\nu), \nu \in \Z \}) \\ \iff
\int_\R f \overline{ E_\kappa^{p^\prime} g} dx =0 \ \forall g \in L_{p^\prime}(\R) \iff
\int_\R E_\kappa^p f \cdot \overline{ g} dx =0 \ \forall g \in L_{p^\prime}(\R) \\ \iff
E_\kappa^p f =0 \iff f \in \Ker E_\kappa^p,
\end{multline*}
что влечет соблюдение остальных равенств в (2.2.5). $ \square $

Предложение 2.2.2

Пусть выполнены условия предложения 2.2.1. Тогда имеют место соотношения:

1) при $ 1 < p < \infty $ для $ \kappa, \kappa^\prime \in \Z_+:
\kappa^\prime \le \kappa, $ выполняются равенства
\begin{equation*} \tag{2.2.10}
E_{\kappa^\prime}^p E_\kappa^p = E_\kappa^p
E_{\kappa^\prime}^p =
E_{\kappa^\prime}^p;
\end{equation*}

2) при $ 1 < p < \infty, \kappa \in \Z_+ $ имеют место включения
\begin{equation*} \tag{2.2.11}
\Im E_\kappa^p \subset \Im E_{\kappa +1}^p; \\
\Ker E_{\kappa +1}^p \subset \Ker E_\kappa^p;
\end{equation*}

3) при $ 1 < p < \infty $ для $ \kappa, \kappa^\prime \in \Z_+ $
соблюдаются равенства
\begin{equation*} \tag{2.2.12}
\mathcal E_\kappa^p \mathcal E_{\kappa^\prime}^p =
\begin{cases} \mathcal E_\kappa^p,
  \text{при $ \kappa = \kappa^\prime $};\\
       0,  \text{при $ \kappa \ne \kappa^\prime $};
\end{cases}
\end{equation*}
и при $ 1 < p < \infty, \kappa, \kappa^\prime \in \Z_+: \kappa \ne \kappa^\prime, $
для $ f \in L_p(\R), g \in L_{p^\prime}(\R) $ имеет место равенство
\begin{equation*} \tag{2.2.13}
\int_{\R} (\mathcal E_\kappa^p f) \cdot \overline {(\mathcal
E_{\kappa^\prime}^{p^\prime} g)} dx = 0;
\end{equation*}

4) при $ 1  p < \infty $ для $ \kappa \in \N $ соблюдаются равенства
$$
\Im \mathcal E_\kappa^p = \Im E_\kappa^p \cap \Ker E_{\kappa -1}^p;
\Ker \mathcal E_\kappa^p = \Im E_{\kappa -1}^p +\Ker E_\kappa^p.
$$

Доказательство.

Убедимся в справедливости (2.2.10).
Поскольку при $ \kappa^\prime, \kappa \in \Z_+: \kappa^\prime \le \kappa, $
вследствие (1.2.30), (1.2.26) для $ f \in L_2(\R) \cap L_p(\R) $ имеет место
равенство
$$
(E_\kappa^p E_{\kappa^\prime}^p) f =
E_\kappa^p (E_{\kappa^\prime}^p f) =
E_\kappa^p (U_{\kappa^\prime} f) =
U_\kappa (U_{\kappa^\prime} f) =
U_{\kappa^\prime} f =
E_{\kappa^\prime}^p f,
$$
то выбирая для $ f \in L_p(\R) $ последовательность
$ \{f_n \in L_2(\R) \cap L_p(\R), n \in \N \}, $ сходящуюся к $ f $ в $ L_p(\R), $
и переходя к пределу при $ n \to \infty $ в равенстве
$$
(E_\kappa^p E_{\kappa^\prime}^p ) f_n = E_{\kappa^\prime}^p f_n,
$$
в силу непрерывности в $ L_p(\R) $ входящих в него операторов получаем
второе равенство в (2.2.10) для $ f \in L_p(\R).$ Точно так же устанавливается
справедливость первого равенства в (2.2.10).

Для доказательства первого включения в (2.2.11) для $ f \in \Im E_\kappa^p $
возьмем функцию $ g \in L_p(\R) $ такую, что $ f = E_\kappa^p g. $
Тогда ввиду (2.2.10) имеем
\begin{multline*}
f = E_{\kappa +1}^p f +f -E_{\kappa +1}^p f =
E_{\kappa +1}^p f +f -E_{\kappa +1}^p (E_\kappa^p g)\\ =
E_{\kappa +1}^p f +f -E_\kappa^p g = E_{\kappa +1}^p f +f -f = E_{\kappa +1}^p f \in
\Im E_{\kappa +1}^p,
\end{multline*}
что завершает доказательство первого включения в (2.2.11).

Чтобы установить соблюдение второго включения в (2.2.11), заметим, что
вследствие (2.2.5), (2.2.4) и первого включения в (2.2.11) с учетом
неравенства Гельдера вытекает, что
\begin{multline*}
\Ker E_{\kappa +1}^p =
\biggl\{ f \in L_p(\R): \int_\R f(x)
\overline {g(x)} dx =0 \ \forall g \in
\span \{\phi(2^{\kappa +1} \cdot -\nu),  \nu \in \Z\} \biggr\} \\ =
\biggl\{ f \in L_p(\R): \int_\R f(x)
\overline {g(x)} dx =0 \ \forall g \in
\close_{L_{p^\prime}(\R)} \span \{\phi(2^{\kappa +1} \cdot -\nu),
\nu \in \Z\} \biggr\} \\ =
\biggl\{ f \in L_p(\R): \int_\R f(x) \overline {g(x)} dx =0 \ \forall g \in
\Im E_{\kappa +1}^{p^\prime} \biggr\}\\ \subset
\biggl\{ f \in L_p(\R): \int_\R f(x) \overline {g(x)} dx =0 \ \forall g \in
\Im E_{\kappa}^{p^\prime} \biggr\} \\ =
\biggl\{ f \in L_p(\R): \int_\R f(x)
\overline {g(x)} dx =0 \ \forall g \in
\close_{L_{p^\prime}(\R)} \span \{\phi(2^{\kappa} \cdot -\nu),
\nu \in \Z\} \biggr\} \\ =
\biggl\{ f \in L_p(\R): \int_\R f(x)
\overline {g(x)} dx =0 \ \forall g \in
\span \{\phi(2^{\kappa} \cdot -\nu),  \nu \in \Z\} \biggr\} \\ =
\Ker E_\kappa^p.
\end{multline*}

Перейдем к проверке (2.2.12).
Пусть $ \kappa, \kappa^\prime \in \Z_+ $ и $ \kappa = \kappa^\prime. $
Тогда в силу (2.2.10) получаем
\begin{multline*}
(\mathcal E_\kappa^p)^2 = (E_\kappa^p)^2 - E_{\kappa
-1}^p E_\kappa^p -E_\kappa^p E_{\kappa -1}^p
+(E_{\kappa -1}^p)^2 =\\
= E_\kappa^p - E_{\kappa -1}^p -E_{\kappa -1}^p
+E_{\kappa -1}^p = E_\kappa^p - E_{\kappa -1}^p =
\mathcal E_\kappa^p.
\end{multline*}

Пусть $ \kappa \ne \kappa^\prime. $ Предположим, что $
\kappa^\prime < \kappa. $ Тогда, снова используя (2.2.10), выводим
\begin{multline*}
\mathcal E_{\kappa^\prime}^p \mathcal E_\kappa^p =
E_{\kappa^\prime}^p E_\kappa^p - E_{\kappa^\prime
-1}^p E_\kappa^p - E_{\kappa^\prime}^p E_{\kappa
-1}^p +
E_{\kappa^\prime -1}^p E_{\kappa -1}^p =\\
= E_{\kappa^\prime}^p - E_{\kappa^\prime -1}^p -
E_{\kappa^\prime}^p + E_{\kappa^\prime -1}^p =0.
\end{multline*}

Аналогично проверяется (2.2.12) при
$ \kappa^\prime > \kappa. $

При проверке (2.2.13) для $ f \in L_p(\R), g \in L_{p^\prime}(\R), $
применяя (2.2.3), (2.2.12), находим
$$
\int_{\R} (\mathcal E_\kappa^p f) \cdot \overline {(\mathcal
E_{\kappa^\prime}^{p^\prime} g)} dx =
\int_{\R} ( \mathcal E_{\kappa^\prime}^p \mathcal E_\kappa^p f) \cdot
\overline g dx  =
\int_\R 0 \cdot \overline g dx =0, \kappa, \kappa^\prime \in \Z_+:
\kappa \ne \kappa^\prime.
$$

Наконец, сопоставляя п. 2) предложения 2.2.1 и соотношения (2.2.11) с (1.2.3),
в соответствии с (1.2.4) получаем равенства п. 4). $ \square $

Теорема 2.2.3

Пусть выполнены условия предложения 2.2.1, а также функция $ \tilde \Phi $
помимо ранее отмеченных свойств обладает тем свойством, что
\begin{equation*} \tag{2.2.14}
\int_{\xi \in \R: | \xi| \ge s}
\tilde \Phi(\xi) d\xi \in L_1((1, \infty)).
\end{equation*}
Тогда при $ 1 < p < \infty $ в $ L_p(\R) $ справедливо равенство
\begin{equation*} \tag{2.2.15}
f = \sum_{ \kappa \in \Z_+} \mathcal E_\kappa^p f, f \in L_p(\R).
\end{equation*}

Доказательство.

Для доказательства (2.2.15) с учетом предложения 1.1.1 покажем, что
при $ 1 < p < \infty $ для $ f \in L_p(\R) $ имеет место соотношение
\begin{equation*} \tag{2.2.16}
\| f -E_\kappa^p f \|_{L_p(\R)} \to 0 \text{ при } \kappa \to \infty.
\end{equation*}

Справедливость (2.2.16) при $ p =2 $ уже установлена (см. (1.2.35), (1.2.26)).
Теперь убедимся в справедливости (2.2.16) при $ 1 < p < \infty. $
Покажем сначала, что (2.2.16) имеет место для любой функции
$ f \in L_\infty(\R), $ имеющей компактный носитель. В этой ситуации,
когда $ f \in L_\infty(\R) $ имеет компактный носитель $ \supp f, $
при $ \kappa \in \Z_+, R \in \R_+: R > \sup_{y \in \supp f} | y|, $ имеем
\begin{multline*} \tag{2.2.17}
\| f -E_\kappa^p f \|_{L_p(\R)}^p =
\int_{\R} | f(x) -(E_\kappa^p f)(x) |^p dx \\ =
\int_{x \in \R: | x| \le R} | f(x) -(E_\kappa^p f)(x) |^p dx +
\int_{x \in \R: | x| > R} | f(x) -(E_\kappa^p f)(x) |^p dx \\ =
\int_{R B^1} | f(x) -(E_\kappa^p f)(x) |^p dx +
\int_{x \in \R: | x| > R} | (E_\kappa^p f)(x) |^p dx.
\end{multline*}

Оценивая второе слагаемое в правой части (2.2.17), для $ f \in L_\infty(\R):
\supp f $ -- компакт, при $ \kappa \in \Z_+, R \in \R_+:
R > \sup_{y \in \supp f} | y|, $ с учетом (1.2.22), (1.2.23)
при $ p = \infty, $ получаем, что
\begin{multline*} \tag{2.2.18}
\int_{x \in \R: | x| > R} | (E_\kappa^p f)(x) |^p dx \le
\int_{x \in \R: | x| > R} \| E_\kappa^p f \|_{L_\infty(\R)}^{p -1} \times
| (E_\kappa^p f)(x) | dx \\ =
\| E_\kappa^\infty f \|_{L_\infty(\R)}^{p -1} \times
\int_{x \in \R: | x| > R} \biggl| \sum_{\nu_\kappa \in \Z}
2^{\kappa} \times \biggl(\int_{\R} f(y) \overline {\tilde \phi(2^\kappa y -\nu_\kappa)} dy\biggr) \times
\phi(2^\kappa x -\nu_\kappa)\biggr| dx\\ \le
(c_1 \| f \|_{L_\infty(\R)})^{p -1} \times
\int_{x \in \R: | x| > R} \biggl| \sum_{\nu_\kappa \in \Z}
2^{\kappa} \times \biggl(\int_{\R} f(y) \overline {\tilde \phi(2^\kappa y -\nu_\kappa)} dy\biggr) \times
\phi(2^\kappa x -\nu_\kappa)\biggl| dx.
\end{multline*}

Чтобы оценить правую часть (2.2.18), для $ f \in L_\infty(\R): \supp f $ --
компакт, при $ \kappa \in \Z_+, R \in \R_+: R > \sup_{y \in \supp f} | y|, $
почти для всех $ x \in \R: | x| > R, $ находим, что
\begin{multline*} \tag{2.2.19}
\biggl| \sum_{\nu_\kappa \in \Z} 2^{\kappa} \times
\biggl(\int_{\R} f(y) \overline {\tilde \phi(2^\kappa y -\nu_\kappa)} dy\biggr) \times
\phi(2^\kappa x -\nu_\kappa)\biggr| \\ =
\biggl| \sum_{\nu_\kappa \in \Z: | x -2^{-\kappa} \nu_\kappa| \le | x|/4} 2^{\kappa} \times
\biggl(\int_{\R} f(y) \overline {\tilde \phi(2^\kappa y -\nu_\kappa)} dy\biggr) \times
\phi(2^\kappa x -\nu_\kappa) \\+
\sum_{\nu_\kappa \in \Z: | x -2^{-\kappa} \nu_\kappa| > | x|/4} 2^{\kappa} \times
\biggl(\int_{\R} f(y) \overline {\tilde \phi(2^\kappa y -\nu_\kappa)} dy\biggr) \times
\phi(2^\kappa x -\nu_\kappa)\biggr|\\ \le
\biggl| \sum_{\nu_\kappa \in \Z: | x -2^{-\kappa} \nu_\kappa| \le | x|/4} 2^{\kappa} \times
\biggl(\int_{\R} f(y) \overline {\tilde \phi(2^\kappa y -\nu_\kappa)} dy\biggr) \times
\phi(2^\kappa x -\nu_\kappa)\biggr| \\+
\biggl| \sum_{\nu_\kappa \in \Z: | x -2^{-\kappa} \nu_\kappa| > | x|/4} 2^{\kappa} \times
\biggl(\int_{\R} f(y) \overline {\tilde \phi(2^\kappa y -\nu_\kappa)} dy\biggr) \times
\phi(2^\kappa x -\nu_\kappa)\biggr|\\ \le
\sum_{\nu_\kappa \in \Z: | x -2^{-\kappa} \nu_\kappa| \le | x|/4} 2^{\kappa} \times
\biggl| \int_{\R} f(y) \overline {\tilde \phi(2^\kappa y -\nu_\kappa)} dy\biggr| \times
| \phi(2^\kappa x -\nu_\kappa)| + \\
\sum_{\nu_\kappa \in \Z: | x -2^{-\kappa} \nu_\kappa| > | x|/4} 2^{\kappa} \times
\biggl| \int_{\R} f(y) \overline {\tilde \phi(2^\kappa y -\nu_\kappa)} dy\biggr| \times
| \phi(2^\kappa x -\nu_\kappa)|.
\end{multline*}
Замечая, что при $ \kappa \in \Z_+, \nu_\kappa \in \Z $ почти для всех $ x \in \R $
имеет место неравенство
\begin{equation*}
| \phi(2^\kappa x -\nu_\kappa) | \le \|  \phi \|_{L_\infty(\R)},
\end{equation*}
а для $ f \in L_\infty(\R): \supp f $ есть компакт, при $ \kappa \in \Z_+,
\nu_\kappa \in \Z, $ благодаря (2.2.1) (см. также (1.2.24)), выполняется
неравенство
\begin{multline*}
\biggl| \int_{\R} f(y) \overline {\tilde \phi(2^\kappa y -\nu_\kappa)} dy\biggr| =
\biggl| \int_{\supp f} f(y) \overline {\tilde \phi(2^\kappa y -\nu_\kappa)} dy\biggr|\\ \le
\| f \|_{L_\infty(\R)} \times \int_{\supp f} | \tilde \phi(2^\kappa y -\nu_\kappa)| dy\\ \le
\| f \|_{L_\infty(\R)} \times \int_{\supp f} 2^\kappa \times
\int_{2^{-\kappa}(\nu_\kappa +(1/2) B^1)}
\tilde \Phi(2^\kappa (y -u)) du dy \\ =
2^\kappa \times \| f \|_{L_\infty(\R)} \times
\int_{\supp f} \int_{2^{-\kappa}(\nu_\kappa +(1/2) B^1)}
\tilde \Phi(2^\kappa (y -u)) du dy,
\end{multline*}
и учитывая, что при $ \kappa \in \Z_+, R \in \R_+: R \ge 4 \sup_{z \in \supp f} | z|, R > 2,
x \in \R: | x| > R, \nu_\kappa \in \Z:
| x -2^{-\kappa} \nu_\kappa| \le | x|/4,
y \in \supp f, u \in 2^{-\kappa}(\nu_\kappa +(1/2) B^1), $ соблюдается неравенство
\begin{multline*}
| y -u| = | -x +y +x -2^{-\kappa} \nu_\kappa +2^{-\kappa} \nu_\kappa -u| \ge
| x| -| y +x -2^{-\kappa} \nu_\kappa +2^{-\kappa} \nu_\kappa -u| \\ \ge
| x| -| y| -| x -2^{-\kappa} \nu_\kappa| -| 2^{-\kappa} \nu_\kappa -u| \ge
| x| -\left(\sup_{z \in \supp f} | z|\right) -| x|/4 -(1/2) 2^{-\kappa}\\ \ge | x| -R/4 -| x|/4
-R/4 > | x| -| x|/4 -| x|/4 -| x|/4 = | x|/4,
\end{multline*}
для первой суммы в правой части (2.2.19), имея в виду (2.2.14), получаем
соотношение
\begin{multline*} \tag{2.2.20}
\sum_{\nu_\kappa \in \Z: | x -2^{-\kappa} \nu_\kappa| \le | x|/4} 2^{\kappa} \times
\biggl| \int_{\R} f(y) \overline {\tilde \phi(2^\kappa y -\nu_\kappa)} dy| \times
| \phi(2^\kappa x -\nu_\kappa)\biggr|\\ \le
\sum_{\nu_\kappa \in \Z: | x -2^{-\kappa} \nu_\kappa| \le | x|/4}
2^{\kappa} \times 2^\kappa \times \| f \|_{L_\infty(\R)}\\ \times
\biggl(\int_{\supp f} \int_{2^{-\kappa}(\nu_\kappa +(1/2) B^1)}
\tilde \Phi(2^\kappa (y -u)) du dy\biggr) \times
\| \phi \|_{L_\infty(\R)} \\ =
c_2(\phi) \| f \|_{L_\infty(\R)} 2^{2 \kappa} \times
\sum_{\nu_\kappa \in \Z: | x -2^{-\kappa} \nu_\kappa| \le | x|/4}
\int_{\supp f} \int_{2^{-\kappa}(\nu_\kappa +(1/2) B^1)}
\tilde \Phi(2^\kappa (y -u)) du dy \\ =
c_2 \| f \|_{L_\infty(\R)} 2^{2 \kappa} \times
\int_{\supp f} \biggl(\sum_{\nu_\kappa \in \Z: | x -2^{-\kappa} \nu_\kappa| \le | x|/4}
\int_{2^{-\kappa}(\nu_\kappa +(1/2) B^1)}
\tilde \Phi(2^\kappa (y -u)) du\biggr) dy \\ =
c_2 \| f \|_{L_\infty(\R)} 2^{2 \kappa} \times
\int_{\supp f} \int_{\cup_{\nu_\kappa \in \Z: | x -2^{-\kappa} \nu_\kappa| \le | x|/4}
(2^{-\kappa}(\nu_\kappa +(1/2) B^1))}
\tilde \Phi(2^\kappa (y -u)) du dy \\ \le
c_2 \| f \|_{L_\infty(\R)} 2^{2 \kappa} \times
\int_{\supp f} \int_{u \in \R: | y -u| \ge | x|/4}
\tilde \Phi(2^\kappa (y -u)) du dy \\ =
c_2 \| f \|_{L_\infty(\R)} 2^{2 \kappa} \times
\int_{\supp f} \int_{z \in \R: | z| \ge | x|/4}
\tilde \Phi(2^\kappa z) dz dy \\ =
c_2 \| f \|_{L_\infty(\R)} (\mes \supp f)
2^{2 \kappa} \times \int_{z \in \R: | z| \ge | x|/4} \tilde \Phi(2^\kappa z) dz \\ =
c_2 \| f \|_{L_\infty(\R)} (\mes \supp f)
2^{2 \kappa} \times \int_{\xi \in \R: | \xi| \ge 2^\kappa | x|/4}
\tilde \Phi(\xi) 2^{-\kappa} d\xi \\ =
c_2 \| f \|_{L_\infty(\R)} (\mes \supp f)
2^{ \kappa} \times \int_{\xi \in \R: | \xi| \ge 2^{\kappa -2} | x|}
\tilde \Phi(\xi) d\xi,\\
\text{ почти для всех } x \in \R: | x| > R,
R \in \R_+: R > 4 \sup_{z \in \supp f} | z|, R > 2, \kappa \in \Z_+.
\end{multline*}

Кроме того, при оценке второй суммы в правой части (2.2.19), учитывая, что
при $ \kappa \in \Z_+, \nu_\kappa \in \Z $ соблюдается неравенство
\begin{multline*}
\biggl| \int_{\R} f(y) \overline {\tilde \phi(2^\kappa y -\nu_\kappa)} dy\biggr| \le
\| f \|_{L_\infty(\R)} \int_{\R} | \tilde \phi(2^\kappa y -\nu_\kappa)| dy \\ =
\| f \|_{L_\infty(\R)} \int_{\R} | \tilde \phi(z)| 2^{-\kappa} dz =
c_3(\tilde \phi) \| f \|_{L_\infty(\R)} 2^{-\kappa},
\end{multline*}
используя (1.2.24), а также, принимая во внимание, что при $ \kappa \in \Z_+, R \in \R_+: R > 4,
x \in \R: | x| > R, $ для $ \nu_\kappa \in \Z: | x -2^{-\kappa} \nu_\kappa| \ge
| x|/4, u \in (2^{-\kappa}(\nu_\kappa +(1/2) B^1)), $ выполняется неравенство
\begin{multline*}
| x -u| = | x -2^{-\kappa} \nu_\kappa  +2^{-\kappa} \nu_\kappa -u| \ge
| x -2^{-\kappa} \nu_\kappa| -| 2^{-\kappa} \nu_\kappa -u| \ge | x|/4 -2^{-\kappa -1}\\
\ge | x|/4 -1/2 \ge | x|/4 -R/8 > | x|/4 -| x|/8 = | x|/8,
\end{multline*}
заключаем, что при $ \kappa \in \Z_+, R \in \R_+: R > 4, $
почти для всех $ x \in \R: | x| > R, $
имеет место неравенство
\begin{multline*} \tag{2.2.21}
\sum_{\nu_\kappa \in \Z: | x -2^{-\kappa} \nu_\kappa| > | x|/4} 2^{\kappa} \times
\biggl| \int_{\R} f(y) \overline {\tilde \phi(2^\kappa y -\nu_\kappa)} dy\biggr| \times
| \phi(2^\kappa x -\nu_\kappa)| \\ \le
\sum_{\nu_\kappa \in \Z: | x -2^{-\kappa} \nu_\kappa| > | x|/4} 2^{\kappa} \times
c_3 \| f \|_{L_\infty(\R)} 2^{-\kappa} \times
2^\kappa \int_{2^{-\kappa}(\nu_\kappa +(1/2) B^1)}
\Phi(2^\kappa (x -u)) du \\ =
c_3 \| f \|_{L_\infty(\R)} 2^\kappa \times
\sum_{\nu_\kappa \in \Z: | x -2^{-\kappa} \nu_\kappa| > | x|/4}
\int_{2^{-\kappa}(\nu_\kappa +(1/2) B^1)}
\Phi(2^\kappa (x -u)) du \\ =
c_3 \| f \|_{L_\infty(\R)} 2^\kappa \times
\int_{\cup_{\nu_\kappa \in \Z: | x -2^{-\kappa} \nu_\kappa| > | x|/4}
(2^{-\kappa}(\nu_\kappa +(1/2) B^1))}
\Phi(2^\kappa (x -u)) du \\ \le
c_3 \| f \|_{L_\infty(\R)} 2^\kappa \times
\int_{u \in \R: | x -u| > | x|/8} \Phi(2^\kappa (x -u)) du \\ =
c_3 \| f \|_{L_\infty(\R)} 2^\kappa \times
\int_{z \in \R: | z| > | x|/8} \Phi(2^\kappa z) dz \\ =
c_3 \| f \|_{L_\infty(\R)} 2^\kappa \times
\int_{\xi \in \R: | \xi| > 2^\kappa | x|/8} \Phi(\xi) 2^{-\kappa} d\xi \\ =
c_3 \| f \|_{L_\infty(\R)} \times
\int_{\xi \in \R: | \xi| > 2^{\kappa -3} | x|} \Phi(\xi) d\xi.
\end{multline*}

Соединяя (2.2.18). (2.2.19), (2.2.20),  (2.2.21), получаем, что
для $ f \in L_\infty(\R): \supp f $ -- компакт, при $ R \in \R_+:
R \ge \max(4, 4 \sup_{z \in \supp f}| z|), \kappa \in \Z_+ $ в виду (2.2.14),
(2.1.1) справедливо соотношение
\begin{multline*} \tag{2.2.22}
\int_{x \in \R: | x| > R} | (E_\kappa^p f)(x) |^p dx \le
(c_1 \| f \|_{L_\infty(\R)})^{p -1}\\ \times
\int_{x \in \R: | x| > R} \biggl(c_2 \| f \|_{L_\infty(\R)} (\mes \supp f)
2^{ \kappa}\\
 \times \int_{\xi \in \R: | \xi| \ge 2^{\kappa -2} | x|}
\tilde \Phi(\xi) d\xi +
c_3 \| f \|_{L_\infty(\R)} \times
\int_{\xi \in \R: | \xi| > 2^{\kappa -3} | x|} \Phi(\xi) d\xi\biggr) dx \\ =
\| f \|_{L_\infty(\R)}^{p} \times
\biggl(c_4 (\mes \supp f) \int_{x \in \R: | x| > R} 2^{\kappa} \\
\times \int_{\xi \in \R: | \xi| \ge 2^{\kappa -2} | x|}
\tilde \Phi(\xi) d\xi dx + c_5 \int_{x \in \R: | x| > R}
\int_{\xi \in \R: | \xi| > 2^{\kappa -3} | x|} \Phi(\xi) d\xi dx\biggr) \\ =
\| f \|_{L_\infty(\R)}^{p} \times
\biggl(c_4 (\mes \supp f) \int_{t \in \R: | t| > 2^{\kappa -2} R} 2^{\kappa}\\
 \times \int_{\xi \in \R: | \xi| \ge | t|}
\tilde \Phi(\xi) d\xi 2^{-\kappa +2} dt +
c_5 \int_{t \in \R: | t| > 2^{\kappa -3} R}
\int_{\xi \in \R: | \xi| > | t|} \Phi(\xi) d\xi 2^{-\kappa +3} dt\biggr) \\ =
\| f \|_{L_\infty(\R)}^{p} \times
\biggl(c_6 (\mes \supp f) \int_{t \in \R: | t| > 2^{\kappa -2} R}
\int_{\xi \in \R: | \xi| \ge | t|}
\tilde \Phi(\xi) d\xi dt \\ +
c_{7} 2^{-\kappa} \int_{t \in \R: | t| > 2^{\kappa -3} R}
\int_{\xi \in \R: | \xi| > | t|} \Phi(\xi) d\xi dt\biggr) \to 0 \text{ при }
\kappa \to \infty.
\end{multline*}

Теперь при $ 1 < p < \infty $ для $ f \in L_\infty(\R): \supp f $ -- компакт,
фиксируем $ R \in \R_+: R \ge \max(4, 4 \sup_{z \in \supp f}| z|)$ и для произвольного
$ \epsilon > 0, $ исходя из (2.2.22), найдем $ \kappa_0 \in \Z_+, $ для
которого при $ \kappa \in \Z_+: \kappa \ge \kappa_0, $ выполняется неравенство
\begin{equation*} \tag{2.2.23}
\int_{x \in \R: | x| > R} | (E_\kappa^p f)(x) |^p dx < \epsilon^p.
\end{equation*}
Затем выбрав $ \delta > 0 $ такое, что $ 2 R \delta^p < \epsilon^p, $
используя (1.2.23) при $ p = \infty, $ а также с учетом (1.2.26) применяя
(1.2.35), получим
\begin{multline*} \tag{2.2.24}
\| f - E_\kappa^p f \|_{L_p(R B^1)}^p =
\int_{R B^1} | f -(E_\kappa^p f) |^p dx \\ =
\int_{x \in R B^1: | f(x) -(E_\kappa^p f)(x) | \le \delta}
| f -(E_\kappa^p f) |^p dx +\int_{x \in R B^1: | f(x) -(E_\kappa^p f)(x) | > \delta} | f -(E_\kappa^p f) |^p dx \\ \le
\delta^p \mes (R B^1) +\int_{x \in R B^1: | f(x) -(E_\kappa^p f)(x) | >
\delta} | f -(E_\kappa^p f) |^p dx \\ \le
\delta^p 2R +\| f -(E_\kappa^p f) \|_{L_\infty(\R)}^p
\mes \{x \in R B^1: | f(x) -(E_\kappa^p f)(x) | > \delta \} \\ \le
2R \delta^p +(\| f \|_{L_\infty(\R)}
+\| E_\kappa^p f \|_{L_\infty(\R)})^p
\mes \{x \in R B^1: | f(x) -(E_\kappa^p f)(x) | > \delta \} \\ \le
2R \delta^p +(\| f \|_{L_\infty(\R)}
+\| E_\kappa^p f \|_{L_\infty(\R)})^p
\delta^{-2} \int_{ \{x \in R B^1: | f(x) -(E_\kappa^p f)(x) | > \delta \}}
| f(x) -(E_\kappa^p f)(x) |^2 dx \\ \le
2R \delta^p +(c_{8} \| f \|_{L_\infty(\R)})^p
\delta^{-2} \int_\R | f(x) -(E_\kappa^p f)(x) |^2 dx \\ =
2R \delta^p +(c_{8} \| f \|_{L_\infty(\R)})^p
\delta^{-2} \| f -(U_\kappa f) \|_{L_2(\R)}^2\\ <
\epsilon^p +(c_{8} \| f \|_{L_\infty(\R)})^p
\delta^{-2} \| f -(U_\kappa f) \|_{L_2(\R)}^2 < 2 \epsilon^p
\end{multline*}
для $ \kappa \in \Z_+ $ таких, что
$$
\| f -(U_\kappa f) \|_{L_2(\R)}^2 <
\epsilon^p \delta^2 /(c_{8} \| f \|_{L_\infty(\R)})^p.
$$
Объединяя (2.2.17), (2.2.23), (2.2.24), видим, что при достаточно больших
$ \kappa \in \Z_+ $ соблюдается неравенство
$$
\| f -E_\kappa^p f \|_{L_p(\R)}^p < 3 \epsilon^p,
$$
т.е. (2.2.16) выполняется для любой функции $ f \in L_\infty(\R), $
имеющей компактный носитель. А поскольку множество таких функций
плотно в $ L_p(\R) $ при $ 1 < p < \infty, $ то отсюда в силу (2.1.5)
заключаем, что (2.2.16) верно для любой $ f \in L_p(\R). $

Равенство (2.2.15) в виду (1.1.1) является следствием
соотношения (2.2.16). $ \square $
\bigskip

2.3. В этом пункте дается описание свойств операторов проектирования
на подпространства всплесков, соответствующие неизотропному кратно-масштабному
анализу, порожденному тензорным произведением гладких достаточно
быстро убывающих на бесконечности функций, и других проекторов, которые
используются в п. 2.4. при доказательстве основных результатов работы.

Предложение 2.3.1

Пусть $ d \in \N $ и наборы функций
$$
\mathcal \phi = \{\phi_1, \ldots, \phi_d\},
\tilde{\mathcal \phi} = \{ \tilde \phi_1, \ldots, \tilde \phi_d \}
$$
таковы, что при $ j =1, \ldots, d $ соблюдаются условия предложения 2.2.1
с функциями $ \phi_j, \tilde \phi_j. $
При $ 1 < p < \infty, \kappa \in \Z_+^d $ рассмотрим операторы
$$
E_{\kappa_j}^{j, p} =
E_{\kappa_j}^{\phi_j, \tilde \phi_j, p}, j =1, \ldots, d, (см. предложение 1.2.2)
$$
и определим операторы
\begin{equation*} \tag{2.3.1}
E_\kappa^p = E_\kappa^{\mathcal \phi, \tilde{\mathcal \phi}, p} = \prod_{j=1}^d V_j^{L_p} (E_{\kappa_j}^{\phi_j,
\tilde \phi_j, p}) = \prod_{j=1}^d V_j (E_{\kappa_j}^{j, p}).
\end{equation*}

Тогда для линейных операторов
\begin{eqnarray*}
\mathcal E_\kappa^{p}&: &L_p(\R^d) \mapsto L_p(\R^d), \\
\mathcal E_{\kappa_j}^{j, p}&: &L_p(\R) \mapsto L_p(\R), j =1,\ldots,d,
\end{eqnarray*}
определяемых соотношениями
\begin{eqnarray*}
\mathcal E_\kappa^{p} &=& \sum_{\epsilon \in \Upsilon^d:
\s(\epsilon) \subset \s(\kappa)} (-\e)^\epsilon E_{\kappa -\epsilon}^{p} \\
\mathcal E_{\kappa_j}^{j, p} &=&
\mathcal E_{\kappa_j}^{\phi_j, \tilde \phi_j, p}, j =1,\ldots,d \text{ (см. определение после (1.2.28)) },
\end{eqnarray*}
имеет место равенство
\begin{equation*} \tag{2.3.2}
\mathcal E_\kappa^{p} = \prod_{j=1}^d V_j(\mathcal E_{\kappa_j}^{j, p}).
\end{equation*}

Доказательство.

Проверим справедливость (2.3.2). Используя (2.3.1),
(1.3.3) и п. 2 леммы 1.3.1, имеем
\begin{multline*}
\mathcal E_\kappa^{p} = \sum_{\epsilon \in \Upsilon^d:
\s(\epsilon) \subset \s(\kappa)} (-\e)^\epsilon E_{\kappa
-\epsilon}^{p} =
\sum_{\epsilon \in \Upsilon^d: \s(\epsilon)
\subset \s(\kappa)}
\biggl(\prod_{j=1}^d (-1)^{\epsilon_j}\biggr) \biggl(\prod_{j=1}^d V_j(E_{\kappa_j -\epsilon_j}^{j, p})\biggr) =\\
=\sum_{\epsilon \in \Upsilon^d: \s(\epsilon) \subset \s(\kappa)}
\biggl(\prod_{j=1,\ldots,d: j \notin \s(\kappa)} (-1)^{\epsilon_j}
V_j(E_{\kappa_j -\epsilon_j}^{j, p})\biggr) \biggl(\prod_{j \in \s(\kappa)}
(-1)^{\epsilon_j} V_j(E_{\kappa_j -\epsilon_j}^{j, p})\biggr) =\\
=\sum_{\epsilon \in \Upsilon^d: \s(\epsilon) \subset \s(\kappa)}
\biggl(\prod_{j=1,\ldots,d: j \notin \s(\kappa)} V_j(E_0^{j, p})\biggr)
\biggl(\prod_{j \in \s(\kappa)} (-1)^{\epsilon_j} V_j(E_{\kappa_j -\epsilon_j}^{j, p})\biggl) =\\
=\biggl(\prod_{j=1,\ldots,d: j \notin \s(\kappa)} V_j(E_0^{j, p})\biggr)
\sum_{\epsilon \in \Upsilon^d: \s(\epsilon) \subset \s(\kappa)}
\biggl(\prod_{j \in \s(\kappa)} (-1)^{\epsilon_j} V_j(E_{\kappa_j -\epsilon_j}^{j, p})\biggr) =\\
=\biggl(\prod_{j=1,\ldots,d: j \notin \s(\kappa)} V_j(E_0^{j, p})\biggr)
\sum_{\epsilon^{\s(\kappa)} \in (\Upsilon^d)^{\s(\kappa)}}
\biggl(\prod_{j \in \s(\kappa)} (-1)^{\epsilon_j} V_j(E_{\kappa_j -\epsilon_j}^{j, p})\biggr) =\\
=\biggl(\prod_{j=1,\ldots,d: j \notin \s(\kappa)} V_j(E_0^{j, p})\biggr)
\biggl(\prod_{j \in \s(\kappa)} (\sum_{\epsilon_j =0,1}
(-1)^{\epsilon_j} V_j(E_{\kappa_j -\epsilon_j}^{j, p}))\biggr) =\\
=\biggl(\prod_{j=1,\ldots,d: j \notin \s(\kappa)} V_j(E_0^{j, p})\biggr)
\biggl(\prod_{j \in \s(\kappa)} (V_j(E_{\kappa_j}^{j, p}) - V_j(E_{\kappa_j -1}^{j, p}))\biggr) =\\
=\biggl(\prod_{j=1,\ldots,d: j \notin \s(\kappa)} V_j(E_0^{j, p})\biggr)
\biggl(\prod_{j \in \s(\kappa)} V_j(E_{\kappa_j}^{j, p} - E_{\kappa_j
-1}^{j, p})\biggr). \square
\end{multline*}

Предложение 2.3.2

В условиях предложения 2.3.1 при $ \kappa \in \Z_+^d, 1 < p < \infty $ справедливы следующие утверждения:

1) для $ \kappa^\prime \in \Z_+^d: \kappa^\prime \le \kappa, $ имеют место
равенства
\begin{equation*} \tag{2.3.3}
E_{\kappa^\prime}^p E_\kappa^p =
E_\kappa^p E_{\kappa^\prime}^p = E_{\kappa^\prime}^p;
\end{equation*}

2) для $ f \in L_p(\R^d), g \in L_{p^\prime}(\R^d) $
выполняются равенства
\begin{equation*} \tag{2.3.4}
\int_{\R^d} (E_\kappa f) \cdot \overline g dx = \int_{\R^d} f \cdot
\overline {(E_\kappa g)} dx;
\end{equation*}
и
\begin{equation*} \tag{2.3.5}
\int_{\R^d} (\mathcal E_\kappa f) \cdot \overline g dx =
\int_{\R^d} f \cdot \overline {(\mathcal E_\kappa g)} dx;
\end{equation*}

3) для $ \kappa^\prime \in \Z_+^d $ соблюдаются равенства
\begin{equation*} \tag{2.3.6}
\mathcal E_\kappa^{p} \mathcal E_{\kappa^\prime}^{p} =
\begin{cases} \mathcal E_\kappa^{p},
  \text{при $ \kappa = \kappa^\prime $}; \\
       0,  \text{при $ \kappa \ne \kappa^\prime $}.
\end{cases}
\end{equation*}

4) оператор $ E_\kappa^p $ является проектором в пространстве
$ L_p(\R^d), $ для которого
\begin{multline*} \tag{2.3.7}
\Ker E_\kappa^p =
\{ f \in L_p(\R^d): \int_{\R^d} f(x) \cdot
\overline {\phi(2^\kappa x -\nu)} dx =0 \ \forall \nu \in \Z^d \} \\ =
\{ f \in L_p(\R^d): \int_{\R^d} f(x) \cdot
\overline {g(x)} dx =0 \ \forall g \in
\span \{\phi(2^\kappa \cdot -\nu),  \nu \in \Z^d\} \};
\end{multline*}
\begin{equation*} \tag{2.3.8}
\Im E_\kappa^p = \close_{L_p(\R^d)} (\span \{ \phi(2^\kappa \cdot -\nu),
\nu \in \Z^d \}),
\end{equation*}
где $ \phi: \R^d \mapsto \C $ -- функция, задаваемая равенством
$ \phi(x_1, \ldots, x_d) = \prod_{j =1}^d \phi_j(x_j). $

Доказательство.

Проверяя первое утверждение предложения, на основании (2.3.1), (1.3.3),
п. 2 лемммы 1.3.1 и (2.2.10) имеем
\begin{multline*}
E_{\kappa^\prime}^{p} E_\kappa^{p} = \biggl(\prod_{j=1}^d V_j (E_{\kappa_j^\prime}^{j, p})\biggr)
\biggl(\prod_{j=1}^d V_j (E_{\kappa_j}^{j, p})\biggr) =
\prod_{j=1}^d (V_j (E_{\kappa_j^\prime}^{j, p})  V_j (E_{\kappa_j}^{j, p})) =\\
\prod_{j=1}^d V_j (E_{\kappa_j^\prime}^{j, p}  E_{\kappa_j}^{j, p}) =
\prod_{j=1}^d V_j (E_{\kappa_j^\prime}^{j, p}) = E_{\kappa^\prime}^{p},
\end{multline*}
аналогично проверяется и второе равенство в (2.3.3).

Для доказательства (2.3.4) сначала установим, что при $ \kappa \in \Z_+^d $
для $ f \in L_p(\R^d), g \in L_{p^\prime}(\R^d) $ для любого непустого
множества $ J \subset \{1,\ldots,d\} $ выполняется равенство
\begin{equation*} \tag{2.3.9}
\int_{\R^d} \biggl((\prod_{j \in J} V_j(E_{\kappa_j}^{j, p}))f\biggr) \overline g dx =
\int_{\R^d} f \cdot \overline {\biggl((\prod_{j \in J}
V_j(E_{\kappa_j}^{j, p^\prime}))g\biggr)} dx.
\end{equation*}

Вывод (2.3.9) проведем по индукции относительно $ \card J. $
При $ \card J =1, $ т.е. при $ J = \{j\}, j =1, \ldots,d, $ в силу теоремы
Фубини, (1.3.1), (2.2.2) имеем
\begin{multline*}
 \int_{\R^d} ((V_j(E_{\kappa_j}^{j, p}))f)(x)\overline {g(x)} dx \\ =
 \int_{\R^{d -1}} \int_\R  ((V_j(E_{\kappa_j}^{j, p}))f)(x_1,\ldots, x_{j -1},
x_j, x_{j +1}, \ldots, x_d)\cdot\\ \overline {g(x_1,\ldots, x_{j -1}, x_j,
x_{j +1}, \ldots, x_d)} dx_j dx_1 \ldots dx_{j -1} dx_{j +1} \ldots dx_d \\ =
 \int_{\R^{d -1}} \int_\R  ( E_{\kappa_j}^{j, p} f(x_1,\ldots, x_{j -1},
\cdot, x_{j +1}, \ldots, x_d))(x_j)\cdot\\ \overline {g(x_1,\ldots, x_{j -1}, x_j,
x_{j +1}, \ldots, x_d)} dx_j dx_1 \ldots dx_{j -1} dx_{j +1} \ldots dx_d \\ =
 \int_{\R^{d -1}} \int_\R f(x_1,\ldots, x_{j -1}, x_j, x_{j +1}, \ldots, x_d)\cdot\\
\overline {( E_{\kappa_j}^{j, p^\prime} g(x_1,\ldots, x_{j -1},\cdot, x_{j +1},
\ldots, x_d))(x_j) } dx_j dx_1 \ldots dx_{j -1} dx_{j +1} \ldots dx_d \\ =
 \int_{\R^{d -1}} \int_\R f(x_1,\ldots, x_{j -1}, x_j, x_{j +1}, \ldots, x_d)\cdot\\
\overline {( V_j( E_{\kappa_j}^{j, p^\prime}) g)(x_1,\ldots, x_{j -1}, x_j,
x_{j +1}, \ldots, x_d) } dx_j dx_1 \ldots dx_{j -1} dx_{j +1} \ldots dx_d \\ =
 \int_{\R^d} f(x) \overline{(V_j( E_{\kappa_j}^{j, p^\prime}) g)(x)} dx,
\end{multline*}
что совпадает с (2.3.9) при $ J = \{j\}, j =1,\ldots,d. $

Предположим, что (2.3.9) справедливо для любого множества $ J \subset
\{1,\ldots,d\}, $  для которого $ \card J \le m \le d -1. $ Покажем, что тогда
(2.3.9) имеет место для любого множества $ \mathcal J \subset
\{1,\ldots,d\}, $ для которого $ \card \mathcal J = m +1. $
Представляя $ \mathcal J $  в виде $ \mathcal J = J \cup \{i\}, $
где $ i \notin J, $ и пользуясь (1.3.3), предположением индукции, получаем
\begin{multline*}
\int_{\R^d} \biggl(\biggl(\prod_{j \in \mathcal J} V_j(E_{\kappa_j}^{j, p})\biggr)f\biggr) \overline g dx =
\int_{\R^d} (V_i(E_{\kappa_i}^{i, p}))
\biggl(\biggl(\prod_{j \in J} V_j(E_{\kappa_j}^{j, p})\biggr)f\biggr) \overline g dx \\ =
\int_{\R^d} \biggl(\biggl(\prod_{j \in J} V_j(E_{\kappa_j}^{j, p})\biggr)f\biggr)
\overline {((V_i(E_{\kappa_i}^{i, p^\prime})) g) } dx \\ =
\int_{\R^d} f \cdot \overline {\biggl(\biggl(\prod_{j \in J}
V_j(E_{\kappa_j}^{j, p^\prime})\biggr) ((V_i(E_{\kappa_i}^{i, p^\prime})) g)\biggr) } dx =
\int_{\R^d} f \cdot \overline {\biggl(\biggl(\prod_{j \in \mathcal J}
V_j(E_{\kappa_j}^{j, p^\prime})\biggr)g\biggr)} dx,
\end{multline*}
что завершает вывод (2.3.9).

Полагая в (2.3.9) $ J = \{1,\ldots,d\}, $ и учитывая (2.3.1),
приходим к (2.3.4).

Проверяя (2.3.5), для $ f \in L_p(\R^d), g \in L_{p^\prime}(\R^d) $
ввиду (2.3.4) получаем
\begin{multline*}
\int_{\R^d} (\mathcal E_\kappa f) \overline {g} dx = \int_{\R^d}
\biggl(\sum_{\epsilon \in \Upsilon^d: \s(\epsilon) \subset \s(\kappa)}
(-\e)^\epsilon E_{\kappa -\epsilon} f\biggr) \overline {g} dx\\
=\int_{\R^d} \sum_{\epsilon \in \Upsilon^d:\s(\epsilon) \subset
\s(\kappa)} (-\e)^\epsilon (E_{\kappa -\epsilon} f) \overline {g} dx =
\sum_{\epsilon \in \Upsilon^d: \s(\epsilon) \subset \s(\kappa)}
(-\e)^\epsilon \int_{\R^d} (E_{\kappa -\epsilon} f) \overline {g} dx\\
=\sum_{\epsilon \in \Upsilon^d: \s(\epsilon) \subset \s(\kappa)}
(-\e)^\epsilon \int_{\R^d} f \cdot \overline
{(E_{\kappa -\epsilon} g)} dx
= \int_{\R^d} \sum_{\epsilon \in \Upsilon^d: \s(\epsilon) \subset
\s(\kappa)}
(-\e)^\epsilon f \cdot \overline {(E_{\kappa -\epsilon} g)} dx \\
=\int_{\R^d} f \biggl(\sum_{\epsilon \in \Upsilon^d:\s(\epsilon)
\subset \s(\kappa)} (-\e)^\epsilon \overline
{(E_{\kappa -\epsilon} g)}\biggr) dx \\
=\int_{\R^d} f \overline {\biggl(\sum_{\epsilon \in \Upsilon^d: \s(\epsilon)
\subset \s(\kappa)} (-\e)^\epsilon
(E_{\kappa -\epsilon} g)\biggr)} dx =
\int_{\R^d} f \overline {(\mathcal E_\kappa g)} dx.
\end{multline*}

Наконец, убедимся в справедливости (2.3.6). При $ d =1 $ соотношение
(2.3.6) совпадает с (2.2.12).
Установим соблюдение (2.3.6) при произвольном $ d \in \N. $
Согласно (2.3.2), (1.3.3), п. 2 леммы 1.3.1 имеем
$$
\mathcal E_\kappa^{p} \mathcal E_{\kappa^\prime}^{p} =
\biggl(\prod_{j=1}^d V_j(\mathcal E_{\kappa_j}^{j, p})\biggr)
\biggl(\prod_{j=1}^d V_j(\mathcal E_{\kappa_j^\prime}^{j, p})\biggr) =
\prod_{j=1}^d (V_j(\mathcal E_{\kappa_j}^{j, p}) V_j(\mathcal
E_{\kappa_j^\prime}^{j, p})) = \prod_{j=1}^d V_j(\mathcal
E_{\kappa_j}^{j, p} \mathcal E_{\kappa_j^\prime}^{j, p}).
$$
Отсюда в силу равенства (2.2.12) с учетом (2.3.2) вытекает справедливость
(2.3.6) при произвольном $ d \in \N. $

В заключение установим справедливость утверждения п. 4). Тот факт, что
$ E_\kappa^p $ является проектором в $ L_p(\R^d) $ следует из (2.3.3).
Проверим выполнение (2.3.7). Для $ f \in \Ker E_\kappa^p, $ учитывая, что в силу
(2.3.1), (1.3.4), (2.2.4) для $ \kappa \in \Z_+^d, \nu \in \Z^d $ справедливы
равенства
\begin{multline*} \tag{2.3.10}
(E_\kappa^p \phi(2^\kappa \cdot -\nu))(x_1, \ldots, x_d) =
\biggl(\biggl(\prod_{j =1}^d V_j(E_{\kappa_j}^{j, p})\biggr)\biggl(\prod_{j =1}^d \phi_j(2^{\kappa_j}
y_j -\nu_j)\biggr)\biggr)(x_1, \ldots, x_d) \\ =
\prod_{j =1}^d (E_{\kappa_j}^{j, p} \phi_j(2^{\kappa_j} \cdot -\nu_j))(x_j) =
\prod_{j =1}^d \phi_j(2^{\kappa_j} x_j -\nu_j) = \phi(2^{\kappa_1} x_1 -\nu_1,
\ldots, 2^{\kappa_d} x_d -\nu_d),
\end{multline*}
а также используя (2.3.4), имеем
\begin{multline*}
0 = \int_{\R^d} 0 \cdot \overline{ \phi(2^\kappa x -\nu) } dx =
\int_{\R^d} (e_\kappa^p f)(x) \cdot \overline{ \phi(2^\kappa x -\nu) } dx \\ =
\int_{\R^d} f(x) \cdot \overline{(e_\kappa^{p^\prime} \phi(2^\kappa \cdot
-\nu))(x) } dx = \int_{\R^d} f(x) \cdot \overline{ \phi(2^\kappa x -\nu) } dx,
\nu \in \Z^d, \kappa \in \Z_+^d,
\end{multline*}
т.е. имеет место включение
\begin{equation*}
\Ker E_\kappa^p \subset \{ f \in L_p(\R^d):
\int_{\R^d} f(x) \cdot \overline{ \phi(2^\kappa x -\nu) } dx =0 \ \forall
\nu \in \Z^d \}.
\end{equation*}

Покажем, что соблюдается и обратное включение
\begin{equation*} \tag{2.3.11}
\{ f \in L_p(\R^d):
\int_{\R^d} f(x) \cdot \overline{ \phi(2^\kappa x -\nu) } dx =0 \ \forall
\nu \in \Z^d \} \subset \Ker E_\kappa^p.
\end{equation*}

Для этого убедимся в том, что справедлива следующая

Лемма 2.3.3

При соблюдении условий предложения 2.3.2, если для $ f \in L_p(\R^d) $
выполняются равенства
\begin{equation*} \tag{2.3.12}
\int_{\R^d} f(x) \cdot \overline{ \phi(2^\kappa x -\nu) } dx =0, \nu \in \Z^d,
\end{equation*}
то для любого набора функций $ \{g_j \in L_{p^\prime}(\R), j =1, \ldots, d\}, $
соблюдается равенство
\begin{equation*} \tag{2.3.13}
\int_{\R^d} f(x) \cdot \biggl(\prod{j =1}^d \overline{(E_{\kappa_j}^{j, p^\prime}
g_j)(x_j)}\biggr) dx =0.
\end{equation*}

Доказательство.

Доказательство леммы проведем по индукции относительно размерности $ d. $
Если для $ f \in L_p(\R) $ имеют место равенства (2.3.12) при $ d =1, $ то
в силу (2.2.5) $ E_\kappa^p f =0, $ а благодаря (2.2.2) для $ g \in L_{p^\prime}(\R) $
соблюдается равенство
$$
\int_{\R} f \cdot \overline {(E_\kappa^{p^\prime} g) } dx =
\int_{\R} (E_\kappa^p f) \cdot \overline g dx = \int_\R 0 \cdot \overline g dx =0,
$$
т.е. утверждение леммы при $ d =1 $ справедливо.

Предположим, что утверждение леммы верно при $ d = n -1. $ Покажем, что
тогда оно верно и при $ d = n. $
Пусть для $ f \in L_p(\R^n) $ выполняются равенства (2.3.12) при $ d = n. $
Тогда для каждого $ \nu_n \in \Z $ при любых $ (\nu_1, \ldots, \nu_{n -1}) \in \Z^{n -1} $
по теореме Фубини имеем
\begin{multline*} \tag{2.3.14}
0 = \int_{\R^n} f(x_1, \ldots, x_{n -1}, x_n) \cdot \overline{\biggl(\prod_{j =1}^n
\phi_j(2^{\kappa_j} x_j -\nu_j)\biggr)} dx_1 \ldots dx_n \\ =
\int_{\R^{n -1} \int_\R} f(x_1, \ldots, x_{n -1}, x_n) \cdot \biggl(\prod_{j =1}^n
\overline{\phi_j(2^{\kappa_j} x_j -\nu_j)}\biggr) dx_n dx_1 \ldots dx_{n -1} \\ =
\int_{\R^{n -1} (\int_\R} f(x_1, \ldots, x_{n -1}, x_n) \cdot
\overline{\phi_n(2^{\kappa_n} x_n -\nu_n)} dx_n)\\ \times
\overline{\biggl(\prod_{j =1}^{n -1}\phi_j(2^{\kappa_j} x_j -\nu_j)\biggr)} dx_1 \ldots dx_{n -1}.
\end{multline*}
Учитывая, что $ \phi_n(2^{\kappa_n} \cdot -\nu_n) \in L_{p^\prime}(\R), \nu_n
\in \Z $ (см. (1.2.20)), а, следовательно, $ \phi_n(2^{\kappa_n} x_n -\nu_n) \in
L_{p^\prime}((R B^{n -1}) \times \R), R \in \R_+, \nu_n \in \Z, $ в виду
неравенства Гельдера заключаем, что функция
$ f(x_1, \ldots, x_{n -1}, x_n) \cdot
\overline{\phi_n(2^{\kappa_n} x_n -\nu_n)} $ суммируема на
$ (R B^{n -1}) \times \R, R \in \R_+, $
и по теореме Фубини функция $ \int_\R f(x_1, \ldots, x_{n -1}, x_n) \cdot
\overline{\phi_n(2^{\kappa_n} x_n -\nu_n)} dx_n $ суммируема (и, значит,
измерима) на $ (R B^{n -1}), R \in \R_+. $ Причем, в силу неравенства Гельдера
и теоремы Фубини почти для всех $ (x_1, \ldots, x_{n -1}) \in \R^{n -1} $ имеет
место соотношение
\begin{multline*}
\biggl| \int_\R f(x_1, \ldots, x_{n -1}, x_n) \cdot
\overline{\phi_n(2^{\kappa_n} x_n -\nu_n)} dx_n|^p \\ \le
\biggl( \int_\R | f(x_1, \ldots, x_{n -1}, x_n)|^p dx_n\biggr) \cdot
\biggl(\int_\R | \overline{\phi_n(2^{\kappa_n} x_n -\nu_n)}|^{p^\prime} dx_n\biggr)^{p /p^\prime} \\
\in L_1(\R^{n -1}), \nu_n \in \Z.
\end{multline*}
Таким образом, видим, что для каждого $ \nu_n \in \Z $ функция
$$
\int_\R f(x_1, \ldots, x_{n -1}, x_n) \cdot
\overline{\phi_n(2^{\kappa_n} x_n -\nu_n)} dx_n \in L_p(\R^{n -1}
$$
и для всех $ (\nu_1, \ldots, \nu_{n -1}) \in \Z^{n -1} $ выполняется (2.3.14).
Поэтому, исходя из предположения индукции, применяя с учетом неравенства
Гельдера теорему Фубини, приходим к выводу, что для каждого $ \nu_n \in \Z $
для любого набора функций $ \{g_j \in L_{p^\prime}(\R), j =1, \ldots, n -1\} $
справедливо равенство
\begin{multline*}
0 = \int_{\R^{n -1}} \biggl(\int_\R f(x_1, \ldots, x_{n -1}, x_n) \cdot
\overline{\phi_n(2^{\kappa_n} x_n -\nu_n)} dx_n\biggr)\\ \times
\biggl(\prod_{j =1}^{n -1} \overline{ (E_{\kappa_j}^{j, p^\prime} g_j)(x_j)}\biggr)
dx_1 \ldots dx_{n -1} \\ =
\int_\R \biggl(\int_{\R^{n -1}} f(x_1, \ldots, x_{n -1}, x_n) \cdot
\biggl(\prod_{j =1}^{n -1} \overline{ (E_{\kappa_j}^{j, p^\prime} g_j)(x_j)}\biggr) dx_1 \ldots
dx_{n -1}\biggr)\\
\times \overline{\phi_n(2^{\kappa_n} x_n -\nu_n)} dx_n,
\end{multline*}
причем, в силу неравенства Гельдера и теоремы Фубини функция
$$
\int_{\R^{n -1}} f(x_1, \ldots, x_{n -1}, x_n) \cdot
(\prod_{j =1}^{n -1} \overline{ (E_{\kappa_j}^{j, p^\prime} g_j)(x_j)}) dx_1 \ldots
dx_{n -1} \in L_p(\R)
$$
(см. аналогичную ситуацию выше). Отсюда, принимая во внимание справедливость леммы при $ d =1, $
получаем, что для любого набора функций $ \{g_j \in L_{p^\prime}(\R),
j =1, \ldots, n\} $ соблюдаетс равенство
\begin{multline*}
\int_\R \biggl(\int_{\R^{n -1}} f(x_1, \ldots, x_{n -1}, x_n) \cdot
\biggl(\prod_{j =1}^{n -1} \overline{ (E_{\kappa_j}^{j, p^\prime} g_j)(x_j)}\biggr) dx_1 \ldots
dx_{n -1}\biggr) \\
\times \overline{(E_{\kappa_n}^{n, p^\prime} g_n)(x_n)} dx_n =0,
\end{multline*}
или, вследствие теоремы Фубини,
$$
\int_{\R^{n}} f(x_1, \ldots, x_{n -1}, x_n) \cdot
\biggl(\prod_{j =1}^{n} \overline{ (E_{\kappa_j}^{j, p^\prime} g_j)(x_j)}\biggr) dx_1 \ldots
dx_{n -1} dx_n =0,
$$
т.е. имеет место (2.3.13) при $ d = n. \square $

Теперь проверим выполнение включения (2.3.11). Пусть для $ f \in L_p(\R^d) $
выполняются равенства (2.3.12). Тогда для любой функции $ g \in L_{p^\prime}(\R^d) $
вида
\begin{equation*}\tag{2.3.15}
g(x_1, \ldots, x_d) = \prod_{j =1}^d g_j(x_j),
\end{equation*}
почти для всех $(x_1, \ldots, x_d) \in \R^d, g_j \in L_{p^\prime}(\R), j =1, \ldots, d$,
в силу (2.3.4), (2.3.1), (1.3.4), (2.3.13) имеем
\begin{multline*}
\int_{\R^d} (E_\kappa^p f) \cdot \overline g dx =
\int_{\R^d} f \cdot \overline {(E_\kappa^{p^\prime} g)} dx \\ =
\int_{\R^{d}} f(x_1, \ldots, x_d) \cdot
\overline{\biggl(\biggl(\prod_{j =1}^{d} V_j(E_{\kappa_j}^{j, p^\prime})\biggr) \biggl(\prod_{j =1}^d
g_j(y_j)\biggr)\biggr)(x_1, \ldots, x_d)} dx_1 \ldots dx_d \\ =
\int_{\R^{d}} f(x_1, \ldots, x_d) \cdot
\biggl(\prod_{j =1}^{d} \overline{ (E_{\kappa_j}^{j, p^\prime} g_j)(x_j)}\biggr) dx_1 \ldots
dx_d =0,
\end{multline*}
или, короче,
\begin{equation*} \tag{2.3.16}
\int_{\R^d} (E_\kappa^p f) \cdot \overline g dx =0.
\end{equation*}
Поскольку линейная оболочка функций $ g \in L_{p^\prime}(\R^d) $ вида (2.3.15)
плотна в $ L_{p^\prime}(\R^d), $ то вследствие неравенства Гельдера равенство
(2.3.16) имеет место для любой функции $ g \in L_{p^\prime}(\R^d), $
а, следовательно, $ E_\kappa^p f =0. $
Тем самым, доказательсто (2.3.11), а вместе с ним и (2.3.7) завершено.

Осталось проверить соблюдение (2.3.8).
Сначала заметим, что в виду (2.3.10) множество
$$
\{ \phi(2^\kappa \cdot -\nu), \nu \in \Z^d \} \subset \Im E_\kappa^p,
$$
и благодаря линейности $ \Im E_\kappa^p, $ линейная оболочка
$$
\span \{ \phi(2^\kappa \cdot -\nu), \nu \in \Z^d \} \subset \Im E_\kappa^p.
$$
Из этого включения, учитывая, что в  силу непрерывности проектора
$ E_\kappa^p $ в $ L_p(\R^d) $ его образ $ \Im E_\kappa^p $ замкнут в $ L_p(\R^d), $
вытекает, что
$$
\close_{L_p(\R^d)} (\span \{ \phi(2^\kappa \cdot -\nu),
\nu \in \Z^d \}) \subset \Im E_\kappa^p.
$$

Предположим, что существует функция $ f_0 \in L_p(\R^d) $ такая, что
$$
E_\kappa^p f_0 \notin \close_{L_p(\R^d)} (\span \{ \phi(2^\kappa \cdot -\nu),
\nu \in \Z^d \}).
$$
Тогда вследствие теоремы Хана-Банаха существует непрерывный линейный
функционал на $ L_p(\R^d), $ аннулирующий замкнутое подпространство
$ \close_{L_p(\R^d)} (\span \{ \phi(2^\kappa \cdot -\nu), \nu \in \Z^d \}), $
значение которого на функции $ E_\kappa^p f_0 $ отлично от нуля.
Учитывая равенство $ (L_p(\R^d))^* = L_{p^\prime}(\R^d), $ получаем, что
существует функция $ g \in L_{p^\prime}(\R^d) $ такая, что для любого
$ \nu \in \Z^d $ выполняется равенство
$$
\int_{\R^d} \phi(2^\kappa x -\nu) \cdot \overline {g(x)} dx =0,
$$
а
$$
\int_{\R^d} (E_\kappa^p f_0)(x) \cdot \overline {g(x)} dx \ne 0,
$$
или
\begin{equation*} \tag{2.3.17}
\int_{\R^d} g(x) \cdot \overline {\phi(2^\kappa x -\nu)} dx =0, \nu \in \Z^d,
\end{equation*}
а
\begin{equation*} \tag{2.3.18}
\int_{\R^d} g(x) \cdot \overline {(E_\kappa^p f_0)(x) } dx \ne 0.
\end{equation*}
Из выполнения (2.3.17) в силу (2.3.7) следует равенство $ E_\kappa^{p^\prime} g =0, $
а, благодаря (2.3.4), имеем
$$
0 = \int_{\R^d} (e_\kappa^{p^\prime} g)(x) \cdot \overline{f_0(x)} dx =
\int_{\R^d} g(x) \cdot \overline{(e_{\kappa^p f_0)(x)}} dx,
$$
что противоречит (2.3.18). Таким образом, сделанное выше предположение
не верно, и, следовательно, имеет место (2.3.8). $ \square $

Предложение 2.3.4

Пусть $ d \in \N $ и наборы функций
$$
\mathcal \phi = \{\phi_1, \ldots, \phi_d\},
\tilde{\mathcal \phi} = \{\tilde \phi_1, \ldots, \tilde \phi_d \}
$$
таковы, что при $ j = 1, \ldots, d $ соблюдаются условия теоремы 2.2.3
с функциями $ \phi_j, \tilde \phi_j, $ а функции $ \phi, \tilde \phi $...
задаются равенствами
$$
\phi(x) = \prod_{j =1}^d \phi_j(x_j),
\tilde \phi(x) = \prod_{j =1}^d \tilde \phi_j(x_j), x \in \R^d.
$$
Тогда при $ 1 < p < \infty $ для любой функции $ f \in L_p(\R^d) $ в
пространстве $ L_p(\R^d) $ имеет место равенство
\begin{equation*} \tag{2.3.19}
f = \sum_{\kappa \in \Z_+^d} \mathcal E_\kappa^p f.
\end{equation*}

Доказательство.

Для доказательства (2.3.19) в силу предложения 1.1.1 достаточно показать,
что в условиях теоремы при $ 1 < p < \infty $ для $ f \in L_p(\R^d) $
справедливо соотношение
\begin{equation*} \tag{2.3.20}
\| f -E_\kappa^p f \|_{L_p(\R^d)} \to 0 \text{при } \mn(\kappa) \to \infty.
\end{equation*}

Для вывода соотношения (2.3.20) для $ f \in L_p(\R^d) $ сначала покажем, что оно
имеет место для любой функции $ f \in L_p(\R^d) $ вида
\begin{equation*} \tag{2.3.21}
f(x_1,\ldots,x_d) = \prod_{j =1}^d f_j(x_j) \text{почти для всех}
(x_1,\ldots,x_d) \in \R^d, f_j \in L_p(\R), j =1,\ldots,d.
\end{equation*}
В самом деле, для функции $ f \in L_p(\R^d) $ вида (2.3.21) в силу
(2.3.1), (1.3.2), (2.1.5), теоремы Фубини, (1.3.1) и (2.2.16) имеем
\begin{multline*} \tag{2.3.22}
\| f -E_\kappa^p f \|_{L_p(\R^d)} = \biggl\| \sum_{m =0}^{d -1} \biggl((\prod_{j =1}^m
V_j(E_{\kappa_j}^{j,p})) f -(\prod_{j =1}^{m +1}
V_j(E_{\kappa_j}^{j,p})) f\biggr) \biggr\|_{L_p(\R^d)} \\ \le
\sum_{m =0}^{d -1} \biggl\| \biggl(\prod_{j =1}^m
V_j(E_{\kappa_j}^{j,p})\biggr) f -\biggl(\prod_{j =1}^m
V_j(E_{\kappa_j}^{j,p})\biggr) ((V_{m +1}(E_{\kappa_{m +1}}^{m +1,p})) f) \biggr\|_{L_p(\R^d)} \\ =
\sum_{m =0}^{d -1} \biggl\| \biggl(\prod_{j =1}^m
V_j(E_{\kappa_j}^{j,p})\biggr) (f -(V_{m +1}(E_{\kappa_{m +1}}^{m +1,p})) f) \biggr\|_{L_p(\R^d)} \\ \le
\sum_{m =0}^{d -1} \biggl\| \prod_{j =1}^m
V_j(E_{\kappa_j}^{j,p}) \biggr\|_{\mathcal B(L_p(\R^d), L_p(\R^d)} \cdot
\| f -(V_{m +1}(E_{\kappa_{m +1}}^{m +1,p})) f \|_{L_p(\R^d)} \\ \le
\sum_{m =0}^{d -1} \biggl(\prod_{j =1}^m
\| E_{\kappa_j}^{j,p} \|_{\mathcal B(L_p(\R), L_p(\R))}\biggr) \cdot
\| f -(V_{m +1}(E_{\kappa_{m +1}}^{m +1,p})) f \|_{L_p(\R^d)} \\ \le
\sum_{m =0}^{d -1} \biggl(\prod_{j =1}^m c_1^j\biggr) \cdot
\| f -(V_{m +1}(E_{\kappa_{m +1}}^{m +1,p})) f \|_{L_p(\R^d)} \\ \le
c_2 \cdot \sum_{m =0}^{d -1}
\| f -(V_{m +1}(E_{\kappa_{m +1}}^{m +1,p})) f \|_{L_p(\R^d)}  =
c_2 \cdot \sum_{j =1}^{d}
\biggl(\int_{\R^d} | f(x) -((V_{j}(E_{\kappa_{j}}^{j,p})) f)(x) |^p dx\biggr)^{1/p} \\ =
c_2 \cdot \sum_{j =1}^{d}
\biggl(\int_{\R^{d -1}} \int_\R \biggl| f(x_1,\ldots,x_{j -1}, x_j,x_{j +1},
\ldots,x_d) \\
-((V_{j}(E_{\kappa_{j}}^{j,p})) f)(x_1, \ldots, x_{j -1},
x_j, x_{j +1}, \ldots, x_d) \biggr|^p dx_j dx_1 \ldots dx_{j -1} dx_{j +1} \ldots dx_d\biggr)^{1/p} \\ =
c_2 \cdot \sum_{j =1}^{d}
\biggl(\int_{\R^{d -1}} \int_\R \biggl| f_j(x_j) \cdot \biggl(\prod_{i =1, \ldots,d: i \ne j} f_i(x_i)\biggr)\\
-\biggl(E_{\kappa_{j}}^{j,p} \biggl(f_j(\cdot) \biggl(\prod_{i =1, \ldots,d: i \ne j}
f_i(x_i)\biggl)\biggr)\biggr)(x_j)\biggr|^p dx_j dx_1 \ldots dx_{j -1} dx_{j +1} \ldots dx_d\biggr)^{1/p} \\ =
c_2 \cdot \sum_{j =1}^{d}
\biggl(\int_{\R^{d -1}} \biggl| \prod_{i =1, \ldots,d: i \ne j} f_i(x_i)\biggr|^p \\
\times
\biggl(\int_\R | f_j(x_j) -(E_{\kappa_{j}}^{j,p} f_j)(x_j)|^p dx_j\biggr) dx_1 \ldots dx_{j -1} dx_{j +1} \ldots dx_d\biggr)^{1/p} \\ =
c_2 \cdot \sum_{j =1}^{d}
\biggl(\int_{\R^{d -1}} \biggl| \prod_{i =1, \ldots,d: i \ne j} f_i(x_i)\biggr|^p \cdot
\| f_j -(E_{\kappa_{j}}^{j,p} f_j) \|_{L_p(\R)}^p dx_1 \ldots dx_{j -1}
dx_{j +1} \ldots dx_d\biggr)^{1/p} \\ =
c_2 \cdot \sum_{j =1}^{d}
\biggl\| \prod_{i =1, \ldots,d: i \ne j} f_i(x_i) \biggr\|_{L_p(\R^{d -1})} \cdot
\| f_j -(E_{\kappa_{j}}^{j,p} f_j) \|_{L_p(\R)} \to 0 \text{ при }
\mn(\kappa) \to \infty.
\end{multline*}
Из (2.3.22) следует, что соотношение (2.3.20) справедливо для любой функции
$ f \in L_p(\R^d), $ принадлежащей линейной оболочке функций вида
(2.3.21). Отсюда, принимая во внимание, что линейная оболочка функций вида
(2.3.21) плотна в $ L_p(\R^d) $ (ибо она содержит функции вида
$ \sum_{i =1}^n a_i \chi_{Q_i}(x), a_i \in \C, Q_i = x_i^0 +\delta_i I^d,
x_i^0 \in \R^d, \delta_i \in \R_+^d, i =1, \ldots, n, n \in \N $),
и учитывая, что вследствие (2.3.1), (1.3.2), (2.1.5) нормы операторов
$ E_\kappa^p: L_p(\R^d) \mapsto L_p(\R^d), \kappa \in \Z_+^d, $ не 
превосходят общей константы, заключаем, что (2.3.20) соблюдается для любой 
функции $ f \in L_p(\R^d). \square $
\bigskip

2.4. В этом пункте устанавливается конечный результат работы теорема 2.4.2.
Но прежде, опираясь на лемму 2.1.1 и последующие утверждения, покажем, что 
имеет место

Теорема 2.4.1

Пусть $ d \in \N, 1 < p < \infty $ и соблюдаются условия предложения 2.3.4.
Тогда существует константа $ c_{1}(d, \mathcal \phi, \tilde{\mathcal \phi},p) >0 $ такая, что для
любого семейства чисел $ \{ \sigma_\kappa: \kappa \in \Z_+^d \} $ вида
$ \sigma_\kappa = \prod_{j=1}^d
\sigma^j_{\kappa_j}, $
где $ \sigma^j_{\kappa_j} \in \{-1, 1\}, \kappa_j \in \Z_+, j =1,\ldots,d, $
для $ f \in L_p(\R^d) $ справедливо неравенство
\begin{equation*} \tag{2.4.1}
\biggl\| \sum_{\kappa \in \Z_+^d} \sigma_\kappa \cdot
(\mathcal E_\kappa f ) \biggr\|_{L_p(\R^d)}
\le c_{1} \| f \|_{L_p(\R^d)},
\end{equation*}
где $ \mathcal E_\kappa = \mathcal E_\kappa^p, \kappa \in \Z_+^d $
определяется в предложении 2.3.1.

Доказательство.

Сначала покажем, что в условиях теоремы для любого непустого
множества $ J \subset \{1, \ldots, d \} $ при любом $ k^J \in (\Z_+^d)^J $
и любых наборах чисел $ \{ \sigma^j_{\kappa_j} \in \{-1, 1\}, \kappa_j
= 0,\ldots, k_j \}, j \in J, $
для $ f \in L_p(\R^d) $ имеет место неравенство
\begin{multline*} \tag{2.4.2}
\biggl\| \sum_{\kappa^J \in \Z_+^m(k^J)}  \biggl(\prod_{j \in J} (\sigma^j_{\kappa_j}
V_j (\mathcal E_{\kappa_j}^{j, p}))\biggr)f \biggr\|_{L_p(\R^d)}\\
\le \biggl(\prod_{j \in J} c_{1}(1,\phi_j, \tilde \phi_j,p)\biggr) \cdot
\| f \|_{L_p(\R^d)},
\end{multline*}
где $ m = \card J, $ а $ \mathcal E_{\kappa_j}^{j, p} =
\mathcal E_{\kappa_j}^{\phi_j, \tilde \phi_j,p}, \kappa_j \in \Z_+,
j =1,\ldots,d. $

Доказательство (2.4.2) проведем по индукции относительно $ m. $
При $ m =1, $ т.е. для $ j =1, \ldots, d, $ используя п. 2) леммы
1.3.1, теорему Фубини, (1.3.1), (2.1.4), имеем
\begin{multline*}
\biggl\| \sum_{\kappa_j =0}^{k_j} \sigma^j_{\kappa_j} \cdot (V_j
(\mathcal E_{\kappa_j}^{j, p})) f \biggr\|_{L_p(\R^d)}^p =
\biggl\| \biggl(V_j\biggl(\sum_{\kappa_j =0}^{k_j} \sigma^j_{\kappa_j} \mathcal
E_{\kappa_j}^{j, p}\biggr)\biggr) f \|_{L_p(\R^d)}^p \\
= \int_{\R^d} \biggl| \biggl(V_j (\sum_{\kappa_j =0}^{k_j} \sigma^j_{\kappa_j}
\mathcal E_{\kappa_j}^{j, p})\biggr) f \biggr|^p dx\\ =
\int_{\R^{d -1}}
\int_\R \biggl| \biggl(\biggl(V_j (\sum_{\kappa_j =0}^{k_j} \sigma^j_{\kappa_j}
\mathcal E_{\kappa_j}^{j, p})\biggr) f\biggr) (x_1,\ldots, x_{j -1}, x_j,
x_{j +1}, \ldots, x_d)\biggr|^p\\
\times dx_j dx_1 \ldots dx_{j -1} dx_{j +1} \ldots dx_d \\
= \int_{\R^{d -1}} \int_\R \biggl| \biggl( \biggl(\sum_{\kappa_j =0}^{k_j}
\sigma^j_{\kappa_j} \mathcal E_{\kappa_j}^{j, p}\biggr) f (x_1,\ldots,
x_{j -1}, \cdot, x_{j +1}, \ldots, x_d)\biggr)(x_j) \biggr|^p\\
\times dx_j dx_1 \ldots
dx_{j -1} dx_{j +1} \ldots dx_d\\
 = \int_{\R^{d -1}} \biggl\|
\biggl(\sum_{\kappa_j =0}^{k_j} \sigma^j_{\kappa_j} \mathcal
E_{\kappa_j}^{j, p}\biggr) f (x_1,\ldots, x_{j -1}, \cdot, x_{j +1},
\ldots, x_d) \biggr\|_{L_p(\R)}^p \\
\times dx_1 \ldots dx_{j -1} dx_{j +1} \ldots
dx_d \\
\le \int_{\R^{d -1}} \biggl( c_{1}(\phi_j, \tilde \phi_j,p) \| f (x_1,\ldots, x_{j -1},
\cdot, x_{j +1}, \ldots, x_d) \|_{L_p(\R)}\biggr)^p\\
\times dx_1 \ldots dx_{j -1} dx_{j +1} \ldots dx_d\\ =
(c_{1}(\phi_j, \tilde \phi_j, p))^p \int_{\R^{d -1}} \int_\R | f (x_1,\ldots, x_{j -1},
x_j, x_{j +1}, \ldots, x_d) |^p\\
\times dx_j
dx_1 \ldots dx_{j -1} dx_{j +1} \ldots dx_d \\
= ( c_{1}(\phi_j, \tilde \phi_j,p))^p
\int_{\R^d} | f(x)|^p dx = \bigl( c_{1}(\phi_j, \tilde \phi_j,p) \| f \|_{L_p(\R^d)}\bigr)^p,
\end{multline*}
откуда
\begin{equation*} \tag{2.4.3}
\biggl\| \sum_{\kappa_j =0}^{k_j} \sigma^j_{\kappa_j} (V_j (\mathcal
E_{\kappa_j}^{j, p})) f \biggr\|_{L_p(\R^d)} \le c_{1}(\phi_j, \tilde \phi_j, p)
\| f \|_{L_p(\R^d)},
\end{equation*}
что совпадает с (2.4.2) при $ m =1, J = \{j\}. $

Предположим теперь, что при некотором $ m: 1 \le m \le d -1, $
оценка (2.4.2) имеет место для любого множества $ J \subset \{1,
\ldots, d \}: \card J = m, $ при любом $ k^J \in (\Z_+^d)^J, $
любых наборах чисел $ \{ \sigma^j_{\kappa_j} \in \{-1, 1\},
\kappa_j = 0, \ldots, k_j \}, j \in J, $ и любой функции $ f \in
L_p(\R^d).$ Покажем, что тогда неравенство (2.4.2) справедливо при
$ m+1 $ вместо $ m $ для любого множества $ \mathcal J \subset
\{1, \ldots, d \} $ вместо $ J, $ у которого $ \card \mathcal J =
m +1, $ при любом $ k^{ \mathcal J} \in (\Z_+^d)^{ \mathcal J}, $
любых наборах чисел $ \{ \sigma^j_{\kappa_j} \in \{-1, 1\},
\kappa_j = 0, \ldots, k_j \}, j \in \mathcal J, $ и любой функции
$ f \in L_p(\R^d). $ Представляя множество $ \mathcal J \subset
\{1, \ldots, d \}: \card \mathcal J = m +1, $ в виде $ \mathcal J
= J \cup \{i\}, i \notin J, $ с учетом (1.3.3) в силу (2.4.3) и
предположения индукции получаем
\begin{multline*}
\biggl\| \sum_{\kappa^{\mathcal J} \in \Z_+^{m+1}(k^{\mathcal J})}
\biggl(\prod_{j \in \mathcal J} (\sigma^j_{\kappa_j} V_j (\mathcal
E_{\kappa_j}^{j, p}))\biggr)f \biggr\|_{L_p(\R^d)} \\
= \biggl\| \sum_{(\kappa_i, \kappa^J): \kappa_i =0, \ldots, k_i,
\kappa^J \in \Z_+^m(k^J)} \sigma^i_{\kappa_i} (V_i (\mathcal
E_{\kappa_i}^{i, p})) \biggl(\biggl(\prod_{j \in J} (\sigma^j_{\kappa_j} V_j
(\mathcal E_{\kappa_j}^{j, p})) \biggr)f\biggr) \biggr\|_{L_p(\R^d)}\\ =
\biggl\| \sum_{\kappa_i =0}^{k_i} \sigma^i_{\kappa_i} (V_i (\mathcal
E_{\kappa_i}^{i, p})) \biggl(\sum_{\kappa^J \in \Z_+^m(k^J)}  \biggl(\prod_{j
\in J} (\sigma^j_{\kappa_j} V_j (\mathcal E_{\kappa_j}^{j, p}))
\biggr)f\biggr) \biggr\|_{L_p(\R^d)} \\
\le c_{1}(\phi_i, \phi_i^*, p) \biggl\| \sum_{\kappa^J \in \Z_+^m(k^J)}  \biggl(\prod_{j
\in J} (\sigma^j_{\kappa_j} V_j (\mathcal E_{\kappa_j}^{j, p})
)\biggr)f \biggr\|_{L_p(\R^d)} \\
\le c_{1}(\phi_i, \phi_i^*,p) \biggl(\prod_{j \in J}
c_{1}(\phi_j, \tilde \phi_j, p)\biggr) \cdot \| f \|_{L_p(\R^d)} = \biggl(\prod_{j \in \mathcal
J} c_{1}(\phi_j, \tilde \phi_j, p)\biggr) \cdot \| f \|_{L_p(\R^d)},
\end{multline*}
что завершает вывод (2.4.2).

В частности, из (2.4.2) при $ m = d $ ввиду (2.3.2) получаем, что
в условиях теоремы при любом $ k \in \Z_+^d $ соблюдается неравенство
\begin{equation*} \tag{2.4.4}
\biggl\| \sum_{\kappa \in \Z_+^d(k)} \sigma_\kappa \cdot (\mathcal
E_\kappa f ) \biggr\|_{L_p(\R^d)} \le c_{1} \| f \|_{L_p(\R^d)},
\sigma_\kappa = \prod_{j=1}^d \sigma^j_{\kappa_j}
\end{equation*}
 где $\sigma^j_{\kappa_j} \in \{-1, 1\}, \kappa_j
\in \Z_+, j =1,\ldots,d,  f \in L_p(\R^d),
c_{1} = \prod_{j=1}^d c_{1}(\phi_j, \tilde \phi_j, p)$.

Теперь убедимся в справедливости (2.4.1).
Для произвольного семейства чисел $ \{ \sigma_\kappa = \prod_{j=1}^d
\sigma^j_{\kappa_j}: \sigma^j_{\kappa_j} \in \{-1, 1\}, j =1,\ldots,d,
\kappa \in \Z_+^d \}, $ функции $ f \in L_p(\R^d) $ рассмотрим последовательность
$$
\biggl\{ \biggl(\sum_{\kappa \in \Z_+^d(k \e)} \sigma_\kappa \cdot
(\mathcal E_\kappa f )\biggr) \in L_p(\R^d), k \in \Z_+ \biggr\}
$$
и, принимая во внимание (2.4.4), секвенциальную компактность
шара $ B(L_p(\R^d)) $ относительно $*$-слабой топологии в пространсте
$ L_p(\R^d) = (L_{p^\prime}(\R^d))^*, $ выберем
подпоследовательность
$$
\biggl\{ \sum_{\kappa \in \Z_+^d(k_n \e)} \sigma_\kappa \cdot
(\mathcal E_\kappa f ): k_n < k_{n+1}, n \in \N \biggr\}
$$
и функцию $ F \in L_p(\R^d), $ обладающие тем свойством, что для
любой функции $ g \in L_{p^\prime}(\R^d) $ выполняется равенство
\begin{equation*} \tag{2.4.5}
\lim_{ n \to \infty} \int_{\R^d} \biggl(\sum_{\kappa \in \Z_+^d(k_n \e)}
\sigma_\kappa (\mathcal E_\kappa f )\biggr) g dx =
\int_{\R^d} F g dx.
\end{equation*}
Заметим, что при любом $ \kappa \in \Z_+^d, $ благодаря (2.3.5),
(2.4.5), (2.3.6), для  $ g \in L_{p^\prime}(\R^d) $ имеет место равенство
\begin{multline*}
\int_{\R^d} (\mathcal E_\kappa F) \cdot \overline g dx = \int_{\R^d} F
\cdot \overline {(\mathcal E_\kappa g)} dx = \lim_{ n \to \infty}
\int_{\R^d} \biggl(\sum_{\kappa^\prime \in \Z_+^d(k_n \e)}
\sigma_{\kappa^\prime} \cdot (\mathcal E_{\kappa^\prime} f)\biggr)
\overline{(\mathcal E_\kappa g)} dx \\
= \lim_{ n \to \infty} \int_{\R^d} \mathcal E_\kappa
\biggl(\sum_{\kappa^\prime \in \Z_+^d(k_n \e)}
\sigma_{\kappa^\prime} \cdot (\mathcal E_{\kappa^\prime} f)\biggr) \cdot
\overline g dx \\
= \lim_{ n \to \infty} \int_{\R^d}
\biggl(\sum_{\kappa^\prime \in \Z_+^d(k_n \e)} \sigma_{\kappa^\prime}
\cdot \mathcal E_\kappa^{p} (\mathcal E_{\kappa^\prime}^{p} f)\biggr)
\overline g dx  = \int_{\R^d} \sigma_{\kappa} \cdot (\mathcal
E_\kappa f) \overline g dx,
\end{multline*}
и, значит,
\begin{equation*} \tag{2.4.6}
\mathcal E_\kappa F =
\sigma_{\kappa} \cdot (\mathcal E_\kappa f).
\end{equation*}
Учитывая (2.4.6), (2.3.19), заключаем, что
\begin{equation*}
\sum_{\kappa \in \Z_+^d(k)} \sigma_\kappa \cdot (\mathcal
E_\kappa f ) = \sum_{\kappa \in \Z_+^d(k)} (\mathcal
E_\kappa F)
\end{equation*}
сходится к $ F $ в $ L_p(\R^d) $ при $ \mn(k) \to \infty. $
Поэтому, переходя к пределу при $ \mn(k) \to \infty $ в
неравенстве (2.4.4), приходим к (2.4.1). $ \square $

Следствие

В условиях теоремы 2.4.1 для любого семейства чисел $ \{ \sigma_\kappa:
\kappa \in \Z_+^d \} $ вида $ \sigma_\kappa = \prod_{j=1}^d
\sigma^j_{\kappa_j}, $ где $ \sigma^j_{\kappa_j} \in \{-1, 1\}, \kappa_j
\in \Z_+, j =1,\ldots,d, $ и любой функции $ f \in L_p(\R^d) $
соблюдается неравенство
\begin{equation*} \tag{2.4.7}
\| f \|_{L_p(\R^d)} \le c_{1}
\biggl\| \sum_{\kappa \in \Z_+^d} \sigma_\kappa \cdot
(\mathcal E_\kappa f ) \biggr\|_{L_p(\R^d)}.
\end{equation*}

Доказательство.

Сначала покажем, что при любом $ k \in \Z_+^d, $ любом наборе чисел
$ \{ \sigma_\kappa = \prod_{j=1}^d \sigma^j_{\kappa_j}: \sigma^j_{\kappa_j} \in \{-1, 1\}, j =1,\ldots,d,
\kappa \in \Z_+^d(k) \}, $ для $ f \in L_p(\R^d) $
справедливо неравенство
\begin{equation*} \tag{2.4.8}
\| E_k f \|_{L_p(\R^d)} \le c_{1}
\biggl\| \sum_{\kappa \in \Z_+^d(k)}
\sigma_\kappa \cdot (\mathcal E_\kappa f ) \biggr\|_{L_p(\R^d)}.
\end{equation*}

В самом деле, ввиду (2.3.6), (1.1.1) имеем
\begin{multline*}
\biggl(\sum_{\kappa \in \Z_+^d(k)} \sigma_\kappa \cdot \mathcal
E_\kappa^{p}\biggr)^2 f = \biggl(\sum_{\kappa, \kappa^\prime \in \Z_+^d(k)}
\sigma_\kappa \sigma_{\kappa^\prime} \cdot \mathcal E_\kappa^{p}
\mathcal E_{\kappa^\prime}^p\biggr) f \\
= \biggl(\sum_{\kappa \in \Z_+^d(k)} \sigma_\kappa^2 \cdot \mathcal E_\kappa^{p}\biggr) f =
\sum_{\kappa \in \Z_+^d(k)} \mathcal E_\kappa f = E_k f.
\end{multline*}
Откуда, применяя (2.4.4), выводим
\begin{multline*}
\| E_k f \|_{L_p(\R^d)} = \biggl\| \biggl(\sum_{\kappa \in \Z_+^d(k)}
\sigma_\kappa \cdot \mathcal E_\kappa\biggr)^2 f \biggr\|_{L_p(\R^d)}\\
 =\biggl\| \sum_{\kappa \in \Z_+^d(k)} \sigma_\kappa \cdot \mathcal
E_\kappa \biggl(\sum_{\kappa^\prime \in \Z_+^d(k)}
\sigma_{\kappa^\prime} \cdot (\mathcal E_{\kappa^\prime} f)\biggr)
\biggr\|_{L_p(\R^d)}
\le c_{1} \| \sum_{\kappa^\prime \in \Z_+^d(k)}
\sigma_{\kappa^\prime} \cdot (\mathcal E_{\kappa^\prime} f)
\|_{L_p(\R^d)}.
\end{multline*}

Как видно из вывода (2.4.1) и (2.3.20), в неравенстве (2.4.8) можно
перейти к пределу при $ \mn(k) \to \infty, $ в результате чего
получим (2.4.7). $ \square $

С помощью теоремы 2.4.1 и следствия из нее, опираясь на схему доказательства
теоремы Литтлвуда-Пэли, изложенную в [3] для
операторов взятия частных сумм кратных рядов Фурье, устанавливается

Теорема 2.4.2

При соблюдении условий теоремы 2.4.1 существуют константы $ c_{2}(d, \mathcal \phi, \tilde{\mathcal \phi}, p) >0,\\
c_{3}(d,\mathcal \phi, \tilde{\mathcal \phi}, p) >0 $ такие, что для любой
функции $ f \in L_p(\R^d) $
выполняются неравенства
\begin{equation*} \tag{2.4.9}
c_{2} \| f \|_{L_p(\R^d)} \le
\biggl( \int_{\R^d} \biggl(\sum_{\kappa \in \Z_+^d} |(\mathcal E_\kappa f)(x)|^2 \biggr)^{p/2}
dx\biggr)^{1/p} \le c_{3} \| f \|_{L_p(\R^d)}.
\end{equation*}

Доказательство.

Рассмотрим систему Радемахера, состоящую из функций
$$
\omega_\kappa(t) = \operatorname{sign} \sin( 2^{\kappa +1} \pi t),
t \in I, \kappa \in \Z_+,
$$
и определим семейство функций $ \omega_\kappa^d, \kappa \in \Z_+^d, $
полагая
$$
\omega_\kappa^d(t) = \prod_{j=1}^d \omega_{\kappa_j}(t_j), t \in I^d.
$$
Как известно (см., например, [3, п. 1.5.2]), существуют константы $ c_{4}(d,p) >0,
c_{5}(d,p) >0 $ такие, что при любом $ k \in \Z_+^d $ для любого набора
чисел $ \{ a_\kappa \in \C, \kappa \in \Z_+^d(k) \} $
имеет место неравенство
\begin{equation*} \tag{2.4.10}
c_{4} ( \sum_{ \kappa \in \Z_+^d(k)} | a_\kappa|^2 )^{1/2}
\le
\| \sum_{ \kappa \in \Z_+^d(k)} a_\kappa \omega_\kappa^d(\cdot) \|_{L_p(I^d)}
\le
c_{5} ( \sum_{ \kappa \in \Z_+^d(k)} | a_\kappa|^2 )^{1/2}.
\end{equation*}

Для $ f \in L_p(\R^d) $ при $ k \in \Z_+^d, $ используя (2.4.8), теорему
Фубини, (2.4.10), выводим
\begin{multline*} \tag{2.4.11}
\| E_k f \|_{L_p(\R^d)}^p = \int_{I^d} \| E_k f \|_{L_p(\R^d)}^p dt \le
\int_{I^d} \biggl(c_{1} \biggl\| \sum_{\kappa \in
\Z_+^d(k)} \omega_\kappa^d(t) \cdot (\mathcal E_\kappa f )
\biggr\|_{L_p(\R^d)}\biggr)^p dt\\
 = (c_{1})^p \int_{I^d} \int_{\R^d} \biggl|
\sum_{\kappa \in \Z_+^d(k)} \omega_\kappa^d(t) \cdot (\mathcal
E_\kappa f )(x) \biggr|^p dx dt \\
= (c_{1})^p \int_{\R^d}
\int_{I^d} \biggl| \sum_{\kappa \in \Z_+^d(k)} \omega_\kappa^d(t) \cdot
(\mathcal
E_\kappa f )(x) \biggr|^p dt dx\\
 = (c_{1})^p \int_{\R^d} \biggl\|
\sum_{\kappa \in \Z_+^d(k)} (\mathcal E_\kappa f )(x) \cdot
\omega_\kappa^d(\cdot) \biggr\|_{L_p(I^d)}^p dx \\
\le (c_{1})^p \int_{\R^d} \biggl(c_{5} \biggl(\sum_{\kappa \in \Z_+^d(k)}
|(\mathcal E_\kappa f )(x)|^2 \biggr)^{1/2}\biggr)^p dx\\
 = (c_{6})^p \int_{\R^d}
\biggl(\sum_{\kappa \in \Z_+^d(k)} |(\mathcal E_\kappa f )(x)|^2
\biggr)^{p/2} dx,
\end{multline*}
и, пользуясь (2.4.10), теоремой Фубини, (2.4.4), получаем
\begin{multline*}
\int_{\R^d} \biggl(\sum_{\kappa \in \Z_+^d(k)} |(\mathcal E_\kappa f )(x)|^2 \biggr)^{p/2} dx
\le \int_{\R^d} \biggl(c_{7} \biggl\| \sum_{\kappa \in \Z_+^d(k)}
(\mathcal E_\kappa f )(x) \cdot \omega_\kappa^d(\cdot) \biggr\|_{L_p(I^d)}\biggr)^p dx\\
= (c_{7})^p \int_{\R^d} \int_{I^d} \biggl| \sum_{\kappa \in \Z_+^d(k)}
\omega_\kappa^d(t) \cdot (\mathcal E_\kappa f )(x) \biggr|^p dt dx\\
= (c_{7})^p \int_{I^d} \int_{\R^d} \biggl| \sum_{\kappa \in
\Z_+^d(k)} \omega_\kappa^d(t) \cdot (\mathcal E_\kappa f)(x) \biggr|^p dx dt\\
= (c_{7})^p \int_{I^d} \biggl\| \sum_{\kappa \in \Z_+^d(k)}
\omega_\kappa^d(t) \cdot (\mathcal E_\kappa f )
\biggr\|_{L_p(\R^d)}^p dt \\
\le (c_{7})^p \int_{i^d} \bigl(c_{1} \| f \|_{L_p(\R^d)}\bigr)^p dt =
(c_{3})^p \| f \|_{L_p(\R^d)}^p,
\end{multline*}
откуда, в частности, имеем
\begin{equation*} \tag{2.4.12}
\int_{\R^d} \biggl(\sum_{\kappa \in \Z_+^d(n \e)}
|(\mathcal E_\kappa f )(x)|^2 \biggr)^{p/2} dx \le (c_{3})^p
\| f\|_{L_p(\R^d)}^p, n \in \Z_+.
\end{equation*}

Для получения второго неравенства в (2.4.9) достаточно применить
теорему Леви о предельном переходе под знаком интеграла к монотонно
возрастающей последовательности функций
$ \{ (\sum_{\kappa \Z_+^d(n \e)} |(\mathcal E_\kappa f)(x)|^2 )^{p/2},
n \in \Z_+\}, $ принимая во внимание (2.4.12) и учитывая, что
почти для всех $ x \in \R^d $ предел
\begin{multline*}
\lim_{n \to \infty} \biggl(\sum_{\kappa \in \Z_+^d(n \e)} |(\mathcal
E_\kappa f)(x)|^2\biggr )^{p/2} =
\biggl(\lim_{n \to \infty} \sum_{\kappa \in \Z_+^d(n \e)}
|(\mathcal E_\kappa f)(x)|^2\biggr )^{p/2} \\ =
\biggl(\sum_{\kappa
\in \Z_+^d} |(\mathcal E_\kappa f)(x)|^2 \biggr)^{p/2}.
\end{multline*}
Последнее равенство справедливо в силу того, что при $ k \in
\Z_+^d $ имеет место соотношение
\begin{multline*}
\sum_{\kappa \in \Z_+^d(\mn(k) \e)} |(\mathcal E_\kappa f)(x)|^2 \le
\sum_{\kappa \in \Z_+^d(k)} |(\mathcal E_\kappa f)(x)|^2 \le
\sum_{\kappa \in \Z_+^d(\mx(k) \e)} |(\mathcal E_\kappa f)(x)|^2, x \in \R^d.
\end{multline*}

Ввиду сказанного и на основании (2.3.20) переходя к пределу в (2.4.11) при
$ \mn(k) \to \infty, $ приходим
к первому неравенству в (2.4.9). $ \square $
\newpage

\centerline{Список литературы}
[1] Темляков В. Н. Приближение функций с ограниченной смешанной производной //
Тр. МИАН СССР, 178 (1986), 3--113.

[2]
Галеев Э. М.
Поперечники классов Бесова $ B_{p,\theta}^r(T^d) $ //
Матем. заметки, 69:5 (2001), 656--665.

[3]
Никольский С. М.
Приближение функций многих переменных и теоремы вложения // М.: Наука. 1977.

[4]
Бесов О. В.
Теорема Литтлвуда-Пэли для смешанной нормы //
Тр. МИАН СССР, 170 (1984), 31--36.

[5]
Кудрявцев С. Н.
Теорема типа Литтлвуда–Пэли и следствие из нее //
Изв. РАН. Сер. матем., 77:6 (2013),  97--138.

[6]
Чуи Ч. К.
Введение в вейвлеты // М.: Мир, 2001.

[7]
Новиков И. Я., Протасов В. Ю., Скопина М.  А.
Теория всплесков // М.: ФИЗМАТЛИТ, 2006.

[8]
Meyer Y.
Wavelets and operators // Cambridge University Press, 1992.

[9]
Hernandez E. and Weiss G. //
A first course of wavalets, CRC Press LLC, Boca Raton, FL. 1996.

[10]
Кудрявцев С. Н.
Обобщённые ряды Хаара и их применение //
Analysis Mathematica, 37:2 (2011), 103 -- 150.

[11]
Кудрявцев С. Н.
Теорема типа Литтлвуда-Пэли для ортопроекторов
на подпространства всплесков //
http://arxiv.org/abs/1204.1830

[12]
Кудрявцев С. Н.
Приближение производных функций конечной гладкости из неизотропных классов //
Изв. РАН. Сер. матем. 68:1 (2004), 79 -- 122.

[13]
Кудрявцев С. Н.
Приближение и восстановление производных для функций, удовлетворяющих
смешанным условиям Гельдера // Изв. РАН. Сер. матем. 71:5 (2007), 37 -- 80.

[14]
Стейн И. М.
Сингулярные интегралы и дифференциальные свойства функций // М.: Мир, 1973.

[15]
Wojtaszczyk P.
Wavelets as unconditional bases in $ L_p(\R) $ //
J. Fourier Anal. Appl. 5:1 (1999), 73 -- 85.

[16]
Кашин Б. С., Саакян А.А.
Ортогональные ряды // М.: Изд-во АФЦ, 1999.

\end{document}